\documentclass[english]{smfart}
\usepackage{hyperref,sabbah_gaussian-type}
\begin{document}
\frontmatter
\title{Differential systems of pure Gaussian type}

\author[C.~Sabbah]{Claude Sabbah}
\address{UMR 7640 du CNRS\\
Centre de Mathématiques Laurent Schwartz\\
École polytechnique\\
F--91128 Palaiseau cedex\\
France}
\email{Claude.Sabbah@polytechnique.edu}
\urladdr{http://www.math.polytechnique.fr/~sabbah}

\thanks{This research was supported by the grants ANR-08-BLAN-0317-01 and ANR-13-IS01-0001-01 of the Agence nationale de la recherche.}

\begin{abstract}
We give the transformation rule for the Stokes data of the Laplace transform of a differential system of pure Gaussian type.
\end{abstract}


\dedicatory{Dedicated to the memory of Andrey Bolibrukh\hspace*{5mm}\vspace*{\baselineskip}}
\maketitle
\bgroup
\let\\\relax\tableofcontents
\egroup
\mainmatter

\section*{Introduction}
Computing the behaviour of the Stokes structure of a differential equation in one complex variable by the Laplace transformation leads in general to difficult combinatorial problems. In this article we make explicit the topological Laplace transformation in a simple case, that of differential systems of pure Gaussian type. It is well-known that the function $t\mto\exp(ct^2/2)$ has a simple behaviour by the Fourier transformation of the real variable $t$. Differential systems of pure Gaussian type are those systems of the complex variable~$t$ whose solutions asymptotically behave like sums of terms $t^\alpha(\log t)^k\exp(ct^2/2)$ when $t\to\infty$, and which have no other singularities. Their Laplace transform satisfy the same property, and the question we address is the computation of the Stokes data at infinity of the Laplace transform of such a system, in terms of the Stokes data at infinity of the original one.

To begin with, we relate various ways to encode the Stokes phenomenon of such a system. The sheaf-theoretic way (filtered local systems in the sense of Deligne \cite{Deligne78}) is suitable for the computation involving higher dimensional underlying spaces (\ie the computation of the Laplace transform by an integral formula). However, it is more common to express the Stokes phenomenon by linear algebra objects, like Stokes matrices. Here we find convenient to express it by a family of pairwise opposite filtrations of a vector space.

The computation is not too difficult, and the result easy to formulate, in the case when all nonzero complex numbers $c$ occurring in the asymptotic expansions of the solutions have the same argument. This is explained in Section \ref{sec:computationLaplace}, which could be regarded as a supplementary exercise in \cite{Malgrange91} illustrating the computation of the topological Laplace transformation. However, in general, the way this subset $C$ is embedded in $\CC^*$ introduces complicated combinatorial problems. In such a case, we reduce the question to the determination of the behaviour of the Stokes data when the subset $C$ varies in $\CC^*$. However, we do not give a precise answer to this latter question.

Our aim is to develop a topological approach for the computation of Stokes data of the Laplace transform of a system of linear differential equations in one variable, following the general method of \cite{Bibi10}. Other techniques are classically developed, mainly analytic ones (\cf \eg \cite{Malgrange91}, in particular Chap.\,XII, and the references therein), but our method gives rationality results in a straightforward way (if the Stokes data of the original system can be defined over $\QQ$, then so are the Stokes data of the Laplace transformed system). A more complicated example is analyzed in \cite{H-S14}. Our approach relies on a fundamental theorem by T.\,Mochizuki \cite{Mochizuki10} (\cf Theorem \ref{th:mochi}), who has also developed a slightly different method for such computations \cite{Mochizuki09b}. Lastly, let us mention a completely different method, still of a topological nature, following from the general Riemann-Hilbert correspondence of d'Agnolo and Kashiwara \cite{D-K13}.

\subsubsection*{Acknowledgements}
I thank Marco Hien and Takuro Mochizuki for stimulating discussions on this subject.

\section{Differential systems of pure Gaussian type}
\subsection{Setting an notation}\label{subsec:settingnotation}
We consider the Riemann sphere $\PP^1$ as covered by two affine charts~$\Afu_t$ with coordinate $t$ and $\Afu_{t'}$ with coordinate $t'$ such that $t'=1/t$ on the intersection of the two charts. We denote by $\Clt$ the ring of differential operators with coefficients in $\CC[t]$. All $\Clt$-modules will be left modules, that we identify with $\CC[t]$-modules with connection. Holonomic $\Clt$-modules (\ie torsion modules over $\Clt$) can be extended as holonomic modules over the sheaf $\cD_{\PP^1}$ of algebraic differential operators on $\PP^1$, in such a way that in the chart $t'$, multiplication by $t'$ is invertible.

Given $\varphi\in\CC[t]$, a basic example of a $\Clt$-module is $E^{\varphi(t)}:=(\CC[t],\rd+\rd\varphi)$. We will also use this notation in two variables $t,\tau$, with $\varphi\in\CC[t,\tau]$. The extension of $E^{\varphi(t)}$ to $\PP^1$, when restricted to the chart $t'$, denoted then $E^{\varphi(1/t')}$, is nothing but $(\CC[t',t^{\prime-1}],\rd+\rd\varphi(1/t'))$. By extending scalars from $\CC[t',t^{\prime-1}]$ to $\CC\lpr t'\rpr:=\CC\lcr t'\rcr[t^{\prime-1}]$, we get $\cE^{\varphi(1/t'})=\CC\lpr t'\rpr\otimes_{\CC[t',t^{\prime-1}]}E^{\varphi(1/t')}=(\CC\lpr t'\rpr,\rd+\rd\varphi(1/t'))$.

Given any morphism $f$ between complex manifolds or smooth algebraic varieties, we denote by $f^+$ the pull-back functor of left $\cD$-modules. This is nothing but $f^*$ for the underlying $\cO$-modules, together with the pull-back connection.

\begin{definition}\label{def:puregaussian}
Let $C$ be a finite subset of $\CC^*$. A differential system of \emph{pure Gaussian type~$C$} is a free $\CC[t]$-module $M$ of finite rank $r$ equipped with a connection $\nabla=\rd+A(t)\rd t$, where
\begin{itemize}
\item
$A(t)$ is a $r\times r$ matrix with entries in $\CC[t]$,
\item
setting $M'=\CC[t,t^{-1}]\otimes_{\CC[t]}M$, then $(\CC\lpr t'\rpr\otimes_{\CC[t',t^{\prime-1}]} M',\nabla)$ is isomorphic to a direct sum $\bigoplus_{c\in C}(\cE^{-c/2t^{\prime2}}\otimes R_c)$, where $\cE^{-c/2t^{\prime2}}=(\CC\lpr t'\rpr,\rd+c/t^{\prime3})$ and $R_c$ is a finite dimensional $\CC\lpr t'\rpr$-vector space with regular singular connection.
\end{itemize}
We will say that $M$ has pure Gaussian type if there exists $C$ such that $M$ has pure Gaussian type $C$.
\end{definition}

Note that, by our assumption, a differential system of pure Gaussian type $C$ is purely irregular at infinity, of slope $2$ and irregularity $2r$. We can regard $M$ as a $\Clt$-module: using a $\CC[t]$-basis $\bmm$, the action of $\partial_t$ is given by $\partial_t\bmm=\bmm\cdot A(t)$.

\begin{proposition}
Any sub or quotient $\Clt$-module of a differential system of pure Gaussian type $C$ has also pure Gaussian type $C$. The full subcategory of the category of holonomic $\Clt$-modules consisting of objects of pure Gaussian type $C$ is abelian.
\end{proposition}

\begin{proof}
Let $M$ be of pure Gaussian type $C$ and let $N\subset M$ be a $\Clt$-submodule. The characteristic variety of~$N$ and $M/N$ is reduced to the zero section on $\Afu_t$, hence both are $\CC[t]$-locally free of finite rank, hence free of finite rank. It remains to check their behaviour at infinity, and the statement reduces to proving that any $\CC\lpr t'\rpr$-subspace or quotient space with connection of a direct sum $\bigoplus_{c\in C}(\cE^{-c/2t^{\prime2}}\otimes R_c)$ takes the same form, a result which is easy by noticing that there is no nonzero morphism $(\cE^{-c/2t^{\prime2}}\otimes R_c)\to(\cE^{-c'/2t^{\prime2}}\otimes R_{c'})$ if $c\neq c'$. The last statement of the proposition is then clear.
\end{proof}

\begin{remarque}[Non rigidity]
A $\Clt$-module $M$ of pure Gaussian type is rigid (\ie its index of rigidity $\mathrm{rig}(M)$ is equal to $2$) if and only if $r=1$. Indeed, the index of rigidity is computed as (\cf\cite{Arinkin08})
\[
\mathrm{rig}(M)=2r^2-\mathrm{irr}_\infty(\End M) -r^2+\eta
\]
with $\eta=\dim\ker\partial_{t'}$ acting on the regular part of $\CC\lpr t'\rpr\otimes (\End M)$. We have $\eta\geq\#C$. On the other hand, $\mathrm{irr}_\infty(\End M)=2\sum_{c\neq c'}r_cr_{c'}$, with $r_c=\rk R_c$. Therefore,
\[
\mathrm{rig}(M)=\eta+\sum_{c\in C}r_c^2\geq\#C+\sum_{c\in C}r_c^2\geq2(\#C).
\]
\end{remarque}

\subsection{Behaviour by Laplace transformation}\label{subsec:opDmod}
We consider the behaviour of differential systems of pure Gaussian type by the Laplace transformation. We will consider the transformation with kernel $\exp(-t\tau)$, that is, we will set $\tau=\partial_t$ and $\partial_\tau=-t$, and we denote by $\wh M$ the $\CC$-vector space $M$ regarded as a $\Cltau$-module through this correspondence. The transformed differential system remains of pure Gaussian type, and the formal behaviour at infinity is made precise in the proof of the following lemma.

\begin{lemme}
Let $M$ be a $\Clt$-module of pure Gaussian type~$C$. Then its Laplace transform $\wh M$ has pure Gaussian type $\wh C:=-1/C=\{-1/c\mid c\in C\}$ and $\rk_{\CC[\tau]}\wh M=\rk_{\CC[t]}M=r$.
\end{lemme}

\begin{proof}
Firstly, the formal stationary phase formula implies that $\wh M$ has singularities at most at $\tau=\infty$ and $\tau=0$, the latter being regular. The decomposition of $\CC\lpr\tau'\rpr\otimes_{\CC[\tau,\tau^{-1}]}\wh M$ is obtained through \cite[(5.10)]{Bibi07a}, which shows (as well-known) that the local Laplace transform of $(\cE^{-c/2t^{\prime2}}\otimes R_c)$ is isomorphic to~$\cE^{1/2c\tau^{\prime2}}\otimes\nobreak R_c$ for $c\neq\nobreak0$. This implies that the formal decomposition of the irregular part of $\wh M$ at $\tau=\infty$ is $\bigoplus_{c\in C}(\cE^{1/2c\tau^{\prime2}}\otimes\nobreak R_c)$, where $R_c$ is now regarded as a $\CC\lpr\tau'\rpr$-vector space with regular singular connection by simply renaming the variable $t'$ to $\tau'$. The regular part of~$\wh M$ at $\tau=\infty$ has rank equal to the dimension of vanishing cycles of the analytic de~Rham complex $\DR^\an M$ at $t=0$. Our assumption implies that this rank is zero. We conclude that $\CC[\tau,\tau^{-1}]\otimes_{\CC[\tau]}\wh M$ is a free $\CC[\tau,\tau^{-1}]$-module of rank $r$.

It remains to show that $\wh M$ is non singular at $\tau=0$. Since the moderate nearby cycles of $M$ at $t=\infty$ are zero, the moderate vanishing cycles of $\wh M$ at $\tau=0$ also vanish, according to the standard correspondence shown in \cite[Prop.\,4.1(iv)]{Bibi05b}. Therefore, since the singularity of $\wh M$ at $\tau=0$ is regular, $\wh M$ has no singularity at $\tau=0$.
\end{proof}

\begin{remarque}\label{rem:Lpm}
The inverse Laplace transformation is the transformation defined by the correspondence $\tau=-\partial_t$, $\partial_\tau=t$. For $M$ of pure Gaussian type, both Laplace and inverse Laplace transformed objects have isomorphic formal models at infinity, but the Stokes structures may be non-isomorphic (\cf Remark \ref{rem:Fpm}).
\end{remarque}

\section{Stokes data of Gaussian type and\texorpdfstring{\\}{cr} the Riemann-Hilbert correspondence}\label{sec:Stokesfil}
In this section we recall the notion of Stokes filtration as defined in \cite{Deligne78} (\cf also \cite{Malgrange83bb}, \cite{B-V89}, \cite{Malgrange91}) in the particular case of Stokes filtrations which are of Gaussian type. We make explicit the correspondence with the more classical approach via Stokes data.

\subsection{Stokes filtration}\label{subsec:Stokesfil}
Let $\kk$ be a field (\eg $\QQ$ or $\CC$). Let~$\cL$ be a local system of finite dimensional $\kk$-vector spaces on the circle $S^1$ with coordinate $\re^{\ri\theta}$ (it should be noted here that, with respect to our original problem, we set $\theta=\arg t'=-\arg t$). We will usually set $r=\rk\cL$. A Stokes filtration of~$\cL$ of \emph{Gaussian type} is a family of subsheaves $\cL_{\leq c}\subset\cL$, with $c\in\CC$, satisfying the following properties:
\begin{enumerate}
\item\label{enum:Stokesfil1}
For each $\theta\in \RR/2\pi\ZZ$, let $\leqtheta$ be the partial order on $\CC$ which is compatible with addition and satisfies
\[
c\leqtheta0\iff c=0\text{ or } \arg c-2\theta\in(\pi/2,3\pi/2)\mod2\pi.
\]
(This means that $\exp(ct^2/2)=\exp(c/2t^{\prime2})$ has moderate growth in the neighbourhood of the point $(|t'|=0,\arg t'=\theta)$ in the real blow-up space $\wt\PP^1_t$ of $\PP^1_t$ at $t=\infty$. We will regard $S^1$ as the circle $S^1_\infty:=\wt\PP^1_t{}_{|\infty}$.) We also set $c\letheta0$ iff $c\neq0$ and $c\leqtheta0$. One requires that, for each $\theta$, the germs $\cL_{\leq c,\theta}$ form an exhaustive increasing filtration of~$\cL_\theta$ with respect to $\leqtheta$.
\item\label{enum:Stokesfil2}
Because the order $\leqtheta$ is open with respect to $\theta$, the germs $\cL_{<c,\theta}\defin\sum_{c'\letheta c}\cL_{\leq c',\theta}$ glue as a subsheaf $\cL_{<c}$ of~$\cL$. One requires that the graded sheaves $\gr_c\cL\defin\cL_{\leq c}/\cL_{<c}$ are locally constant sheaves on $S^1$.
\item\label{enum:Stokesfil3}
Near any $\re^{\ri\theta}\in S^1$, one requires that there are local isomorphisms $(\cL,\cL_\bbullet)\simeq(\gr\cL,(\gr\cL)_\bbullet)$, where the Stokes filtration on $\gr\cL\defin\bigoplus_{c\in\CC}\gr_c\cL$ is the natural one, that is, $(\gr\cL)_{\leq c,\theta}=\bigoplus_{c'\leqtheta c}\gr_{c'}\cL_\theta$. In particular, $\gr_c\cL=0$ except for~$c$ in a finite set $C\subset\CC$, called the set of exponential factors of the Stokes filtration $(\cL,\cL_\bbullet)$.
\end{enumerate}

We will say that the Stokes filtration is of \emph{pure Gaussian type} if it satisfies moreover:
\begin{enumerate}\setcounter{enumi}{3}
\item\label{enum:Stokesfil4}
The local system $\cL$ is \emph{constant} and $C\subset\CC\moins\{0\}$
\end{enumerate}

\begin{remarques}\label{rem:Stokesdir}\mbox{}
\begin{enumerate}
\item\label{rem:Stokesdir1}
The general definition of a Stokes filtration is somewhat more complicated, taking into account ramification, which does not occur in the present setting. In particular, the Stokes filtrations of Gaussian type are non ramified and their set $\Phi$ of exponential factors takes the form $C/2t^{\prime2}$ near at the origin of the chart $\Afu_{t'}$. Such a Stokes filtration may have a ``regular component'', while a Stokes filtration of pure Gaussian type does not.

\item\label{rem:Stokesdir4}
It is easy to check from the local grading property \eqref{enum:Stokesfil3} above that the family $(\cL_{\leq c})_{c\in\CC}$ can be recovered from the family $(\cL_{<c})_{c\in\CC}$ by the formula:
\[
\cL_{\leq c,\theta}=\bigcap_{c\letheta c'}\cL_{<c',\theta},\quad\theta\in S^1.
\]
In such a way, we could as well define a Stokes filtration of Gaussian type as a family of subsheaves $(\cL_{<c})_{c\in\CC}$ of $\cL$ such that, defining $\cL_{\leq c}$ by the formula above, the corresponding family $(\cL_{\leq c})_{c\in\CC}$ satisfies Properties \eqref{enum:Stokesfil1}--\eqref{enum:Stokesfil3} above. It will be simpler, when computing Laplace transforms of Stokes-filtered local systems, to adopt this point of view, due to Lemma \ref{lem:DRmodv/u}.

\item\label{rem:Stokesdir2}
For each pair $c\neq c'\in\CC$, there are exactly four values of $\theta\bmod2\pi$, say $\theta_{c,c'}^{(\nu)}$ ($\nu\in\ZZ/4\ZZ$), such that~$c$ and $c'$ are not comparable at $\theta$. We have $\theta^{(\nu+1)}_{c,c'}=\theta^{(\nu)}_{c,c'}+\pi/2$. These values are called the Stokes directions of the pair $(c,c')$. For any $\theta$ in one component of $\RR/2\pi\ZZ\moins\{\theta^{(\nu)}_{c,c'}\mid \nu\in\ZZ/4\ZZ\}$, we have $c\letheta c'$, and the reverse inequality for any $\theta$ in the next component. We denote the images of these intervals in $S^1$ via $\theta\mto \re^{\ri\theta}$ by $S^1_{c\leq c'}$ and $S^1_{c'\leq c}$ respectively. If $c=c'$, we set $S^1_{c\leq c}\defin S^1$.
\item\label{rem:Stokesdir3}
For each pair $c,c_o\in\CC$, the inclusion $j_{c\leq c_o}:S^1_{c\leq c_o}\hto S^1$ is open. We will denote by $\beta_{c\leq c_o}$ the functor $j_{c\leq c_o,!}j_{c\leq c_o}^{-1}$, consisting in restricting a sheaf to this open set and extending its restriction by $0$ to get a new sheaf on $S^1$. The filtration condition \eqref{enum:Stokesfil1} above implies that, for each pair $c,c_o$, there is a natural monomorphism $\beta_{c\leq c_o}\cL_{\leq c}\hto\cL_{\leq c_o}$.

\item\label{rem:Stokesdir5}
Let $\cF^*$ be the constant sheaf on $\Afu_t$ of rank $r$ and let $\cF=\wtj_*\cF^*$, if $\wtj$ denotes the open inclusion $\Afu_t\hto\wt\PP^1_t$ complementary to $S^1_\infty\hto\wt\PP^1_t$. Set $\cL=\cF_{|S^1_\infty}$. A~Stokes filtration $\cL_\bbullet$ of $\cL$ determines a family of subsheaves $\cF_\bbullet$ of $\cF$ by gluing $\cF^*$ with $\cL_\bbullet$. It~will be convenient to set $\cF^*_{\leq0}=\cF^*$ and $\cF^*_{<0}=0$, so that $\cF_{\leq0}$ restricts to $\cF^*$ on~$\Afu_t$, while $\cF_{<0}$ is supported on $S^1_\infty$. We call $(\cF,\cF_\bbullet)$ a Stokes-filtered sheaf.
\end{enumerate}
\end{remarques}

A morphism $\lambda:(\cL,\cL_\bbullet)\to(\cL',\cL'_\bbullet)$ of Stokes-filtered local systems is a morphism of local systems satisfying $\lambda(\cL_{\leq c})\subset\cL'_{\leq c}$ for each $c\in\CC$.

By a \emph{$C$-good closed interval $I\subset\RR/2\pi\ZZ$}, we mean a closed interval containing \emph{exactly one} Stokes direction for \emph{each pair $c\neq c'$ in $C$}, and such that each such Stokes direction belongs to the interior of $I$. Below, we will only use $C$-good closed intervals which are the image in $\RR/2\pi\ZZ$ of an interval $[\theta_o, \theta_o+\pi/2]$, where $\theta_o$ is not a Stokes direction for any pair $c\neq c'$ in~$C$.

\begin{proposition}[{\cf\cite[Prop\ptbl2.2]{H-S09}}]\label{prop:strict}\mbox{}
\begin{enumerate}
\item\label{prop:strict1}
On any $C$-good closed interval $I\subset \RR/2\pi\ZZ$, there exists a unique splitting $\cL_{|I}\simeq\bigoplus_c\gr_c\cL_{|I}$ compatible with the Stokes filtrations. With respect to this splitting, we have $\cL_{\leq c_o|I}=\bigoplus_{c\in C}\beta_{c\leq c_o}\gr_c\cL_{|I}$.
\item\label{prop:strict2}
Let $\lambda:(\cL,\cL_\bbullet)\to(\cL',\cL'_\bbullet)$ be a morphism of Stokes-filtered local systems of type~$C$. Then, for any $C$-good closed interval $I\subset \RR/2\pi\ZZ$, the morphism $\lambda_{|I}$ is graded with respect to the splittings in \eqref{prop:strict1}.
\end{enumerate}
\end{proposition}

(See also \cite[Chap.\,3]{Bibi10} for the arguments.)

\begin{remarque}\label{rem:GammaI}
On any $C$-good closed interval $I\subset \RR/2\pi\ZZ$, we have $\Gamma(I,\beta_{c\leq c_o}\gr_c\cL_{|I})=\nobreak0$ if $c\neq c_o$. If $c=c_o$, we have $\beta_{c\leq c_o}\gr_c\cL_{|I}=\gr_c\cL_{|I}$. We conclude from Proposition \ref{prop:strict}\eqref{prop:strict1} that $\Gamma(I,\cL_{\leq c_o})=\Gamma(I,\gr_{c_o}\cL_{|I})$.
\end{remarque}

\begin{proposition}\label{prop:catStokesfiltabelian}
The category of Stokes-filtered local systems $(\cL,\cL_\bbullet)$ of Gaussian type, and the full sub-category of Stokes-filtered local systems of pure Gaussian type, are abelian.\qed
\end{proposition}

\subsection{Stokes data}\label{subsec:Stokesdata}
These are linear data which provide a description of a~Stokes-filtered local system of Gaussian type. Let~$C$ be a non-empty finite subset of~$\CC$. We say that $\theta_o\in\RR/2\pi\ZZ$ is \emph{generic} with respect to $C$ if it is not a Stokes direction (\cf Remark \ref{rem:Stokesdir}\eqref{rem:Stokesdir2}) for any pair $c\neq c'\in C$. Once~$\theta_o$ generic with respect to~$C$ is chosen, there is a unique numbering of the set~$C$ in such a way that \hbox{$c_1\lethetao c_2\lethetao\cdots\lethetao c_n$}. We will~set
\begin{equation}\label{eq:thetaonu}
\theta^{(\nu)}_o=\theta_o+\nu\pi/2\quad (\nu\in\ZZ/4\ZZ).
\end{equation}
When $\theta$ varies in the good closed interval $[\theta^{(\nu)}_o,\theta^{(\nu+1)}_o]$, the order between~$c$ and~$c'$ changes exactly once for each pair $c\neq c'$, so the order on $C$ at $\theta^{(\nu+1)}_o$ is exactly reversed with respect to the order at $\theta^{(\nu)}_o$. In the following, we will use the reference order of~$C$ at $\theta_o=\theta_o^{(0)}$.

\begin{definition}[First definition]\label{def:catStokesdata}
Let~$C$ be a non-empty finite subset of $\CC$ and let $\theta_o\in\RR/2\pi\ZZ$ be generic with respect to $C$.
\begin{itemize}
\item
The category of Stokes data of Gaussian type~$C$ totally ordered by~$\theta_o$ (we also say \emph{of Gaussian type $(C,\theta_o)$}) has objects consisting of four families of finite dimensional $\kk$-vector spaces $(G_c^{(\nu)})_{c\in C}$ ($\nu\in\ZZ/4\ZZ$) and a diagram of morphisms
\begin{starequation}\label{eq:catStokesdata}
\begin{array}{c}
\xymatrix@C=1.5cm{
\bigoplus_{c\in C}G_c^{(1)}\ar[d]_{S^{(2,1)}}&\bigoplus_{c\in C}G_c^{(0)}\ar[l]_-{S^{(1,0)}}\\
\bigoplus_{c\in C}G_c^{(2)}\ar[r]_{S^{(3,2)}}&\bigoplus_{c\in C}G_c^{(3)}\ar[u]_-{S^{(0,3)}}
}
\end{array}
\end{starequation}%
such that, for the numbering $C=\{c_1,\dots,c_n\}$ defined by~$\theta_o$, $S^{(\nu,\nu-1)}=(S^{(\nu,\nu-1)}_{ij})_{i,j=1,\dots, n}$ is block-upper triangular for $\nu=0,2$ (\resp block-lower triangular for $\nu=1,3$), \ie $S^{(\nu,\nu-1)}_{ij}:G_{c_j}^{(\nu-1)}\to G_{c_i}^{(\nu)}$ is zero unless $i\leq j$ (\resp $i\geq j$), and $S^{(\nu,\nu-1)}_{ii}$ is invertible (so $\dim G_{c_i}^{(\nu-1)}=\dim G_{c_i}^{(\nu)}$, and $S^{(\nu,\nu-1)}$ itself is invertible).
\item
We say that it is of \emph{pure Gaussian type} if $C\subset\CC\moins\{0\}$ and the \emph{monodromy} $T:=S^{(0,3)}S^{(3,2)}S^{(2,1)}S^{(1,0)}$ is equal to the identity.

\item
The \emph{formal monodromy} on the $c$-component is the isomorphism
\[
T_c=S^{(0,3)}_{c,c}S^{(3,2)}_{c,c}S^{(2,1)}_{c,c}S^{(1,0)}_{c,c}:G_c^{(0)}\isom G_c^{(0)}.
\]

\item
A morphism of Stokes data of type $(C,\theta_o)$ consists of a family of morphisms of $\kk$-vector spaces $\lambda_c^{(\nu)}:G_c^{(\nu)}\to G_c^{\prime(\nu)}$, $c\in C$, $\nu\in\ZZ/4\ZZ$, which are compatible with the corresponding diagrams \eqref{eq:catStokesdata}.
\end{itemize}
\end{definition}

One checks that for such a morphism, one has in particular,
\[
S^{(\nu,\nu-1)}_{c,c}\lambda_c^{(\nu-1)}=\lambda_c^{(\nu)}S^{(\nu,\nu-1)}_{c,c}.
\]

The category of Stokes data of (pure) Gaussian type $(C,\theta_o)$ is clearly abelian. Fixing bases in the spaces $G_c^{(\nu)}$, $c\in C$, $\nu\in\ZZ/4\ZZ$, allows one to present Stokes data by matrices $(\Sigma^{(\nu,\nu-1)})_{\nu\in\ZZ/4\ZZ}$ where $\Sigma^{(\nu,\nu-1)}=(\Sigma^{(\nu,\nu-1)}_{ij})_{i,j=1,\dots, n}$ is block-upper (\resp \nobreakdash-lower) triangular for $\nu=0,2$ (\resp $\nu=1,3$) and each $\Sigma^{(\nu,\nu-1)}_{ii}$ is invertible. A set $(\Sigma^{\prime(\nu,\nu-1)})_{\nu\in\ZZ/4\ZZ}$ is equivalent to a set $(\Sigma^{(\nu,\nu-1)})_{\nu\in\ZZ/4\ZZ}$ if there exist invertible block-diagonal matrices $(\Lambda^{(\nu)})_{\nu\in\ZZ/4\ZZ}$ such that, for each pair $(i,j)$,
\[
\Sigma^{\prime(\nu,\nu-1)}_{ij}=\Lambda^{(\nu)}_i\Sigma^{(\nu,\nu-1)}_{ij}(\Lambda^{(\nu-1)}_j)^{-1}\qquad\forall\nu\in\ZZ/4\ZZ.
\]
In particular, up to equivalence, we may assume that $\Sigma^{\prime(\nu,\nu-1)}_{ii}=\id$ for all $i=1,\dots, n$ and $\nu=1,2,3$. We then have $\Sigma^{\prime(0,3)}_{ii}=T_i$ (matrix of the ``formal'' monodromy of~$G^{(0)}_{c_i}$). This leads to a variant of the first definition.

\begin{definition}[Variant of the first definition]\label{def:catStokesdatavariant}
Let~$C$ be a non-empty finite subset of~$\CC$ and let $\theta_o\in\RR/2\pi\ZZ$ be generic with respect to $C$. The category of Stokes data of Gaussian type ~$C$ totally ordered by~$\theta_o$ has objects consisting of a family of finite dimensional $\kk$-vector spaces $(G_c)_{c\in C}$, an automorphism~$T_c$ of~$G_c$ for each $c\in C$, and a diagram of morphisms
\begin{starequation}\label{eq:catStokesdatavariant}
\begin{array}{c}
\xymatrix@C=1.5cm{
\bigoplus_{c\in C}G_c\ar[d]_{S^{(2,1)}}&\bigoplus_{c\in C}G_c\ar[l]_-{S^{(1,0)}}\\
\bigoplus_{c\in C}G_c\ar[r]_{S^{(3,2)}}&\bigoplus_{c\in C}G_c\ar[u]_-{S^{(0,3)}}
}
\end{array}
\end{starequation}
such that, for the numbering $C=\{c_1,\dots,c_n\}$ defined by~$\theta_o$,
\begin{enumerate}
\item
$S^{(\nu,\nu-1)}_{ii}=\id$ for each $i\in\{1,\dots,n\}$ and $\nu\in\ZZ/4\ZZ$,
\item
$S^{(\nu,\nu-1)}=(S^{(\nu,\nu-1)}_{ij})_{i,j=1,\dots, n}$ is block-upper triangular for $\nu=0,2$ (\resp block-lower triangular for $\nu=1,3$), \ie $S^{(\nu,\nu-1)}_{ij}:G_{c_j}^{(\nu-1)}\to G_{c_i}^{(\nu)}$ is zero unless $i\leq j$ (\resp $i\geq j$),
\item
In the case of \emph{pure Gaussian type}, we moreover ask that $0\notin C$ and that the \emph{monodromy} $T:=\diag(T_1,\dots,T_n)\cdot S^{(0,3)}S^{(3,2)}S^{(2,1)}S^{(1,0)}$ is equal to the identity.
\end{enumerate}
\end{definition}

As we have seen at the matrix level, each Stokes data as in Definition \ref{def:catStokesdata} is isomorphic to one with all $G_c^{(\nu)}$ identified to the same $G_c$ and $S^{(\nu,\nu-1)}_{ii}=\id$ for $\nu=1,2,3$. In order to get Stokes data as in Definition \ref{def:catStokesdatavariant}, we then set $T_i=S^{(0,3)}_{ii}$ and the new $S^{(0,3)}$ is $\diag(T_1,\dots,T_n)^{-1}\cdot S^{(0,3)}$.

A morphism of Stokes data as in Definition \ref{def:catStokesdatavariant} consists of a family $(\lambda_c^{(\nu)})$ which is compatible with the diagrams \eqref{eq:catStokesdatavariant} and such that $T_c\lambda_c^{(3)}=\lambda_c^{(0)}T_c$ for each $c\in C$.

In the case where we assume that the monodromy is equal to the identity (but maybe $0\in C$), we can present the definition in another way.

\begin{definition}[Second definition]\label{def:catStokesdatabis}
The category of Stokes data of Gaussian type $(C,\theta_o,T=\id)$ has objects consisting of a finite dimensional $\kk$-vector space $L$ and, for each $\nu$, of an exhaustive filtration $L_{\leqnu\cbbullet}$ indexed by $\{1,\dots,n\}$ equipped with the order $\leqnu$ (order defined by $\theta_o$ if $\nu$ is even and reverse order at $\theta_o$  if $\nu$ is odd) which satisfy:
\begin{itemize}
\item
for each $\nu\in\ZZ/4\ZZ$, the filtrations $L_{\leqnu\cbbullet}$ and $L_{\leqnup\cbbullet}$ are \emph{opposite}, that is,
\begin{starequation} \label{eq:catStokesdatabis}
\forall\nu\in\ZZ/4\ZZ,\quad L= \bigoplus_{i=1,\dots,n}L_{\leqnu i}\cap L_{\leqnup i}.
\end{starequation}
\end{itemize}
A morphism between Stokes data of type $(C,\theta_o,T=\id)$ is a morphism between the corresponding vector spaces which is compatible with all filtrations.
\end{definition}

In this definition, exhaustive means that $L_{\leqnu\max_\nu}=L$, where $\max_\nu=n$ for $\nu$ even and $1$ for $\nu$ odd. If we denote by $L_{\lenu i}=\sum_{j\leqnu i,\; j\neq i}L_{\leqnu j}$, then we set $L_{\lenu\min_\nu}=0$ by using the convention that the sum over the empty set is zero.

\subsubsection*{Comparison of Definitions \ref{def:catStokesdata} and \ref{def:catStokesdatabis}}
The correspondence with Definition \ref{def:catStokesdata} goes as follows. For each $\nu$, the opposite filtrations define a grading $\bigoplus_iG_i^{(\nu)}$ with $G_i^{(\nu)}=L_{\leqnu i}\cap L_{\leqnup i}$, together with a canonical isomorphism $\bigoplus_iG_i^{(\nu)}\isom L$. The morphism $S^{(\nu,\nu-1)}$ is the morphism induced by $\id_L$ through the two consecutive isomorphisms at the levels $\nu-1$ and $\nu$. The morphism~$S^{(\nu,\nu-1)}$ is compatible, by definition, with the filtrations induced by $L_{\leqnu\cbbullet}$ on both sides. If $\nu$ is even (\resp odd), this filtration is increasing (\resp decreasing), which means that $S^{(\nu,\nu-1)}$ is block-upper (\resp -lower) triangular. The product of the morphisms $S^{(\nu,\nu-1)}$ is conjugate to $\id_L$, hence is equal to the identity. It remains to check that $S^{(\nu,\nu-1)}$ is a strict isomorphism (with respect to the filtration $L_{\leqnu\cbbullet}$). This is clear because $\id_L$ is obviously strict, and by definition the isomorphism $L\simeq\bigoplus_iG_i^{(\nu)}$ is strict with respect to both filtrations $L_{\leqnu\cbbullet}$ and $L_{\leqnup\cbbullet}$ for each $\nu$.

Conversely, given Stokes data of type $(C,\theta_o)$ as in Definition \ref{def:catStokesdata} with the supplementary assumption $T=\id$, we set $\theta_o^{(0)}=\theta_o$ and we define $L$ as $\bigoplus_{i=1,\dots,n}G_i^{(0)}$. The natural increasing (\resp decreasing) filtration on $L$ induced by this gradation is denoted by $L_{\leqnu\cbbullet}$ with $\nu=0$ (\resp $\nu=1$). The increasing (\resp decreasing) filtration $L_{\leqnu\cbbullet}$ with $\nu=2$ (\resp $\nu=3$) is the filtration induced by the increasing (\resp decreasing) filtration attached to the grading $\bigoplus_{i=1,\dots,n}G_i^{(2)}$ through any of the isomorphisms $(S^{(2,1)}S^{(1,0)})^{-1}$, $S^{(0,3)}S^{(3,2)}$ from $\bigoplus_{i=1,\dots,n}G_i^{(2)}$ to $L$.

\begin{exemple}[Trivial Stokes data]\label{exem:trivialStokesdata}
Let $c_o\in\CC$ and let $L$ be a finite dimensional $\kk$\nobreakdash-vector space. Then the trivial Stokes data of Gaussian type $(\{c_o\},\theta_o,T=\id)$ and of exponent $c_o$ are the Stokes data defined by
\[
L_{\leqnu c}=\begin{cases}
L&\text{if }c_o\leqnu c,\\
0&\text{if }c\lenu c_o.
\end{cases}
\]
\end{exemple}

\begin{exemple}[Adding trivial Stokes data]
Let $(L,(L_{\leqnu\cbbullet})_{\nu\in\ZZ/4\ZZ}$ be Stokes data of Gaussian type $(C,\theta_o,T=\id)$. Let $c_0\in\CC\moins C$ be such that $c_0\leqnu c$ for all $c\in C$ and for $\nu$ even, while $c\leqnu c_0$ for all $c\in C$ and $\nu$ odd. We will then set $C':=C\cup\{c_0\}=\{c_0,c_1,\dots,c_n\}$ with respect to the ordering at $\theta_o$. For example, if $C\subset\RR$ and $\cos2\theta_o>0$, one can choose $c_0\in\RR$ such that $c_0<c$ for all $c\in\CC$. Let $L_o$ be a finite dimensional $\kk$-vector space. One obtains Stokes data of Gaussian type $(C',\theta_o,T=\id)$ in the following way:
\begin{itemize}
\item
the vector space is $L\oplus L_o$,
\item
for $\nu$ even,
\[
(L\oplus L_o)_{\leqnu c_0}=0\oplus L_o,\quad(L\oplus L_o)_{\leqnu c_1}=L_{\leqnu c_1}\oplus L_o,\dots, (L\oplus L_o)_{\leqnu c_n}=L_{\leqnu c_n}\oplus L_o,
\]
\item
and for $\nu$ odd,
\[
(L\oplus L_o)_{\leqnu c_0}=L\oplus L_o,\quad(L\oplus L_o)_{\leqnu c_1}=L_{\leqnu c_1}\oplus 0,\dots, (L\oplus L_o)_{\leqnu c_n}=L_{\leqnu c_n}\oplus 0.
\]
\end{itemize}
We have an exact sequence
\begin{starequation}\label{eq:extStokes}
0\to (L,L_{\leqnu\cbbullet})\to(L\oplus L_o,(L\oplus L_o)_{\leqnu\cbbullet})\to (L_o,L_{o,\leqnu\cbbullet})\to0,
\end{starequation}%
where $(L_o,L_{o,\leqnu\cbbullet})$ are the trivial Stokes data of Gaussian type $(\{c_0\},\theta_o,T=\id)$.
\end{exemple}

\subsection{Stokes data attached to a Stokes-filtered local system}
We will now define a functor (depending on~$\theta_o$) from the category of Stokes-filtered constant local systems of pure Gaussian type~$C$ to the category of Stokes data of pure Gaussian type $(C,\theta_o)$, and we will show that it is an equivalence.

Let us also fix intervals~$I^{(\nu)}=[\theta^{(\nu)}_o,\theta^{(\nu+1)}_o]$ of length $\pi/2$ on $\RR/2\pi\ZZ$ so that the intersection $I^{(\nu)}\cap I^{(\nu+1)}$ consists of the point $\theta^{(\nu+1)}_o$, which is not a Stokes direction for any pair $c\neq c'\in C$ (recall that $\theta^{(\nu)}_o$ is defined by \eqref{eq:thetaonu}).

To a constant local system~$\cL$ on $S^1$ we attach the following ``monodromy data'':
\begin{enumerate}
\item
a vector space $L=\Gamma(S^1,\cL)$,
\item
vector spaces $L^{(\nu)}=\Gamma(I^{(\nu)},\cL)$ ($\nu\in\ZZ/4\ZZ$),
\item
vector spaces $L_{\theta^{(\nu)}_o}=\Gamma(I^{(\nu-1)}\cap I^{(\nu)},\cL)\simeq\cL_{\theta^{(\nu)}_o}$,
\item
a diagram of natural restriction isomorphisms $a_\nu^{(\nu)},a_\nu^{(\nu+1)},b^{(\nu)}$,
\begin{equation}\label{eq:diagL}
\begin{array}{c}
\xymatrix@R=3mm@C=5mm{
&&&L^{(1)}\ar[dll]^-{a_1^{(2)}}\ar[drr]_-{a_1^{(1)}}\ar@/_2pc/[lllddd]_-{S^{(2,1)}}&&&\\
&L_{\theta_o^{(2)}}
&&&&L_{\theta_o^{(1)}}
&\\
&&&&&&\\
L^{(2)}\ar[uur]_-{a_2^{(2)}}\ar[ddr]^-{a_2^{(3)}}\ar@/_2pc/[rrrddd]_-{S^{(3,2)}}&&&L\ar[uurr]_{b^{(1)}}\ar[uull]^{b^{(2)}}\ar[ddll]_{b^{(3)}}\ar[ddrr]^{b^{(0)}}&&&L^{(0)}\ar[uul]^-{a_0^{(1)}}\ar[ddl]_-{a_0^{(0)}}\ar@/_2pc/[llluuu]_-{S^{(1,0)}}\\
&&&&&&\\
&L_{\theta_o^{(3)}}
&&&&L_{\theta_o^{(0)}}
&\\
&&&L^{(3)}\ar[urr]^-{a_3^{(0)}}\ar[ull]_-{a_3^{(3)}}\ar@/_2pc/[rrruuu]_-{S^{(0,3)}}&&&
}
\end{array}
\end{equation}
\end{enumerate}

We will use the following description: $(L^{(\nu)},S^{(\nu,\nu-1)})_{\nu\in\ZZ/4\ZZ}$, with isomorphisms $S^{(\nu,\nu-1)}:L^{(\nu-1)}\isom L^{(\nu)}$ and monodromy $T^{(0)}=\id:L^{(0)}\isom L^{(0)}$, where
\[
S^{(\nu,\nu-1)}=(a_\nu^{(\nu)})^{-1}a_{\nu-1}^{(\nu)},\quad T^{(0)}=S^{(0,3)}S^{(3,2)}S^{(2,1)}S^{(1,0)}.
\]

Assume now that $(\cL,\cL_\bbullet)$ is a Stokes-filtered local system with associated graded local system $\gr\cL=\bigoplus_{i=1}^n\gr_{c_i}\cL$.

\begin{definition}[Stokes data attached to $(\cL,\cL_\bbullet)$]\label{def:stokescL}
The filtration $\cL_{\leq c,\theta_o^{(\nu)}}$ of the germ $\cL_{\theta_o^{(\nu)}}$ induces a filtration on~$L_{\theta_o^{(\nu)}}$ which is $\theta_o^{(\nu)}$-increasing. Through $b^{(\nu)}$, $L$ comes equipped with a filtration $L_{\leqnu \cbbullet}$.
\end{definition}

Let us check that the pair of filtrations $L_{\leqnu \cbbullet}$ and $L_{\leqnup \cbbullet}$ are \emph{opposite}. It is enough to check this on $L^{(\nu)}:=\Gamma(I^{(\nu)},\cL)\simeq L$, which is similarly equipped with filtrations $L^{(\nu)}_{\leqnu \cbbullet}$ and $L^{(\nu)}_{\leqnup \cbbullet}$. We identify $L^{(\nu)}_{\leqnu c_o}\cap L^{(\nu)}_{\leqnup c_o}$ with $\Gamma(I^{(\nu)},\cL_{\leq c_o})$, and Remark \ref{rem:GammaI} identifies this space with $\Gamma(I^{(\nu)},\gr_{c_o}\cL)$. We thus obtain Stokes data of type $(C,\theta_o)$ as in Definition \ref{def:catStokesdatabis}.

In such a way, we have defined the desired functor (to check the compatibility with morphisms, use Proposition \ref{prop:strict}\eqref{prop:strict2}).

As a consequence of the previous discussion we can state the following classical result (the bijection at the level of $\Hom$ follows from Proposition \ref{prop:strict}\eqref{prop:strict2}):

\begin{proposition}\label{prop:datafiltered}
The previous functor is an equivalence between the category of constant Stokes-filtered local systems of pure Gaussian type~$C$ and the category of Stokes data of pure Gaussian type $(C,\theta_o)$.\qed
\end{proposition}

\Subsection{A quasi-inverse functor: the sheaves \texorpdfstring{$\cL_{\leq c_o}$}{cL} in terms of Stokes data}\label{subsec:LleqStokes}
Let us fix $c_o\in\CC$. We will make explicit the description of the subsheaves $\cL_{<c_o}$ and $\cL_{\leq c_o}$ of $\cL$ in terms of Stokes data.

Let us start with Stokes data of type $(C,\theta_o)$ as in Definition \ref{def:catStokesdata}. The (constant) sheaf $\cL$ is obtained by the following gluing procedure, with respect to the covering~$(I^{(\nu)})_{\nu\in\ZZ/4\ZZ}$. We set
\[
\cL_{I^{(\nu)}}=\bigoplus_{c\in C}\Big(\CC_{I^{(\nu)}}\otimes_\CC G_c^{(\nu)}\Big),
\]
and we define the gluing isomorphisms
\[
g^{(\nu,\nu-1)}:\cL_{I^{(\nu-1)}}\big|_{I^{(\nu-1)}\cap I^{(\nu)}}\isom \cL_{I^{(\nu)}}\big|_{I^{(\nu-1)}\cap I^{(\nu)}}
\]
as $\id_{\CC_{I^{(\nu-1)}\cap I^{(\nu)}}}\otimes S^{(\nu,\nu-1)}$.

With this presentation, the subsheaf $\cL_{\leq c_o}$ is defined by the data of the subsheaves
\[
\cL_{I^{(\nu)},\leq c_o}=\bigoplus_{c\in C} \Big((\beta_{c\leq c_o}\CC_{I^{(\nu)}})\otimes_\CC G_c^{(\nu)} \Big),
\]
and the gluing isomorphisms induced by $g^{(\nu,\nu-1)}$ (and similarly for $\cL_{<c_o}$ by replacing $\leq$ with $<$). The only point to check is that the gluing isomorphisms preserve these subsheaves. We will check this for each $(c,c')$-component of $g^{(\nu,\nu-1)}$.

\begin{itemize}
\item
If $c\leqnu c'$, we have the implication $c'\leqnu c_o\implique c\leqnu c_o$, and the identity morphism $\id_{\CC_{I^{(\nu-1)}\cap I^{(\nu)}}}$ sends $\beta_{c'\leq c_o}\CC_{I^{(\nu)}\cap I^{(\nu-1)}}$ to $\beta_{c\leq c_o}\CC_{I^{(\nu)}\cap I^{(\nu-1)}}$, since either both sheaves are equal to $\CC_{I^{(\nu)}\cap I^{(\nu-1)}}$ if $c'\leqnu c_o$, or the first one is zero otherwise.
\item
Otherwise, we have $S^{(\nu,\nu-1)}_{c,c'}=0$, and thus the $(c,c')$-component of $g^{(\nu,\nu-1)}$ is zero.
\end{itemize}

\begin{proposition}\mbox{}
\begin{enumerate}
\item
For each $c_o\!\in\!\CC$, we have $H^k(S^1,\cL_{<c_o})\!=\!0$ for $k\!\neq\!1$ and $\dim H^1(S^1,\cL_{<c_o})\!=\nobreak\!\nobreak2r$.
\item
For $c_o\notin C$, we also have $H^k(S^1,\cL_{\leq c_o})=H^k(S^1,\cL_{<c_o})$ for each $k$.
\item
For $c_o\in C$, we have $H^k(S^1,\cL_{\leq c_o})=0$ for $k\geq2$ and
\[
-\chi(S^1,\cL_{\leq c_o})=2r+\dim\coker(T_{c_o}-\id)-\dim\ker(T_{c_o}-\id),
\]
and the dimension of the space of global sections $H^0(S^1,\cL_{\leq c_o})$ is equal to the dimension of the subspace of $G_{c_o}^{(0)}$ defined as
\begin{multline*}
\bigcap_{c\neq c_o}\ker S_{c,c_o}^{(1,0)}\cap \bigcap_{c\neq c_o}\ker S_{c,c_o}^{(2,1)}S_{c_o,c_o}^{(1,0)}\cap \bigcap_{c\neq c_o}\ker S_{c,c_o}^{(3,2)}S_{c_o,c_o}^{(2,1)}S_{c_o,c_o}^{(1,0)}\\
\cap \bigcap_{c\neq c_o}\ker S_{c,c_o}^{(0,3)}S_{c_o,c_o}^{(3,2)}S_{c_o,c_o}^{(2,1)}S_{c_o,c_o}^{(1,0)},
\end{multline*}
that is, the intersection of the kernels of the $S_{c,c_o}^{(\nu,\nu-1)}$ ($c\neq c_o$, $\nu\in\ZZ/4\ZZ$) when regarded in a natural way in the same $G_{c_o}^{(0)}$.
\end{enumerate}
\end{proposition}

Up to equivalence, we can assume that $S^{(\nu,\nu-1)}_{c_o,c_o}=\id$ for each $c_o\in C$ and $\nu=1,2,3$, so that for $c$ fixed the spaces $G_c^{(\nu)}$ ($\nu\in\ZZ/4\ZZ$) are all canonically identified to a space that we denote by~$G_c$. We the regard $S_{c,c_o}^{(\nu,\nu-1)}$ as a morphism $G_{c_o}\to G_c$. With such a representative, we have
\[
\dim H^0(S^1,\cL_{\leq c_o})=\dim\bigcap_{\nu\in\ZZ/4\ZZ}\bigcap_{c\neq c_o}\ker S_{c,c_o}^{(\nu,\nu-1)}.
\]

\begin{proof}
We will assume for the sake of simplicity that $\theta_o$ is generic with respect to $C\cup\{c_o\}$, that is, $\theta_o$ is not a Stokes direction for any pair $c\neq c'$ in $C\cup\{c_o\}$. We then have
\[
\cL_{\leq c_o,|I^{(\nu)}\cap I^{(\nu-1)}}=\bigoplus_{\substack{c\in C\\ c\leqnu c_o}} \Big(\CC_{I^{(\nu)}\cap I^{(\nu-1)}}\otimes_\CC G_c^{(\nu)} \Big).
\]

The closed covering $(I^{(\nu)})_{\nu\in\ZZ/4\ZZ}$ is a Leray covering for $\cL_{\leq c_o}$ or $\cL_{<c_o}$ since
\begin{itemize}
\item
on $I^{(\nu)}$, $H^k(I^{(\nu)},\beta_{c\leq c_o}\CC)=0$ for all $k$ if $c\neq c_o$ and for all $k\geq1$ if $c=\nobreak c_o$; similarly, $H^k(I^{(\nu)},\beta_{c<c_o}\CC)=0$ for all $k$; this is because there is exactly one Stokes direction for $(c,c_o)$ in $I^{(\nu)}$ if $c\neq c_o$ and, if $c=\nobreak c_o$, because $\beta_{c<c_o}\CC=0$ and $\beta_{c\leq c_o}\CC=\nobreak\CC_{I^{(\nu)}}$;
\item
on $I^{(\nu-1)}\cap I^{(\nu)}$, the corresponding $H^k$ are zero for $k\geq1$ because $\cL_{\leq c_o}$ and $\cL_{<c_o}$ are sheaves.
\end{itemize}
We can compute the cohomology by using the \v{C}ech complex relative to this covering \cite[Cor\ptbl of Th\ptbl5.2.4, p\ptbl209]{Godement64}. Except for $\cL_{\leq c_o}$ with $c_o\in C$, this complex reduces to $\bigoplus_\nu\Gamma(I^{(\nu-1)}\cap I^{(\nu)},\cL_\bbullet)$ placed in degree one, with $\cL_\bbullet=\cL_{<c_o}$ (any $c_o$) or $\cL_{\leq c_o}$ ($c_o\notin C$). Since $C\cup\{c_o\}$ is totally ordered on each $I^{(\nu-1)}\cap I^{(\nu)}$, and since the order is reversed when we change $\nu$ to $\nu+1$, the sum of the dimensions of two consecutive terms in this direct sum is equal to~$r$.

Assume now that $c_o\in C$. Then the \v{C}ech complex for $\cL_{\leq c_o}$ has two terms:
\begin{equation}\label{eq:diffcechLco}
\bigoplus_\nu G_{c_o}^{(\nu)}\to\bigoplus_{c,\nu} G_c^{(\nu)}
\end{equation}
where the component $G_{c_o}^{(\nu-1)}\to\bigoplus_{c,\nu'} G_c^{(\nu')}$ of the differential takes values in $G_{c_o}^{(\nu-1)}\oplus\bigoplus_{c\leqnu c_o} G_c^{(\nu)}$ and is equal to $\id\oplus\bigoplus_{c\leqnu c_o}S_{c,c_o}^{(\nu,\nu-1)}$.

This implies that $H^2(S^1,\cL_{\leq c_o})=0$ and, by using the exact sequence
\begin{multline*}
0\to H^0(S^1,\cL_{\leq c_o})\to H^0(S^1,\gr_{c_o}\cL)\to H^1(S^1,\cL_{<c_o})\\
\to H^1(S^1,\cL_{\leq c_o})\to H^1(S^1,\gr_{c_o}\cL)\to0
\end{multline*}
we find the desired formula for $\chi(S^1,\cL_{\leq c_o})$. Let us now compute the kernel of \eqref{eq:diffcechLco}. If $(x_\nu)\in\bigoplus_\nu G_{c_o}^{(\nu)}$ belongs to the kernel then, considering the component of the image on $\bigoplus_\nu G_{c_o}^{(\nu)}$, we find $x_\nu=-S^{(\nu,\nu-1)}_{c_o,c_o}(x_{\nu-1})$ and, since $T=\id$, the kernel is isomorphic to the subspace of $G_{c_o}^{(0)}$ defined in the last point of the proposition.
\end{proof}

\subsection{The Riemann-Hilbert correspondence for differential systems of pure Gaussian type}
Let $M$ be of pure Gaussian type. The analytic de~Rham complex $\DR^\an M$ on $\Afu_t$ can be extended as a complex on $\wt\PP^1_t$ by considering the \emph{rapid decay} de~Rham complex $\DR^{\rrd\infty}M$ obtained by replacing (in the definition of $\DR$) holomorphic forms on $\Afu_t$ with holomorphic forms having rapid decay at infinity.

The Riemann-Hilbert-Deligne  correspondence associates to $M$ the (constant) sheaf $\cF^*=\cH^0\DR^\an M$ on $\Afu_t$ and the subsheaves $\cL_{<c}:=\cH^0\DR^{\rrd\infty}(M\otimes E^{ct^2/2})_{|S^1_\infty}$ of $\cL:=\cF_{|S^1_\infty}$ (\cf Remarks \ref{rem:Stokesdir}\eqref{rem:Stokesdir4} and \eqref{rem:Stokesdir5}). From \cite{Deligne78} (\cf also \cite{Malgrange83bb}, \cite{B-V89}, \cite{Malgrange91}), we get:

\begin{proposition}\label{prop:equivconnStokesloc}
There is an equivalence between the category of differential systems on $\Afu_t$ of pure Gaussian type at $t=\infty$, and the category of $\CC$-Stokes-filtered sheaves of pure Gaussian type on $\wt\PP^1_t$.\qed
\end{proposition}

We also have the following theorem (recall the notation in \S\ref{subsec:settingnotation}).

\begin{theoreme}
Let $(X,x_o)$ be a pointed simply connected complex manifold and let $C=(c_1,\dots,c_n):X\to(\CC^*)^n\moins\textup{diagonals}$ be a holomorphic map. Let $M_o$ be a differential system of pure Gaussian type with formal decomposition at infinity given by $\bigoplus_{i=1}^n(\cE^{-c_i(x_o)/2t^{\prime2}}\otimes R_{c_i(x_o)})$. Then there exists a unique locally free $\cO_X[t]$-module $(\cM,\nabla)$ with flat connection such that, denoting by $i_x$ the inclusion $\{x\}\hto X$,
\begin{enumerate}
\item
$i_x^+(\cM,\nabla)$ is a $\Clt$-module of Gaussian type with formal decomposition at infinity given by $\bigoplus_{i=1}^n(\cE^{-c_i(x)/2t^{\prime2}}\otimes R_{c_i(x_o)})$ (in particular, the $R_{c_i(x_o)}$ remain constant),
\item
$i_{x_o}^+(\cM,\nabla)=M_o$.
\end{enumerate}
\end{theoreme}

\begin{proof}
The formal $\cO_X\lpr t'\rpr$-module $\bigoplus_{i=1}^n(\cE^{-c_i(x)/2t^{\prime2}}\otimes R_{c_i(x_o)})$, together with its connection $\nabla$, is a formal meromorphic flat bundle which is \emph{good} since $c_i(x)\neq c_j(x)$ for each $i\neq j$ and each $x\in X$. From \cite[Cor.\,II.6.4]{Bibi00b} we deduce the exis\-tence and uniqueness of a meromorphic flat bundle $(\cM_\infty,\nabla)$ in an analytic neighbourhood of $X\times\infty\subset X\times\PP^1$. Similarly, the local system it induces on $X\times S^1_\infty$ is constant, and this allows us to glue $(\cM_\infty,\nabla)$ with $(\cO_{X\times\Afuan_t}^{\rk M_o},\rd)$. Lastly, by choosing a lattice, one can use a GAGA argument to make the construction algebraic with respect to~$t$.
\end{proof}
\begin{remarque}
One proves in a similar way that any morphism $\varphi_o:M_o\to N_o$ between differential systems of pure Gaussian type $C(x_o)$ extends in a unique way as a morphism $\varphi:(\cM,\nabla)\to(\cN,\nabla)$.
\end{remarque}

\begin{corollaire}\label{cor:equivCCprime}
Let $X$ be a simply connected open subset of $(\CC^*)^n\moins\textup{diagonals}$ and let $C,C'$ be two points of $X$. Then $X$ defines an equivalence of categories between the category of differential systems of pure Gaussian type $C$ and that of pure Gaussian type $C'$. 
\end{corollaire}

\begin{proof}
We apply the above theorem and remark to the inclusion map $X\hto(\CC^*)^n\moins\textup{diagonals}$, and both categories considered in the corollary are equivalent, via the restriction functor, to the category of flat holomorphic families parametrized by $X$ of differential systems of pure Gaussian type $X$.
\end{proof}

\section{Topological Laplace transformation}

\subsection{The theorem of Mochizuki}
Our objective is to compute the Stokes data of~$\wh M$ at $\tau=\infty$ in terms of the Stokes data of $M$ at $t=\infty$. More precisely, we look for a purely topological computation.

In other words we want to express $\wh\cF^*=\DR^\an\wh M$ and, for any $\gamma\in\CC^*$ (and particularly for $\gamma\in-1/C$), the rapid decay de~Rham complex $\FcL_{<\gamma}:=\DR^{\rrd\wh\infty}(\wh M\otimes\nobreak E^{\gamma\tau^2/2})_{S^1_{\wh\infty}}$ uniquely in terms of the Stokes-filtered sheaf $(\cF,\cF_\bbullet)$ attached to $M$. We consider the natural embedding of diagrams of projections
\[
\begin{array}{c}
\xymatrix{
&\Afu_t\times\Afu_\tau\ar[dl]_p\ar[dr]^{\wh p}\\
\Afu_t&&\Afu_\tau
}
\end{array}
\quad\hto\quad
\begin{array}{c}
\xymatrix{
&\wt\PP^1_t\times\wt\PP^1_\tau\ar[dl]_{\wt p}\ar[dr]^{\wt{\wh p}}\\
\wt\PP^1_t&&\wt\PP^1_\tau
}
\end{array}
\]
where $\wt\PP^1_t$ (\resp $\wt\PP^1_\tau$) is the circle completion of $\Afu_t$ (\resp $\Afu_\tau$), that is, the real oriented blow-up of $\PP^1_t$ at $\infty$ (\resp $\PP^1_\tau$ at $\wh\infty$). We set
\[
D_\infty:=\{\infty\}\times\PP^1_\tau,\quad D_{\wh\infty}:=\PP^1_t\times\{\wh\infty\}.
\]
The sheaves of holomorphic functions with moderate growth and rapid decay are thus well-defined on the space $\wt\PP^1_t\times\wt\PP^1_\tau$, and there are corresponding moderate and rapid decay de~Rham complexes. We will use the following notation. Let $X$ be a complex manifold and let $D$ be a divisor with normal crossings in $X$. We denote by $\wt X(D)$ or simply $\wt X$ the real oriented blow-up of $X$ along the irreducible components of $D$. On this space are defined the sheaves of holomorphic forms on $X\moins D$ having rapid decay (\resp moderate growth) along the pull-back of $D$ in $\wt X$. Correspondingly, one defines the de~Rham complex with rapid decay (\resp moderate growth) along $D$ for a $\cD_X$-module, that we denote by $\DR^\rdD$ (\resp $\DR^\modD$). These are complexes of sheaves on $\wt X$ (\cf \cite[Chap.\,8]{Bibi10}). The main tool for our purpose is the following theorem.

\begin{theoreme}[{Mochizuki \cite[Cor.\,4.51]{Mochizuki10}}]\label{th:mochi}
Let $e:X\to\PP^1_t\times\PP^1_\tau$ be a sequence of point blowing-ups above $(\infty,\wh\infty)$ and set $D_X=e^{-1}(D_\infty\cup D_{\wh\infty})$. Let $\wt X$ be the real blow-up space of $X$ along the irreducible components of $D_X$, and let $\wt e:\wt X\to\wt\PP^1_t\times\wt\PP^1_\tau$ be the corresponding map. Then, for each $\gamma\in\CC^*$ we have
\[
\FcL_{<\gamma}=\bR(\wt{\wh p}\circ\wt e)_*\DR^{\rdD_X}e^+(p^+M\otimes E^{-t\tau+\gamma\tau^2/2})_{|S^1_{\wh\infty}}[1].
\]
\end{theoreme}

\begin{remarque}
The use of this analytic theorem could be avoided in the present article, in order to obtain Theorem \ref{th:mainLaplace} below. However, in more complicated cases, this theorem allows a general definition of the topological Laplace transformation, and its use simplifies the presentation (\cf \cite{H-S14}). Note that this theorem implies that, for each $\gamma\in\CC$, the complex $\bR(\wt{\wh p}\circ\wt e)_*\DR^{\rdD_X}e^+(p^+M\otimes  E^{-t\tau+\gamma\tau^2/2})$ has cohomology in degree one at most. Let us consider the open inclusion $\wtj:X\moins D_X\hto\wt X$. Then the natural morphism from this complex to $\bR\wtj_*\wtj^{-1}$ of itself (which also has cohomology in degree one at most) is injective at the cohomology level.
\end{remarque}

This theorem, with $X=\PP^1_t\times\PP^1_\tau$, $D=D_\infty\cup D_{\wh\infty}$, and $e=\id$, reduces the initial problem to expressing $\DR^\rdD(p^+M\otimes  E^{-t\tau+\gamma\tau^2/2})$ in terms of $(\cF,\cF_\bbullet)$. It is not too difficult to express $\cH^0$ of this complexes in terms of $(\cF,\cF_\bbullet)$, as shown in Proposition~\ref{prop:Ggamma} below. However, one cannot apply Majima's asymptotic analysis to obtain the vanishing of the sheaves $\cH^k$, $k\geq1$, since $p^+M\otimes  E^{-t\tau+\gamma\tau^2/2}$ is not good along $D$ at $(\infty,\wh\infty)$, in the sense of \cite{Bibi97}. We will thus apply Theorem \ref{th:mochi} to a suitable modification~$X$.

We continue to use the notation as in Theorem \ref{th:mochi} and we set $q=p\circ e$. We also use the notation of Remark \ref{rem:Stokesdir}\eqref{rem:Stokesdir5}.

\begin{proposition}\label{prop:Ggamma}
Let $e:X\to\PP^1_t\times\PP^1_\tau$ be as in Theorem \ref{th:mochi}. We set, for $\gamma\in\CC$,
\[
\cG_{<\gamma}:=\cH^0\DR^{\rdD_X}e^+(p^+M\otimes E^{-t\tau+\gamma\tau^2/2})
\]
that we regard as a subsheaf of
\[
\cG:=\wtj_*\cH^0\DR^\an e^+(p^+M\otimes E^{-t\tau})
\]
in a natural way ($\wtj:X\moins D_X\hto\wt X$). Then
\begin{starequation}\label{eq:Ggammagen}
\wtj^{-1}\cG_{<\gamma}=\wtj^{-1}\cG=\wtj^{-1}\wt q^{-1}\cF_{\leq0}\quad\text{and}\quad\cG=\wtj_*\wtj^{-1}\wt q^{-1}\cF_{\leq0}.
\end{starequation}%
Moreover, for each $\wt x\in\wt X_{|e^{-1}(D_{\wh\infty})}$, setting $\wt\theta:=\wt{\wh q}(\wt x)$, we have
\begin{starstarequation}\label{eq:Ggamma}
\cG_{<\gamma,\wt x}=
\begin{cases}
\begin{cases}
\cF_t&\text{if $\gamma/2\tau^{\prime2}-t/\tau'\lethetah0$},\\
0&\text{otherwise},
\end{cases}&\text{if $t:=\wt q(\wt x)\in \Afu_t$},\\[15pt]
\dpl\sum_{\substack{c\in C\\ \gamma/2\tau^{\prime2}-1/t'\tau'-c/2t^{\prime2}<_{(\theta,\wh\theta)}0}}\hspace*{-10mm}\cL_{\leq c,\theta}&\text{if $\theta:=\wt q(\wt x)\in S^1_\infty$},
\end{cases}
\end{starstarequation}%
where the sum is taken in $\cL_{\theta}$.
\end{proposition}

\begin{proof}
This is a straightforward computation from the definition of horizontal sections of a flat meromorphic connection. Let us only indicate why $\wtj^{-1}\cG=\wtj^{-1}\wt q^{-1}\cF_\leq0$. We thus forget about $\gamma$. The twist by $E^{-t\tau}$ consists in changing the set of exponential factors \hbox{$\{ct^2/2\mid c\in C\}$} to \hbox{$\{(c+2\tau t')t^2/2\mid c\in C\}$}. We then have $\wtj^{-1}\cG=\sum_{c\in\CC,\;ct^2/2+t\tau\leq0}\cL_{\leq c}$ above $S^1_\infty$ and $\wtj^{-1}\cG=\cF^*$ above $\Afu_t$. Near any point of $S^1_\infty=\{|t|=\infty\}$ and for $\tau$ near $\tau_o\in\Afu_\tau$, we have $t\tau+ct^2/2=ct^2(1+2t'\tau/c)/2$ and $(1+2t'\tau/c)\sim1$. Therefore, at any $\theta\in S^1_\infty$, we have $t\tau+ct^2/2\leqtheta0$ iff $ct^2/2\leqtheta0$, and we conclude that $\wtj^{-1}\cG=\wt p^{-1}\cF_{\leq0}$.
\end{proof}

\subsection{Choice of the blow-up space}
\begin{figure}[htb]
\begin{center}
\includegraphics[scale=.5]{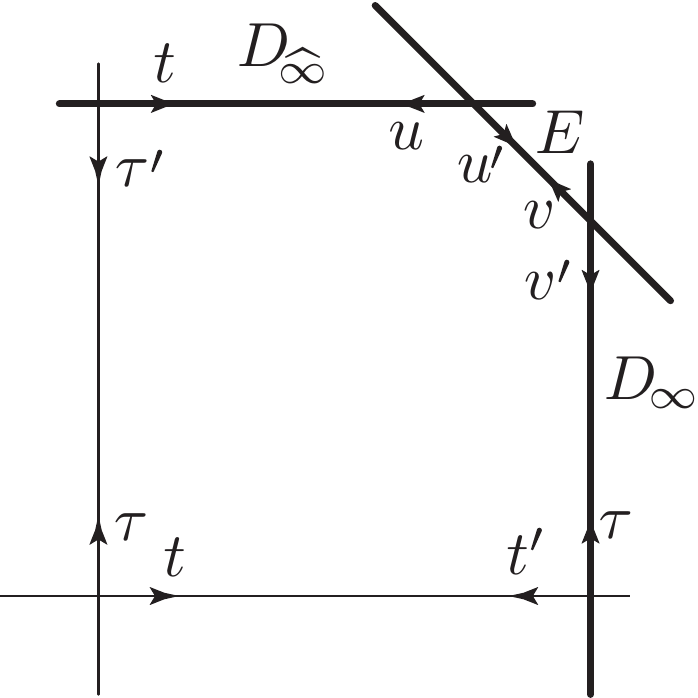}
\caption{The space $Z$}\label{fig:Z}
\end{center}
\end{figure}

We denote by $\epsilon:Z\to\PP^1_t\times\PP^1_\tau$ the blowing-up map with center $(\infty,\wh\infty)$, by $E$ the exceptional divisor and by $D$ the divisor $E\cup D_\infty\cup D_{\wh\infty}$, with $D_\infty:=(\{\infty\}\times\PP^1_\tau)$ and $D_{\wh\infty}:=(\PP^1_t\times\{\wh\infty\})$. We will also set
\[
E^*=E\moins((D_\infty\cap E)\cup(D_{\wh\infty}\cap E),\quad D_{\wh\infty}^*=D_{\wh\infty}\moins(D_{\wh\infty}\cap E),\quad D_\infty^*=D_\infty\moins(D_\infty\cap E).
\]
We notice that $\Afu_t\times\PP^1_\tau$ and $\PP^1_t\times\Afu_\tau$ are naturally included in $Z$ as open subsets, whose complementary sets are respectively $D_\infty\cup E$ and $D_{\wh\infty}\cup E$ in $Z$ (\cf Figure \ref{fig:Z}).

By definition, the pull-back of the affine chart $\Afu_{t'}\times\nobreak\Afu_{\tau'}$ of $\PP^1_t\times\PP^1_\tau$ centered at $(\infty,\wh\infty)$ consists of two affine charts, with coordinates denoted by $(u,u')$ and $(v,v')$, such that the projection is given by
\begin{itemize}
\item
$t'=u$, $\tau'=uu'$ in the first chart,
\item
$t'=vv'$, $\tau'=v'$ in the second chart.
\end{itemize}

\begin{lemme}\label{lem:zeros}
The meromorphic bundle with flat connection $\epsilon^+(p^+M\otimes\nobreak E^{-t\tau+\gamma\tau^2/2})$ is pseudo-good, in the sense that it is good on $D_Z$ minus a finite set in $E^*$, where the exponential factors which are not good behave like $v/u^2$ or $v^2/u^2$.
\end{lemme}

\begin{proof}
The possible exponential factors of $\epsilon^+(p^+M\otimes E^{-t\tau+\gamma\tau^2/2})$ are the following rational functions ($c\in C$)
\[
\epsilon^*(\gamma\tau^2/2-t\tau-ct^2/2)=\frac{\gamma v^2-2v-c}{2u^2}
\]
whose numerator has only simple zeros if $\gamma\not\in\wh C=-1/C$ or if $\gamma=\wh c_o:=-1/c_o$ and $c\neq c_o$, and has a double zero at $v=-c_o$ if $\gamma=\wh c_o$.
\end{proof}

\begin{corollaire}
For each $\gamma\in\CC$, the complex $\DR^{\rdD_Z}\epsilon^+(p^+M\otimes E^{-t\tau+\gamma\tau^2/2})$ has cohomology in degree zero at most.
\end{corollaire}

\begin{proof}
Away from the zeros of $\gamma v^2-2v-c$ ($c\in C$) regarded as points on $E$, the meromorphic flat bundle $\epsilon^+(p^+M\otimes E^{-t\tau+\gamma\tau^2/2})$ is good, and even very good in the sense of \cite[\S7]{Bibi93}, hence one can apply Majima's theorem \cite{Majima84} (\cf also \cite[\S7]{Bibi93}). On the other hand, at the remaining points one can apply Lemma \ref{lem:DRmodv/u}, according to Lemma \ref{lem:zeros}.
\end{proof}

\begin{definition}[Topological Laplace transformation]\label{def:topLaplace}
Let $(\cF,\cF_\bbullet)$ be a Stokes-filtered sheaf on $\wt\PP^1_t$ of Gaussian type $C$. Define $\cG$ and $\cG_{<\gamma}$ by \eqref{eq:Ggammagen} and \eqref{eq:Ggamma} with the blow-up $\epsilon:Z\to\PP^1_t\times\PP^1_\tau$. Set $q=p\circ\epsilon$ and $\wh q=\wh p\circ\epsilon$. The topological Laplace transform $(\wh\cF,\wh\cF_\bbullet)$ of $(\cF,\cF_\bbullet)$ is the Stokes-filtered sheaf defined by the data
\[
\wh\cF^*:=R^1\wt{\wh q}_*\wtj^{-1}\cG,\quad \wh\cL_{<\gamma}=\cH^1\bR\wt{\wh q}_*\cG_{<\gamma}\quad(\gamma\in\CC)
\]
and the gluing morphism $\FcL_{<\gamma}\to\FcL$, where $\FcL$ is the constant sheaf of rank~$r$ regarded as the restriction to~$S^1_{\wh\infty}$ of the push-forward of the constant sheaf $\FcF^*$, is induced by the push-forward of the natural morphism $\cG_{<\gamma}\to\cG$.
\end{definition}

As a direct consequence of Theorem \ref{th:mochi} together with Proposition \ref{prop:Ggamma} we obtain:

\begin{theoreme}\label{th:mainLaplace}
The definition above produces a Stokes-filtered sheaf of Gaussian type $\wh C$ which is isomorphic to the Stokes-filtered sheaf attached to $\wh M$.
\end{theoreme}

\begin{remarque}
It is not obvious, a~priori, that the objects considered in Definition \ref{def:topLaplace} form a Stokes-filtered sheaf. It could be proved in a purely topological way in the present setting, without referring to the analytic Theorem \ref{th:mochi}.
\end{remarque}

\subsection{Topological Laplace transformation on \texorpdfstring{$\Afu_\tau$}{A1}}\label{subsec:toplaplaceAtau}
In this section we make a little more explicit the topological expression of the restriction $\FcF^*$ to $\Afu_\tau$ of the Laplace transform of~$\wh\cF$. In particular the blowing-up $Z$ will not be used here. We consider the pull-back $\wt p^{-1}(\cF,\cF_\bbullet)$ on the open subset $\wt\PP^1_t\times\Afu_\tau$ of~$\wt\PP^1_t\times\wt\PP^1_\tau$. 

We then have $\FcF^*=R^1\wt{\wh p}_*(\wtj^{-1}\cG)=R^1\wt{\wh p}_*(\wtj^{-1}\wt p^{-1}\cF_{\leq0})$. Therefore, $\FcF^*$ is the constant sheaf with fibre $H^1(\wt\PP^1_t,\cF_{\leq0})$, which has dimension $r$ (rank of $\cF^*$), as follows from the lemma below. Let us check directly the vanishing of all other $R^k\wt{\wh p}_*(\wtj^{-1}\cG)$.

\begin{lemme}\label{lem:H1distancefinie}
For a filtered local system of pure Gaussian type $(\cF,\cF_\bbullet)$, we have $H^j(\wt\PP^1_t,\cF_{\leq0})=0$ for $j\neq1$ and $\dim H^1(\wt\PP^1_t,\cF_{\leq0})=r$.
\end{lemme}

\begin{proof}
We consider the closed covering $(\wt\PP^1_{t,\mu})_{\mu\in\ZZ/4\ZZ}$ of $\wt\PP^1_t$ consisting of points $t$ with $\arg t\in I^{(\mu)}$. In $\CC_{t,\mu}$, $\cF_{\leq0}$ is the constant sheaf. We first prove that $(\wt\PP^1_{t,\mu})_{\mu\in\ZZ/4\ZZ}$ is a Leray covering for $\cF_{\leq0}$ and that $H^0(\wt\PP^1_{t,\mu},\cF_{\leq0})=0$ for each $\mu$.

On each $\wt\PP^1_{t,\mu}$, we have an isomorphism between $(\cF,\cF_\bbullet)$ and $(\gr\cF,\gr_\bbullet\cF)$, so we can assume that $(\cF,\cF_\bbullet)$ is graded, and we are easily reduced to the case where $(\cF,\cF_\bbullet)$ has exactly one exponential factor $c\in C$. Note now that, on $I^{(\mu)}\subset S^1_\infty$, the function $\re^{ct^2/2}$ changes its asymptotic behaviour (from exponential growth to rapid decay) exactly once, and this occurs in the interior of $I^{(\mu)}$. Therefore, when restricted to~$I^{(\mu)}$, $\cF_{\leq0}$ is zero on some nonempty closed interval and is the maximal extension of~$\cF_{|\CC_{t,\mu}}$ on the complementary nonempty open interval. It follows that $\cF_{\leq0}$ is acyclic on $\wt\PP^1_{t,\mu}$.

The two-by-two (or more) intersections of the covering are either segments or the origin, and it is easy to check that either $\cF_{\leq0}$ is acyclic on such a segment, or is the constant sheaf on such a segment, hence has no $H^k$ with $k\geq1$. This ends the proof of our assertion.

Next, we prove that $H^2(\wt\PP^1_t,\cF_{\leq0})=0$. This is obtained by Poincaré-Verdier duality, by checking that the Verdier dual of $\cF_{\leq0}$ is a sheaf up to a shift, which satisfies properties similar to that of $\cF_{\leq0}$ (\cf the proof of Lemma 4.16 in \cite{Bibi10}).

Lastly, the rank of $H^1(\wt\PP^1_t,\cF_{\leq0})$ is obtained by a simple Euler characteristic computation.
\end{proof}

\Subsection{The topological space for computing the Laplace transform on \texorpdfstring{$\wt\PP^1_\tau$}{P1}}\label{subsec:topspaceLaplace}

We now make more explicit the computation of $\cG$ and $\cG_{<\gamma}$ ($\gamma\in\CC$). We denote by $\wt Z$ the real blow-up space of the components of $D$ in $Z$, so that the map $\epsilon$ lifts as $\wt\epsilon:\wt Z\to\wt\PP^1_t\times\wt\PP^1_\tau$. We notice that $\Afu_t\times\wt\PP^1_\tau$ and $\wt\PP^1_t\times\Afu_\tau$ are naturally included in~$\wt Z$ as open subsets, which complement respectively $\partial\wt Z_{|D_\infty\cup E}$ and $\partial\wt Z_{|D_{\wh\infty}\cup E}$. We will mainly consider the complementary inclusions
\begin{equation}\label{eq:wtij}
\partial\wt Z_{|D_{\wh\infty}\cup E}\Hto{\wti}\wt Z\Hfrom{\wtj}\wt\PP^1_t\times\Afu_\tau.
\end{equation}

The restriction $\partial\wt Z_{|E}$ is described in this way:
\begin{align*}
\partial\wt Z_{|E}&\simeq S^1_u\times S^1_{u'}\times[0,\infty]_{u'}\\
&\simeq S^1_{t'}\times S^1_{\tau'}\times[0,\infty]_{u'},\\
&\simeq S^1_u\times S^1_v\times[0,\infty]_v,
\end{align*}
where the isomorphisms on the arguments are given by
\begin{align*}
(\arg u,\arg u')&\mto(\arg t'=\arg u,\arg\tau'=\arg u+\arg u')\\
(\arg u,\arg u')&\mto(\arg u, \arg v=-\arg u')
\end{align*}
and on the absolute values by $|v|=1/|u'|$. The restriction $\partial\wt Z_{|E^*}$ is obtained by replacing $[0,\infty]$ with $(0,\infty)$ in the formulas above. On the other hand, $\partial\wt Z_{|D_{\wh\infty}}$ (\resp $\partial\wt Z_{|D_\infty}$) is identified to the space $\wt\PP^1_t\times S^1_{\tau'}$ (\resp $S^1_{t'}\times\wt\PP^1_\tau$), and the gluing with $\partial\wt Z_{|E}$ is obtained by identifying $\partial\wt\PP^1_t\times S^1_{\tau'}$ with $S^1_{t'}\times S^1_{\tau'}\times\{|u'|=0\}$ (\resp $S^1_{t'}\times\partial\wt\PP^1_\tau$ to $S^1_{t'}\times S^1_{\tau'}\times\{|u'|=\infty\}$). Notice that we will also use the notation $S^1_\infty$ for $\partial\PP^1_t$ and $S^1_{\wh\infty}$ for $\partial\wt\PP^1_\tau$, and it will be clear from the context whether we use $\arg t$ or $\arg t'$ (\resp $\arg\tau$ or $\arg\tau'$) as a coordinate.

We will use the following diagrams:
\begin{equation}\label{eq:diagZZtilde}
\begin{array}{ccc}
\xymatrix{
&Z\ar[d]_-\epsilon&\\
&\PP^1_t\times\PP^1_\tau\ar[dl]_p\ar[dr]^{\wh p}&\\
\PP^1_t&&\PP^1_\tau
}
&\qquad&
\xymatrix{
&\wt Z\ar[d]_-{\wt\epsilon}&\\
& \wt\PP^1_t\times \wt\PP^1_\tau\ar[dl]_{\wt p}\ar[dr]^{\wt{\wh p}}&\\
\wt\PP^1_t&& \wt\PP^1_\tau
}
\end{array}
\end{equation}

\begin{lemme}\label{lem:structwtZ}
The space $\wt Z$ is homeomorphic to a product of two closed discs.
\end{lemme}

\begin{proof}
Firstly, $\wt\epsilon$ induces a diffeomorphism $\wt\epsilon^{-1}(\wt\PP^1_t\times\nobreak\Afu_\tau)\simeq\wt\PP^1_t\times\Afu_\tau$, which is a product of a closed disc by an open disc. We will now identify $\partial\wt Z_{|D_{\wh\infty}\cup E}$ with a product of a closed disc with $S^1_{\tau'}$.

We regard $\partial\wt Z_{|E}$ as the product $S^1_u\times[0,\infty]_{u'}\times S^1_{\tau'}$. Recall that $\partial\wt Z_{|D_{\wh\infty}^*}$ is naturally identified with $\Afu_t\times S^1_{\tau'}$. We now identify $\Afu_t\cup(S^1_u\times[0,\infty]_{u'})$ with the closed disc having coordinate $w$ ($|w|\leq\infty$) in the following way:
\begin{itemize}
\item
for $|w|<1$, we have $w=\dfrac{|t|}{1+|t|}\,\re^{\ri\arg t}$,
\item
for $|w|\geq1$, we have $w=(1+|u'|)\re^{-\ri\arg u}$.
\end{itemize}
We will then set $\partial\wt Z_{|D_{\wh\infty}\cup E}\simeq\ov\Delta_w\times S^1_{\tau'}$, where $\Delta_w\simeq\CC_w$.
\end{proof}

\subsection{Behaviour of the function \texorpdfstring{$\exp\bigl((c/2)t^2+t\tau-(\gamma/2)\tau^2\bigr)$}{exp} on \texorpdfstring{$\partial\wt Z$}{dZ}}
For the sake of simplicity, we will denote by $\partial\wt Z$ the restriction of $\wt Z$ above $\{t=0\}\cup D_{\wh\infty}\cup E$ (so~we do not take into account the restriction over $D_\infty^*$).

Let us fix an nonzero complex number $c\in\CC^*$ and let us consider the pull-back of the function $\exp\bigl((c/2)t^2+t\tau-(\gamma/2)\tau^2\bigr)$ to $Z$ and $\wt Z$, where $\gamma\in\CC^*$ is regarded as a parameter that we vary.

\begin{definition}[of $\ov\Delta_w^\rrd(\wh\theta,c,\gamma)$]
For $\wh\theta\in S^1_{\tau'}$ fixed, we define
\[
\ov\Delta_w^\rrd(\wh\theta,c,\gamma)\subset \ov\Delta_w\times\{\wh\theta\}\subset \ov\Delta_w\times S^1_{\tau'}=\partial\wt Z
\]
as the subset consisting of points in the neighbourhood of which the function $\exp\bigl((c/2)t^2+t\tau-(\gamma/2)\tau^2\bigr)$ has rapid decay.
\end{definition}

\subsubsection*{Determination of $\ov\Delta_w^\rrd(\wh\theta,c,\gamma)$ in $\partial\wt Z_{|D_{\wh\infty}}$}
Since $\gamma\neq0$, the behaviour for $t$ finite and $\tau\to\infty$ of the function $\exp\bigl((c/2)t^2+t\tau-(\gamma/2)\tau^2\bigr)$ is governed by the sign of $\reel(-(\gamma/2)\tau^2)$. We find that
\begin{equation}
\ov\Delta_w^\rrd(\wh\theta,c,\gamma)\cap\partial\wt Z_{|D_{\wh\infty}}=
\begin{cases}
\Delta_w^{\leq1}&\text{if }0\lethetah \gamma,\\
\emptyset&\gamma\leqthetah 0.
\end{cases}
\end{equation}

\subsubsection*{Determination of $\ov\Delta_w^\rrd(\wh\theta,c,\gamma)$ in $\partial\wt Z_{|E}$}
Let us first work with the coordinate $w'=w-1=|u'|\re^{-\ri\arg u}$ with $|u'|\in(0,+\infty]$ (since we already know the behaviour at $|u'|=0$ by the previous computation). Using the expression
\begin{equation}\label{eq:expressionu'}
(c/2)t^2+t\tau-(\gamma/2)\tau^2=\frac{cu^{\prime2}+2u'-\gamma}{2u^2u^{\prime2}}
\end{equation}
in the chart $(u,u')$, which becomes $(cu^{\prime2}+2u'-\gamma)\re^{-2\ri\wh\theta}/2$ when we restrict to $\arg\tau'=\arg u+\arg u'=\wh\theta$, and since $u'=\re^{-\ri\wh\theta}w'$, we have
\[
\ov\Delta_w^\rrd(\wh\theta,c,\gamma)\cap \partial\wt Z_{|E}=\{w'\mid\reel(c\,w^{\prime2}+2\re^{-\ri\wh\theta}w'-\gamma\,\re^{-2\ri\wh\theta})<0\},
\]
that is,
\[
\ov\Delta_w^\rrd(\wh\theta,c,\gamma)\cap \partial\wt Z_{|E}=\{w'\mid\reel\big(c(w'-\re^{-\ri\wh\theta}\wh c)^2-(\gamma-\wh c)\re^{-2\ri\wh\theta}\big)<0\}.
\]
The argument (in the coordinate $w'$) of the center $\wh c(\wh\theta):=\re^{-\ri\wh\theta}\wh c$ of the real hyperbola having equation $\reel\big(c(w'-\re^{-\ri\wh\theta}\wh c)^2-(\gamma_{\wh\theta}-\wh c)\re^{-2\ri\wh\theta}\big)=0$ is equal to $\arg\wh c(\wh\theta)=\pi-\wh\theta-\arg c$.

\subsubsection*{Pictures of $\ov\Delta_w^\rrd(\wh\theta,c,\gamma)$ in $\partial\wt Z$}
\begin{figure}[htb]
\begin{center}
\begin{tabular}{ccccccc}
\includegraphics[scale=.22]{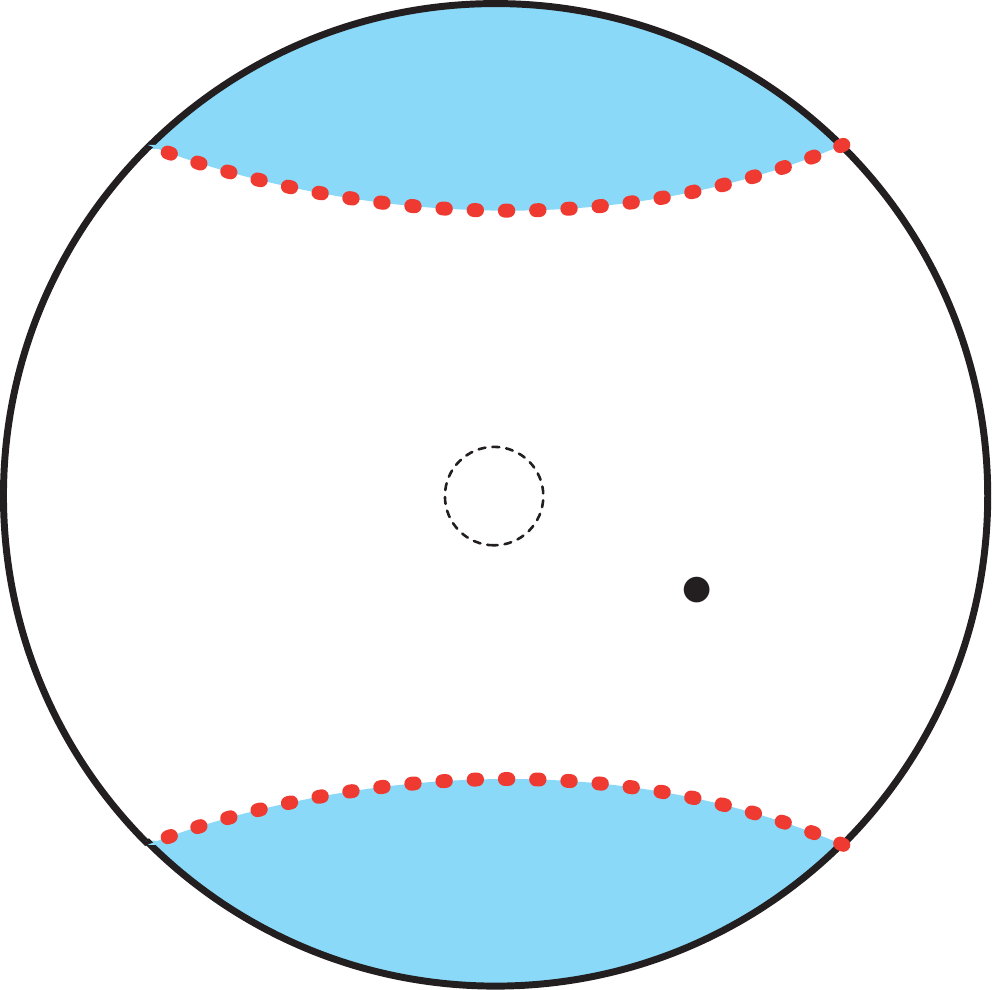}&&\includegraphics[scale=.22]{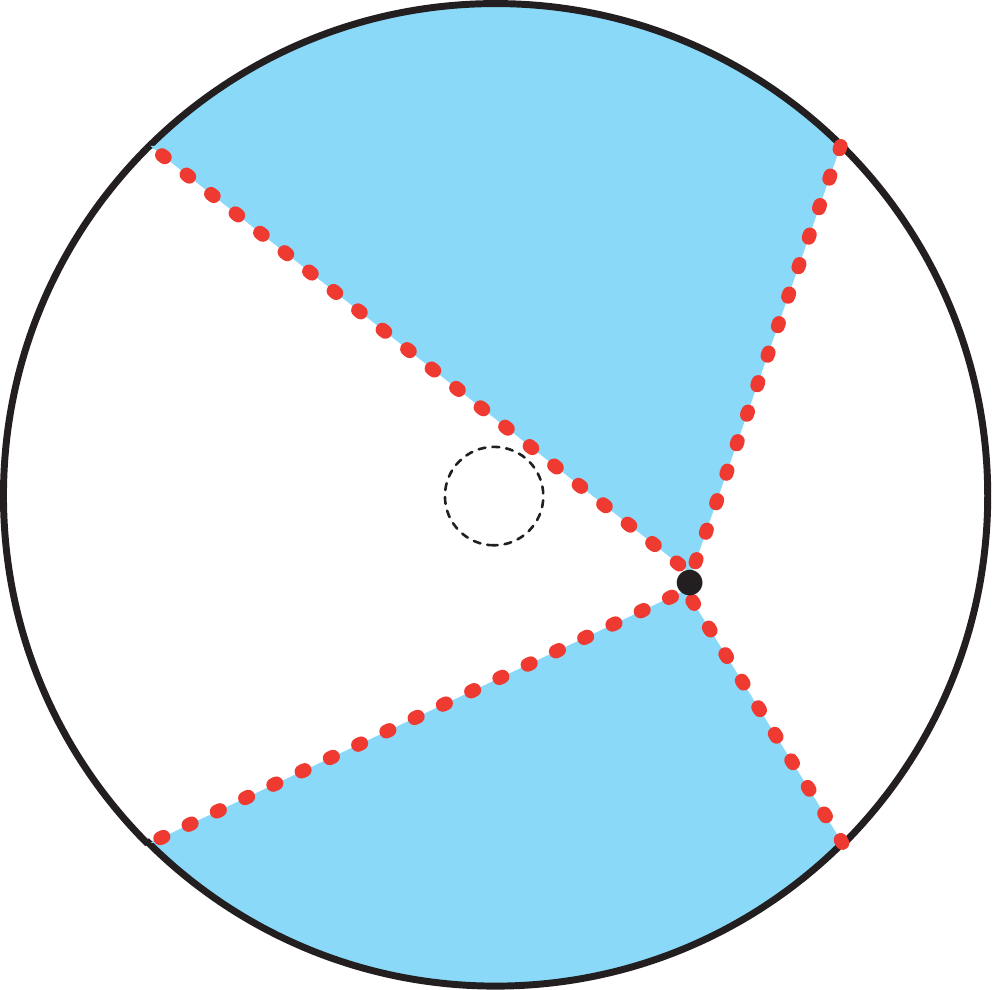}&&\includegraphics[scale=.22]{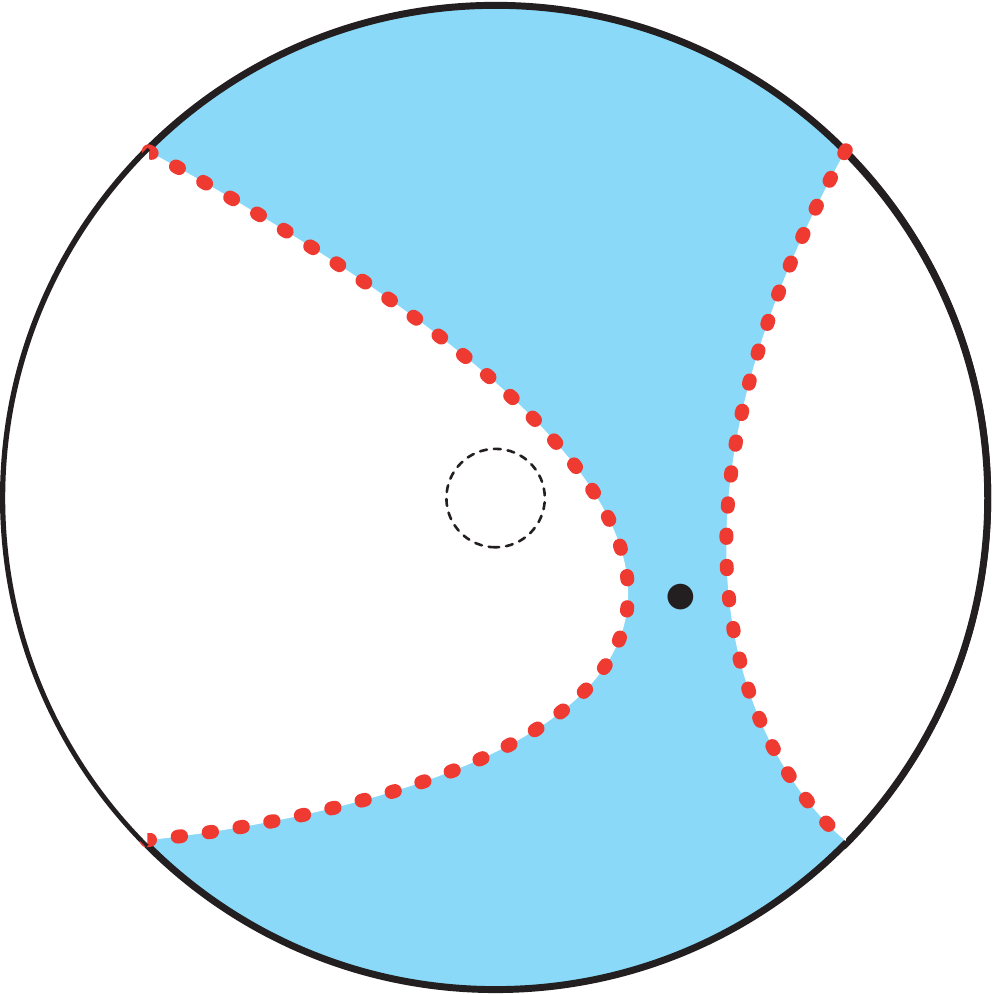}&&\includegraphics[scale=.22]{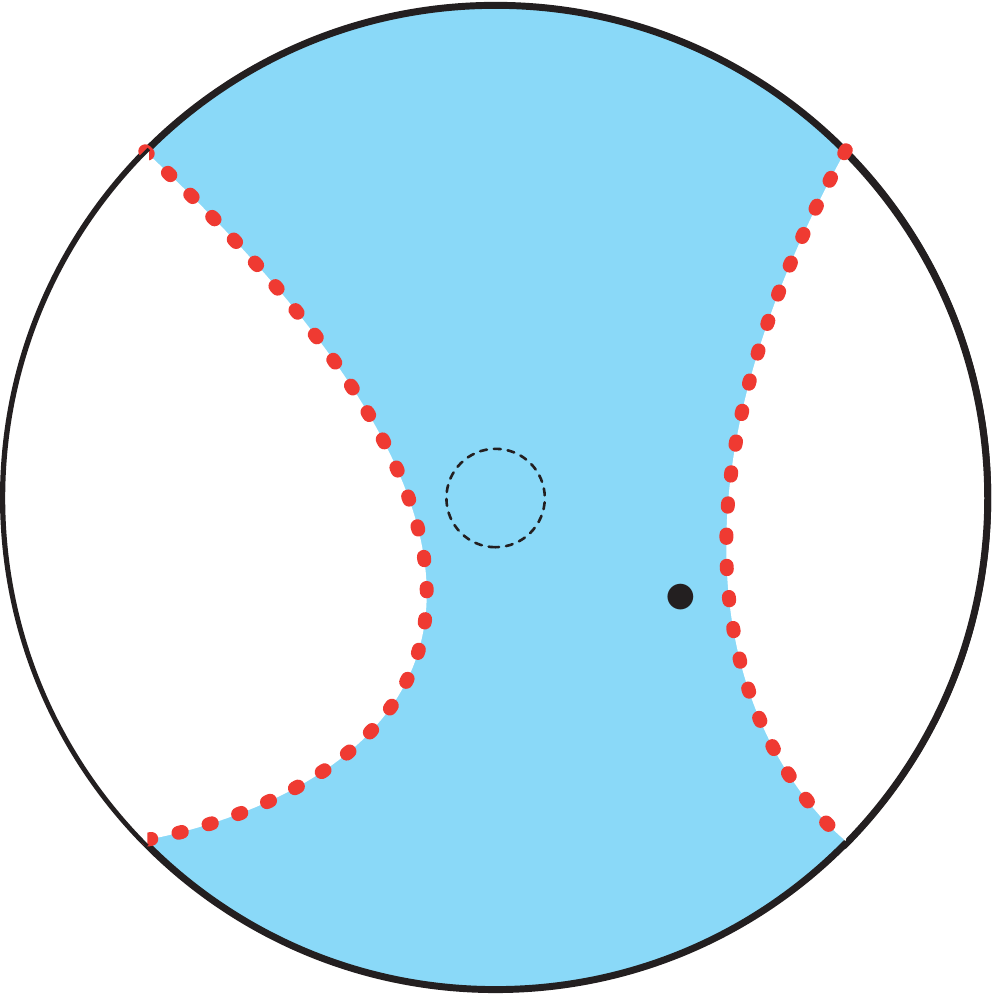}\\
$\gamma\lethetah\wh c\lethetah0$&&$\gamma\neq\wh c$&&$\wh c\lethetah \gamma\lethetah0$&&$\wh c\lethetah0\lethetah \gamma$
\end{tabular}
\caption{$\ov\Delta_w^\rrd(\wh\theta,c,\gamma)$ if $\wh c\lethetah0$}\label{fig:rdnueven}
\end{center}
\end{figure}

\begin{figure}[htb]
\begin{center}
\hspace*{-.8cm}\begin{tabular}{ccccccc}
\includegraphics[scale=.22]{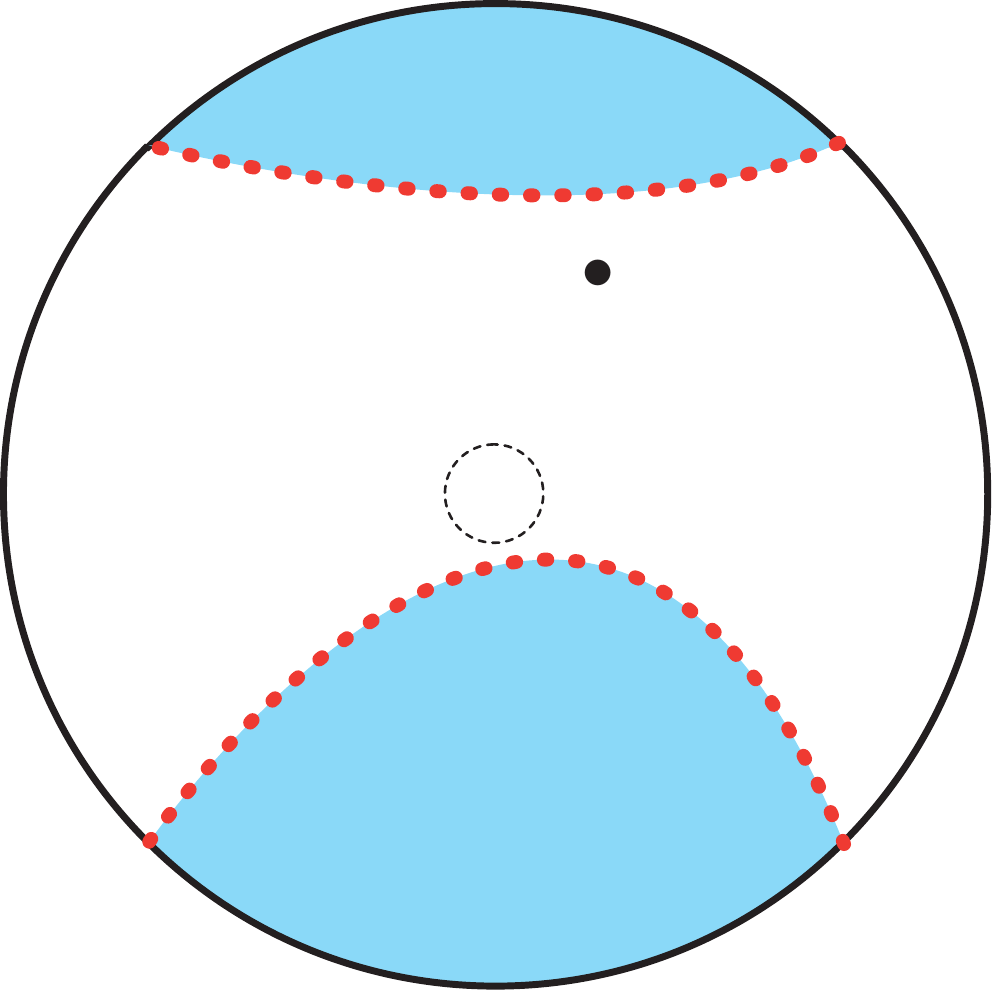}&&\includegraphics[scale=.22]{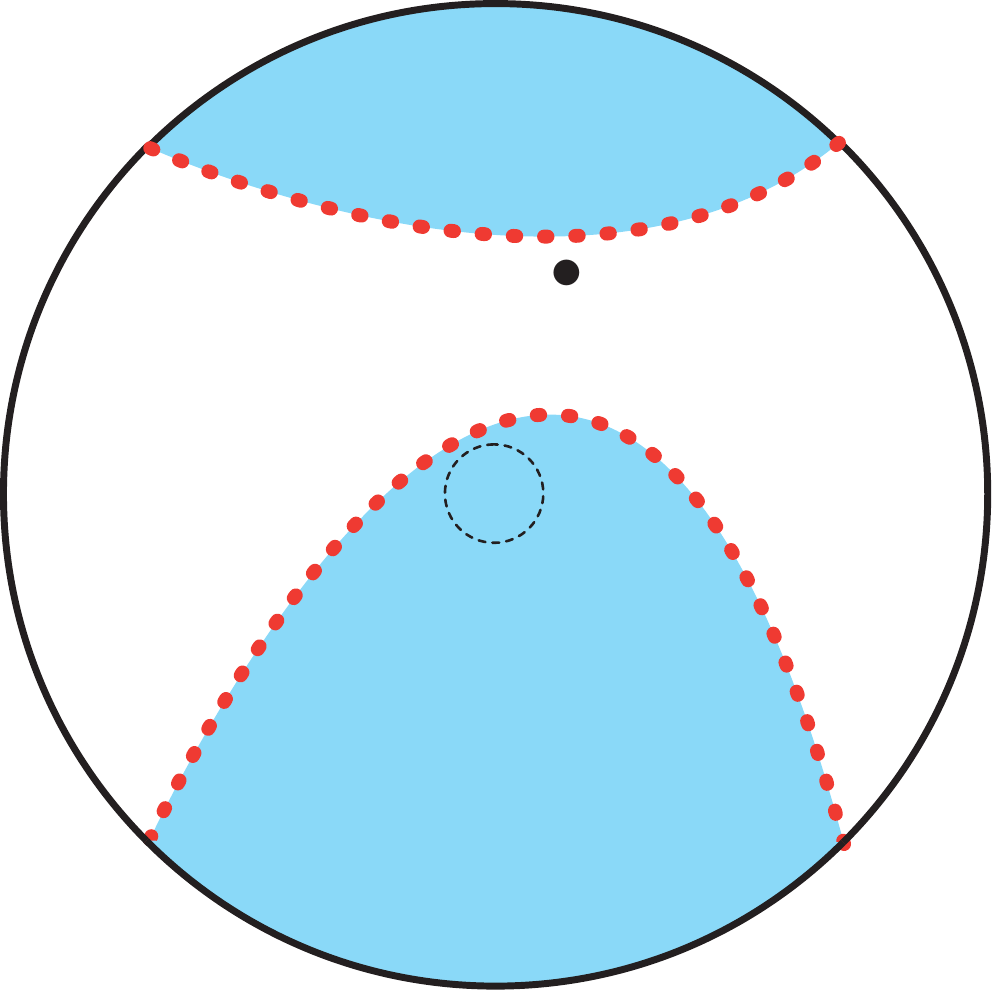}&&\includegraphics[scale=.22]{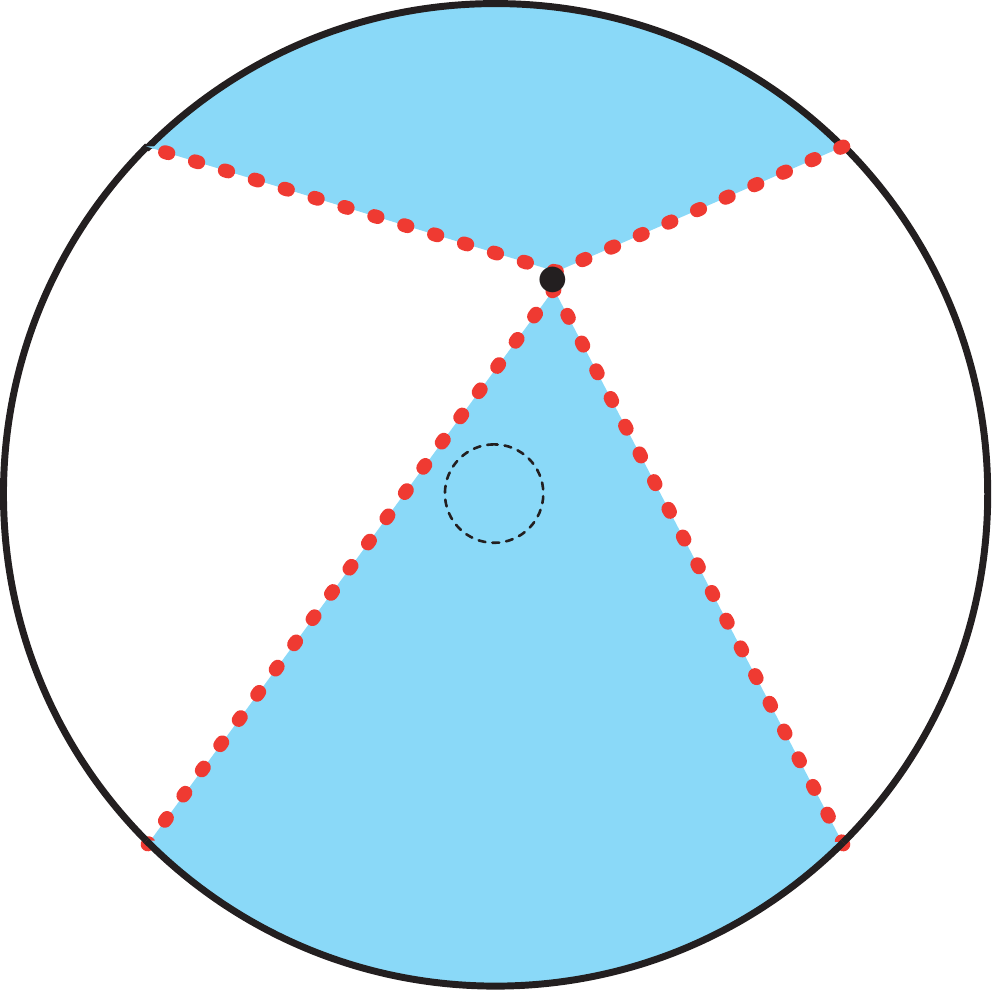}&&\includegraphics[scale=.22]{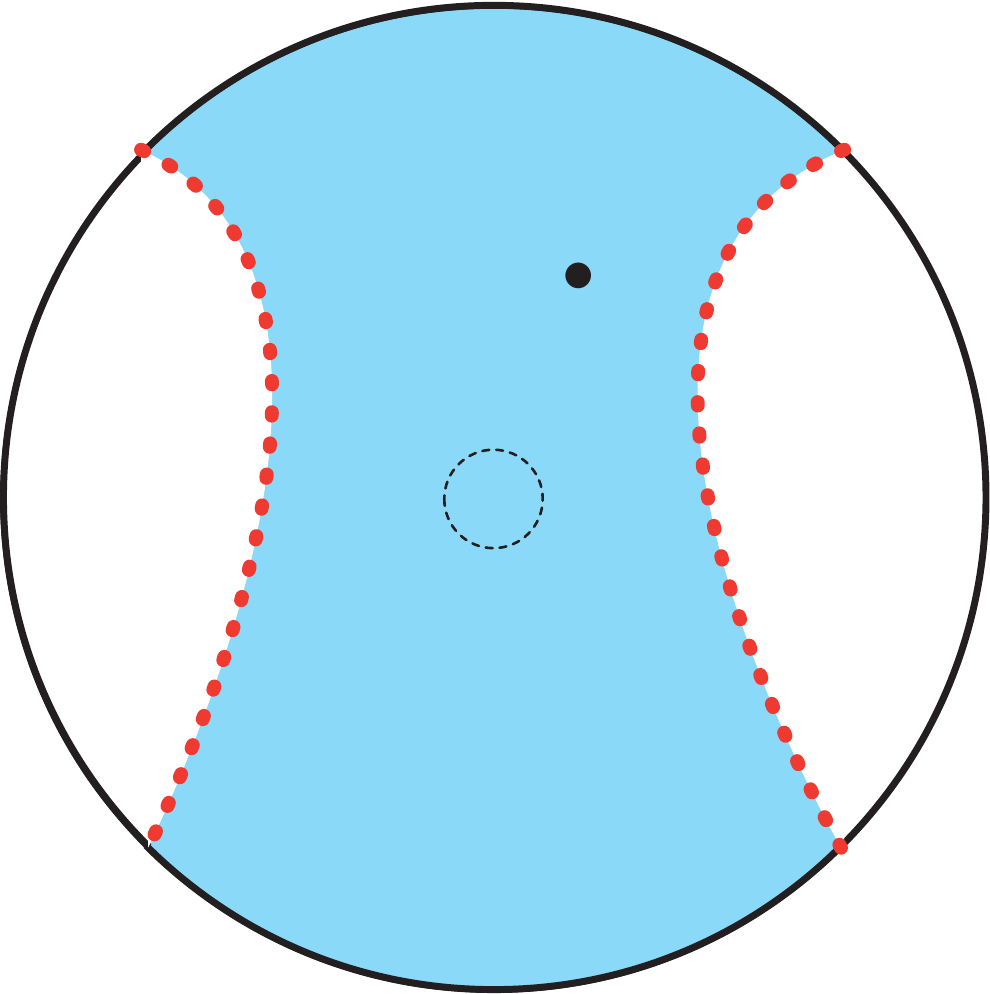}\\
$\gamma\lethetah0\lethetah \wh c$&&$0\lethetah \gamma\lethetah\wh c$&&$\gamma=\wh c$&&$0\lethetah\wh c\lethetah \gamma$
\end{tabular}
\caption{$\ov\Delta_w^\rrd(\wh\theta,c,\gamma)$ if $0\lethetah\wh c$}\label{fig:rdnuodd}
\end{center}
\end{figure}

For these pictures, we assume that $c$ is real and positive. On the pictures \ref{fig:rdnueven} and \ref{fig:rdnuodd}, the black dot is the center $\wh c(\wh\theta)$, the middle circle is the circle $|w|=1$ and the external circle is the circle $|w|=\infty$. The regions $\ov\Delta_w^\rrd(\wh\theta,c,\gamma)$ are the colored regions and the boundaries are excluded if they are red. The disc $\Delta_w^{\leq1}$ is also indicated on the picture.

It will be useful to remark the following:

\pagebreak[2]
\begin{lemme}\label{lem:positionrelative}\mbox{}
\begin{enumerate}
\item\label{lem:positionrelative1}
For any $c$ and $\wh\theta$ fixed, we have
\[
\gamma'\lethetah \gamma\implique \ov\Delta_w^\rrd(\wh\theta,c,\gamma')\subset \ov\Delta_w^\rrd(\wh\theta,c,\gamma).
\]
\item\label{lem:positionrelative2}
We have $\wh c(\wh\theta)\in \ov\Delta_w^\rrd(\wh\theta,c,\gamma)$  if and only if $\wh c\lethetah \gamma$.
\item\label{lem:positionrelative3}
Let us fix two positive numbers $c<c'$.
\begin{itemize}
\item
Assume $\cos2\wh\theta>0$. Then $\wh c'(\wh\theta)\in \ov\Delta_w^\rd(\wh\theta,c',\gamma)\implique \wh c'(\wh\theta)\in \ov\Delta_w^\rd(\wh\theta,c,\gamma)$.
\item
Assume $\cos2\wh\theta<0$. Then $\wh c(\wh\theta)\in \ov\Delta_w^\rd(\wh\theta,c,\gamma)\implique \wh c(\wh\theta)\in \ov\Delta_w^\rd(\wh\theta,c',\gamma)$.\qed
\end{itemize}
\end{enumerate}
\end{lemme}

\subsection{The sheaves \texorpdfstring{$\cG$}{G} and \texorpdfstring{$\cG_{<\gamma}$}{Gg} (\texorpdfstring{$\gamma\in\CC$}{gC})}\label{subsec:Gleqb}
We now make more explicit the expression of the sheaves $\cG$ and $\cG_{<\gamma}$. Recall (\cf Lemma \ref{lem:structwtZ} and its proof) that we identify $\partial\wt Z_{|D_{\wh\infty}\cup E}$ with the product $\ov\Delta_w\times S^1_{\tau'}$.

Let $\rho:[0,\infty]\to[0,1]$ be a decreasing homeomorphism such that $\rho(0)=0$ and $\rho(\infty)=1$. It induces a homeomorphism, still denoted by $\rho:w\mto \rho(|w|)\re^{\ri\arg w}$, from~$\ov\Delta_w$ onto $\Delta_w^{\leq1}$, which sends $\Delta_w$ homeomorphically onto $\Delta_w^{<1}=\Afu_t$ and $\{|w|=\infty\}$ onto~$\partial\wt\PP^1_t$.

\begin{enumerate}
\item
The sheaf $\wcF_{\leq0}$ on $\ov\Delta_w$ is the pull-back of the sheaf $\cF_{\leq0}$ by the map \hbox{$\rho:\ov\Delta_w\to\wt\PP^1_t$}. It is a subsheaf of the constant sheaf $\wcF$, pull-back of $\cF$.
\item
The sheaf $\cG$ on $\partial\wt Z_{|D_{\wh\infty}\cup E}\simeq\ov\Delta_w\times S^1_{\tau'}$ is the pull-back of $\wcF_{\leq0}$ by the projection \hbox{$\ov\Delta_w\times S^1_{\tau'}\to\ov\Delta_w$}. It is a subsheaf of the constant sheaf $\cG'$, pull-back of $\wcF$.
\end{enumerate}

Note that the quotient sheaf $\cG'/\cG$ is supported on $\{|w|=\infty\}\times S^1_{\tau'}$ and $\cG$ is constant in the interior of $\ov\Delta_w\times S^1_\tau$.

\begin{lemme}
The push-forward of $\cG$ by the projection $\ov\Delta_w\times S^1_{\tau'}\to S^1_{\tau'}$ is the constant local system of rank $r$  on~$S^1_{\tau'}$.
\end{lemme}

\begin{proof}
Identical to that of Lemma \ref{lem:H1distancefinie}, since there is no topological difference between~$\cF_{\leq0}$ on $\wt\PP^1_t$ and $\wcF_{\leq0}$ on $\ov\Delta_w$.
\end{proof}

\subsubsection*{The sheaves $\cG_{<\gamma}$ on $\partial\wt Z_{|D_{\wh\infty}}\simeq\{|w|\leq1\}\times S^1_{\tau'}$}
Let $j_{0<\gamma}$ denote the open inclusion
\[
(\partial\wt Z_{|D_{\wh\infty}})_{0<\gamma}:=\wt\PP^1_t\times\{\arg \gamma-2\arg\tau'\in(-\pi/2,\pi/2)\bmod2\pi\}\hto\wt\PP^1_t\times S^1_{\tau'}.
\]
(Above $\Afu_t$, this is the domain of rapid decay, equivalently of moderate growth, of the function $\exp\bigl((c/2)t^2+t\tau-(\gamma/2)\tau^2\bigr)$ for $t$ remaining at finite distance). As in Remark \ref{rem:Stokesdir}\eqref{rem:Stokesdir3}, we denote by $\beta_{0<\gamma}$ the functor $j_{0<\gamma,!}j_{0<\gamma}^{-1}$.

On $\partial\wt Z_{|D_{\wh\infty}}$, the first item of \eqref{eq:Ggamma} amounts to (recall that $\cG=\cG'$ is the constant sheaf of rank~$r$ on $\partial\wt Z_{|D_{\wh\infty}}$):
\begin{align*}
\cG_{<\gamma}&:=\beta_{0<\gamma}\,\cG.
\end{align*}
We have a natural inclusion
\[
\cG_{<\gamma}\hto\cG.
\]

\subsubsection*{The sheaves $\cG_{<\gamma}$ on $\partial\wt Z_{|E}$}
In order to understand the computation below, we will regard for a moment $\partial\wt Z_{|E}$ as the product $\ov{E^*}\times S^1_u$ with the identification $S^1_u=S^1_{t'}$ and $\ov{E^*}:=S^1_v\times[0,\infty]_v$ (while $E^*=\CC^*_v=S^1_v\times(0,\infty)_v$). In these coordinates the limit of the expression $(c/2)t^2+t\tau-(\gamma/2)\tau^2$ is $(c+2v-bv^2)\re^{-2\ri\arg u}/2$ (this expression holds on $E^*\times S^1_u$, and one should replace $v$ with $(|v|,\re^{\ri\arg v})$ to extend it to $\ov{E^*}\times S^1_u$).
Above~$E^*$, we regard $\cG_{<\gamma}$ as the family, parametrized by $v\in E^*$, of the sheaves $\cL_{<c_\gamma(v)}$ with
\begin{equation}\label{eq:cbv}
c_\gamma(v)=bv^2-2v.
\end{equation}
For $v=v_o$ fixed, these sheaves have been analyzed in \S\ref{subsec:LleqStokes}. Their definition needs the use of the functor $\beta_{c<c_\gamma(v_o)}$. For $v$ varying, we will similarly consider the open subset $(\partial\wt Z_{|E})_{c<c_\gamma(v)}$ of $\partial\wt Z_{|E}$. The corresponding open inclusion will be denoted by $j_{c<c_\gamma(v)}$, with associated functor $\beta_{c<c_\gamma(v)}:=j_{c<c_\gamma(v),!}\,j_{c<c_\gamma(v)}^{-1}$. The definition of the sheaves $\cL_{<c_\gamma(v_o)}$ from the Stokes data also uses the given covering $S^1_u=\bigcup_{\mu\in\ZZ/4\ZZ}I^{(\mu)}(\theta_o)$, for some choice of $\theta_o$ generic with respect to $C$ (\cf\S\ref{subsec:Stokesdata}). We will therefore have to use the corresponding covering
\[
\partial\wt Z_{|E}=\bigcup_{\mu\in\ZZ/4\ZZ}\ov{E^*}\times I^{(\mu)}(\theta_o)=:\bigcup_{\mu\in\ZZ/4\ZZ}(\partial\wt Z_{|E})^{(\mu)}.
\]

However, we will present $\cG_{<\gamma}$ on $\partial\wt Z_{|E}$ by regarding now $\partial\wt Z_{|E}$ as the product \hbox{$\{1\leq|w|\leq\infty\}\times S^1_{\tau'}$} which will be better suited to the computation of the push-forward by the projection map to $S^1_{\tau'}$. We will work with the coordinates $(|u'|,\re^{-\ri\arg u},\re^{\ri\arg\tau'})$ and we set $w'=|u'|\re^{-\ri\arg u}$ if $|u'|\neq0,\infty$. The set $(\partial\wt Z_{|E})_{c<c_\gamma(v)}$ is then defined as
\begin{equation}\label{eq:Zc<c_gamma(v)}
(\partial\wt Z_{|E})_{c<c_\gamma(v)}:=\{(w',\re^{\ri\arg\tau'})\mid\reel(c\,w^{\prime2}+2w'\re^{-\ri\arg\tau'}-\gamma\,\re^{-2\ri\arg\tau'})<0\}.
\end{equation}
Strictly speaking, this definition holds on $|w'|\in(0,\infty)$, and we implicitly understand that $(\partial\wt Z_{|E})_{c<c_\gamma(v)}$ is the closure of the subset defined by this relation, that is,
\begin{enumerate}\renewcommand{\theenumi}{\ref{eq:Zc<c_gamma(v)}.\alph{enumi}}
\item\label{enum:a}
the restriction of $(\partial\wt Z_{|E})_{c<c_\gamma(v)}$ to $|w'|=\infty$ is defined by $\reel(c\,\re^{-2\ri\arg u})<\nobreak0$, so is identified with the product $(S^1_u)_{c<0}\times S^1_{\tau'}$,
\item\label{enum:b}
the restriction of $(\partial\wt Z_{|E})_{c<c_\gamma(v)}$ to $|w'|\!=\!0$ is defined by \hbox{$\reel(-\gamma\,\re^{-2\ri\arg\tau'})\!<\!0$}, \ie coincides there with the restriction of the subset $(\partial\wt Z_{D_{\wh\infty}})_{0<\gamma}$ considered above.
\end{enumerate}
We will also use the closed covering $\partial\wt Z_{|E}=\bigcup_{\mu\in\ZZ/4\ZZ}(\partial\wt Z_{|E})^{(\mu)}$, where $(\partial\wt Z_{|E})^{(\mu)}$ is defined by $\re^{-\ri\arg w'}:=\re^{\ri\arg u}\in I^{(\mu)}(\theta_o)$.

\begin{lemme}\label{lem:cG'}
On $\partial\wt Z_{|E}$, the constant sheaf $\cG'$ can be obtained by gluing the sheaves
\[
\cG^{\prime(\mu)}=\bigoplus_{c\in C}\CC_{(\partial\wt Z_{|E})^{(\mu)}}\otimes_\CC G_c^{(\mu)}
\]
with gluing morphisms $g^{(\mu,\mu-1)}=\id_{\CC_{(\partial\wt Z_{|E})^{(\mu)}\cap(\partial\wt Z_{|E})^{(\mu-1)}}}\otimes S^{(\mu,\mu-1)}$.
\end{lemme}

\begin{proof}
This follows from the similar property for $\cL$ on $S^1_\infty$, which is equivalent to the property that the product of the Stokes matrices is equal to the identity.
\end{proof}

\begin{lemme}\label{lem:cGgamma}
On $\partial\wt Z_{|E}$ the subsheaves $\cG_{<\gamma}\subset\cG'$ are obtained by gluing the subsheaves
\[
\cG^{(\mu)}_{<\gamma}=\bigoplus_{c\in C}(\beta_{c< c_\gamma(v)}\CC_{(\partial\wt Z_{|E})^{(\mu)}})\otimes_\CC G_c^{(\mu)}=:\bigoplus_{c\in C}\cG^{(\mu)}_{<\gamma,c}
\]
with the isomorphisms induced by $g^{(\mu,\mu-1)}$.
\end{lemme}

\begin{proof}
This is straightforward from the pointwise definition of $\cG_{<\gamma}$ given by the second item of \eqref{eq:Ggamma}. One can check the preservation of the subsheaves by the gluing isomorphisms by using the same argument as that given in \S\ref{subsec:LleqStokes}.
\end{proof}

\begin{lemme}
For each $\gamma\in\CC^*$, $\cG_{<\gamma}$ is a subsheaf of $\cG$, and coincides with the previously defined $\cG_{<\gamma}$ on $\{|w|=1\}\times S^1_{\tau'}$.
\end{lemme}

\begin{proof}
This is immediate from Properties \eqref{enum:a} and \eqref{enum:b} above.
\end{proof}

\begin{remarque}[Restriction of $\cG_{<\gamma}$ to a fiber of $\wh q$]\label{rem:restrGgammafibre}
We will denote (with some abuse) by $\cG_{<\gamma,\wh\theta}$ etc.\ the restriction of $\cG_{<\gamma}$ etc.\ to the fiber $\wh q^{-1}(\wh\theta)$, where~$\wh q$ is the projection $\partial\wt Z_{|D_{\wh\infty}\cup E}\to S^1_{\tau'}$. It follows then from Lemma \ref{lem:positionrelative}\eqref{lem:positionrelative1} that for $\gamma'\lethetah \gamma$ we have natural morphisms $\cG_{<\gamma',\wh\theta}\to\cG_{<\gamma,\wh\theta}\to\cG_{\wh\theta}$.

If we fix $\wh\theta$ in $S^1_{\tau'}$, the restriction $\cG_{<\gamma,\wh\theta}$ to the fibre $[1,\infty]_w\times S^1_w\times\{\wh\theta\}=[0,\infty]_{w'}\times\nobreak S^1_{w'}\times\nobreak\{\wh\theta\}$ is described as follows:
\begin{enumerate}
\item
We first consider the subsets \eqref{eq:Zc<c_gamma(v)}$_{\wh\theta}$ for $\arg\tau'=\wh\theta$ and $c$ varying in $C$ (completed with the corresponding \eqref{enum:a}$_{\wh\theta}$ and \eqref{enum:b}$_{\wh\theta}$). They look like those on Figures \ref{fig:rdnueven} and~\ref{fig:rdnuodd}.
\item
We may use the simplified version of the Stokes data of $(\cL,\cL_\bbullet)$ where $S^{(\mu,\mu-1)}_{c,c}=\id$ for all $c$ and $\mu=1,2,3$, so that $G_c^{(\mu)}=G_c$ for all $c$. On the subset \eqref{eq:Zc<c_gamma(v)}$_{\wh\theta}$ indexed by $c$ we consider the constant sheaf with fibre $G_c$, extended by zero.
\item
We use the covering $((\partial\wt Z_{|E})^{(\mu)})_{\mu\in\ZZ/4\ZZ}$ and the gluing morphisms $g^{(\mu,\mu-1)}$ to replace the direct sum of the previous sheaves to a new sheaf $\cG_{<\gamma,\wh\theta}$.
\end{enumerate}
\end{remarque}

\section{Computation of topological Laplace transforms}\label{sec:computationLaplace}

Our aim in this section is to express the Stokes data attached to the topological Laplace transform $(\FcF,\FcF_\bbullet)$ in terms of those attached to $(\cF,\cF_\bbullet)$. According to Theorem \ref{th:mainLaplace}, this is equivalent to the computation of the Stokes data attached to $\wh M$ in terms of those attached to $M$.

We start with a Stokes-filtered sheaf $(\cF,\cF_\bbullet)$ of type $C\subset\CC^*$. We will make the computation with the following simplifying assumption: $\arg c$ is independent of $c\in C$. We will denote by $\arg C$ the corresponding value. Note that $\wh C=-1/C$ satisfies then the same property and we have $\arg\wh C=\pi-\arg C$.

\begin{remarque}
Corollary \ref{cor:equivCCprime} can be used to reduce the computation to the case where this assumption is fulfilled, in a non explicit way however, hence this is of no use for an explicit computation of Stokes data.

Indeed, if $C$ does not satisfy the previous assumption, we can find $C'\subset\CC^*$ defining a point in $(\CC^*)^n\moins\text{diagonals}$ with $\arg c'$ constant for $c'\in C'$, and a simply connected open subset $X$ of $(\CC^*)^n\moins\text{diagonals}$ containing $C$ and $C'$: choose nonintersecting paths from each $c\in C$ to distinct points of $\RR_+^*$, defining thus a path from $C$ to $C'$ in $(\CC^*)^n\moins\text{diagonals}$, and take for $X$ a simply connected open neighbourhood of this path in $(\CC^*)^n\moins\text{diagonals}$.

Then the equivalence of Corollary \ref{cor:equivCCprime} Laplace-transforms into an equivalence with respect to $\wh X=-1/X$ and if we know the transformation rule of Stokes data for the pair $(C',\wh C')$, we can use this equivalence to obtain the transformation rule for $(C,\wh C)$. Unfortunately, the equivalence of Corollary \ref{cor:equivCCprime} given by $X$ is not explicit in terms of Stokes data.
\end{remarque}

We will express the Stokes data (as in Definition \ref{def:catStokesdatabis}) of type $(\wh C,\wh\theta_o)$ attached to $(\FcF,\FcF_\bbullet)$ in terms of those of type $(C,\theta_o)$ attached to $(\cF,\cF_\bbullet)$. A suitable choice of~$\theta_o$ and $\wh\theta_o$ will simplify the computation.

Let us fix a choice of $\frac12\arg C$. Then the Stokes directions attached to $C$ are $\frac12\arg C+\nobreak k\pi/4\bmod2\pi$, with $k=1,3,5,7$. We can therefore choose $\theta_o=\frac12\arg C$. The order on $C$ at $\theta_o$ is then the order $c\lethetao c'\iff|c|<|c'|$. The numbering $c_1,\dots,c_n$ is by increasing absolute values. Recall that we set $\theta_o^{(\nu)}=\theta_o+\nu\pi/2$, so that the order on $C$ at $\theta_o^{(\nu)}$ is the usual order between the absolute values if $\nu$ is even, and the reverse order if $\nu$ is odd.

We also choose
\[
\wh\theta_o=\pi/2+\frac12\arg\wh C=\pi-\frac12\arg C=\pi-\theta_o
\]
and we set $\wh\theta_o^{(\nu)}=\wh\theta_o+\nu\pi/2$. The order on $\wh C$ at $\wh\theta_o^{(\nu)}$ is by the absolute values of the elements of $\wh C$ if $\nu$ is \emph{odd}, and the reverse order if $\nu$ is \emph{even}. Note that the numbering $\wh c_1,\dots,\wh c_n$ of $\wh C$ at $\wh\theta_o$ is thus equal to the numbering induced by that of~$C$ at $\theta_o$.

In the coordinate $w'$, the centers $\wh c(\wh\theta)$ of the hyperbolas corresponding to a fibre at~$\wh\theta$ are defined by $w'=\wh c\,\re^{-\ri\wh\theta}$, hence satisfy
\[
\arg_{w'}\wh c(\wh\theta)=\pi-\arg c-\wh\theta=\pi-2\theta_o-\wh\theta.
\]
It follows that, for each $\nu\in\ZZ/4\ZZ$, we have
\[
\arg_{w'}\wh c(\wh\theta_o^{(\nu)})=-\theta_o^{(\nu)}
\]
(recall that $\arg_{w'}=-\arg_{t'}$ and we use $\arg_{t'}$ to parametrize $S^1_{t'}=\partial\ov\Delta_w$). Our aim is to compute the filtrations $\wh L_{\leqnu\cbbullet}$ at the points $\wh\theta_o^{(\nu)}$ for each $\nu\in\ZZ/4\ZZ$.

\begin{theoreme}\label{th:main}
Under the previous assumptions, let $(L,L_{\leqnu\cbbullet})$ be the Stokes data of pure Gaussian type $(C,\theta_o)$ attached to $(\cF,\cF_\bbullet)$. Then the Stokes data $(\wh L,\wh L_{\leqnu\cbbullet})$ of pure Gaussian type $(\wh C,\wh\theta_o)$ attached to $(\FcF,\FcF_\bbullet)$ are equal to $(L,L_{\leqnu\cbbullet})$.
\end{theoreme}

We will now fix $\nu\in\ZZ/4\ZZ$ and $\gamma\in\CC$, and we set $k\in\{1,\dots,n\}$ in such a way that
\[
\begin{cases}
\wh c_k\lenu \gamma\leqnu\wh c_{k+1}&(\nu\text{ even}),\\
\wh c_{k+1}\lenu \gamma\leqnu\wh c_k&(\nu\text{ odd}).
\end{cases}
\]
We first describe a Leray covering of $\ov\Delta_w\times\{\wh\theta_o^{(\nu)}\}$ suited for such a computation.

\begin{figure}[htb]
\begin{center}\tabcolsep0pt
\begin{tabular}{cccc}
\includegraphics[scale=.4]{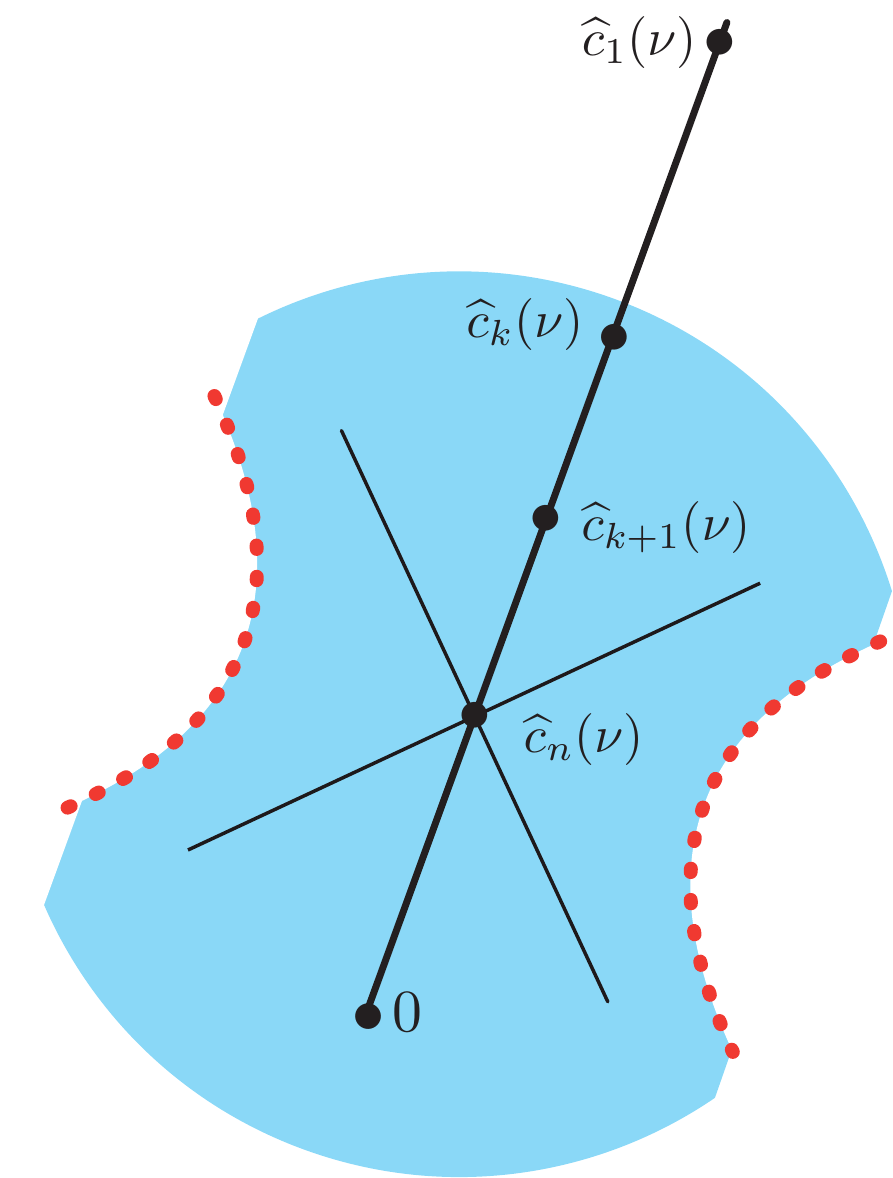}&\includegraphics[scale=.4]{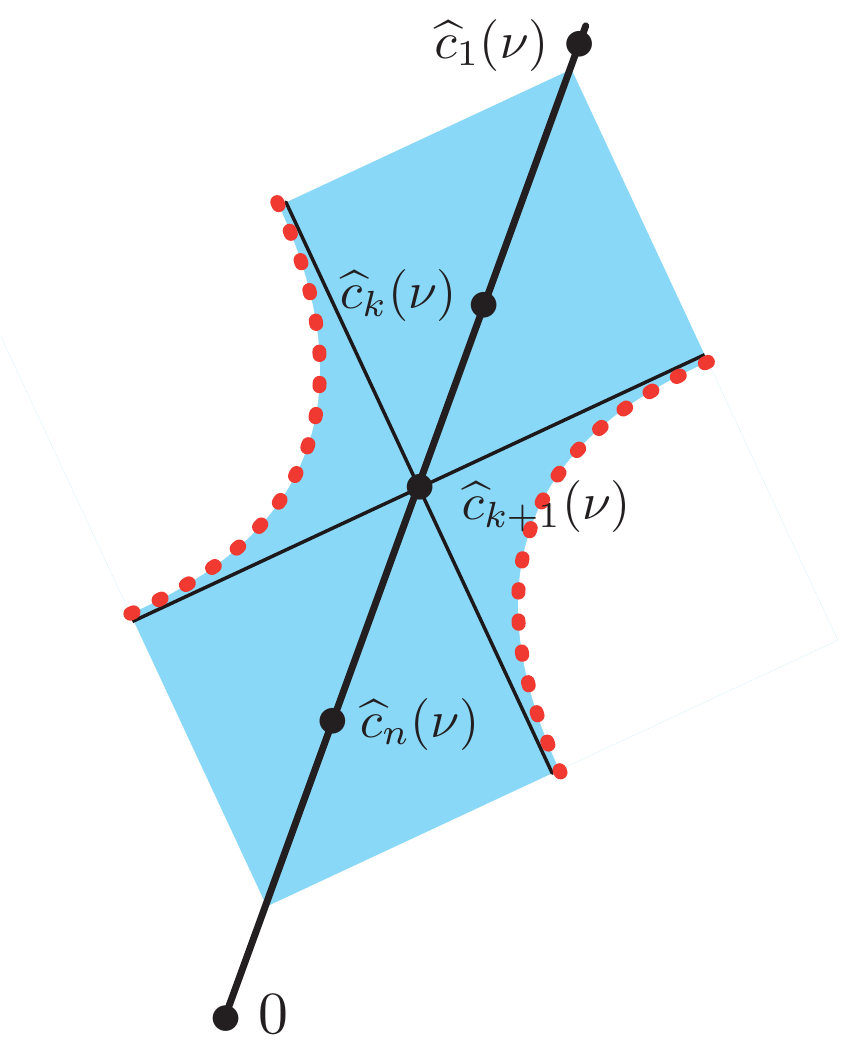}&\includegraphics[scale=.4]{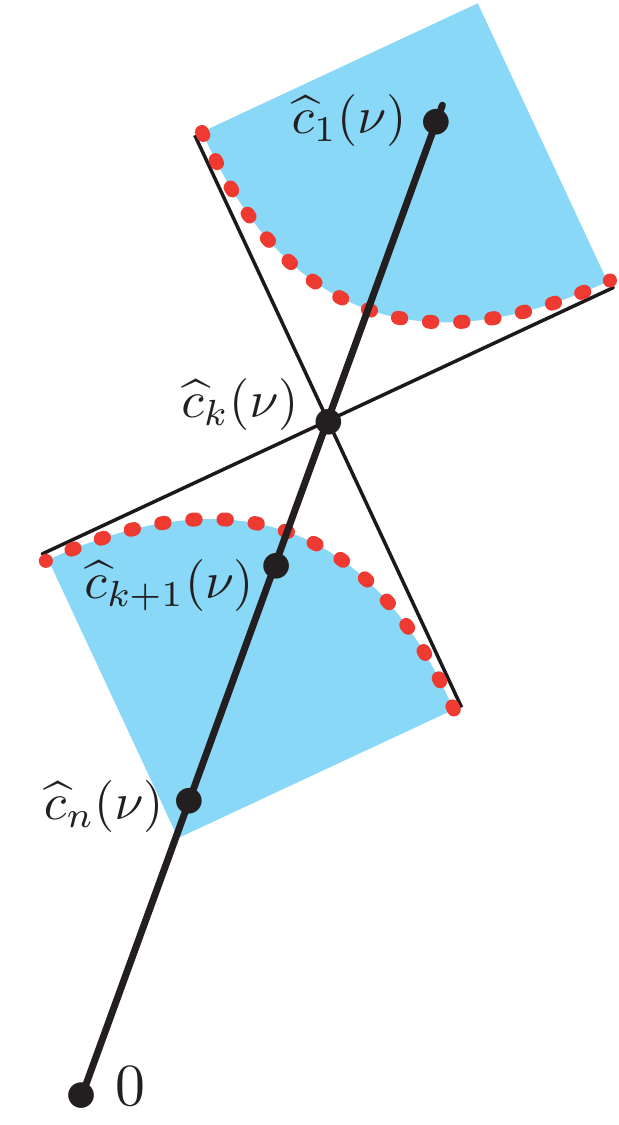}&
\includegraphics[scale=.4]{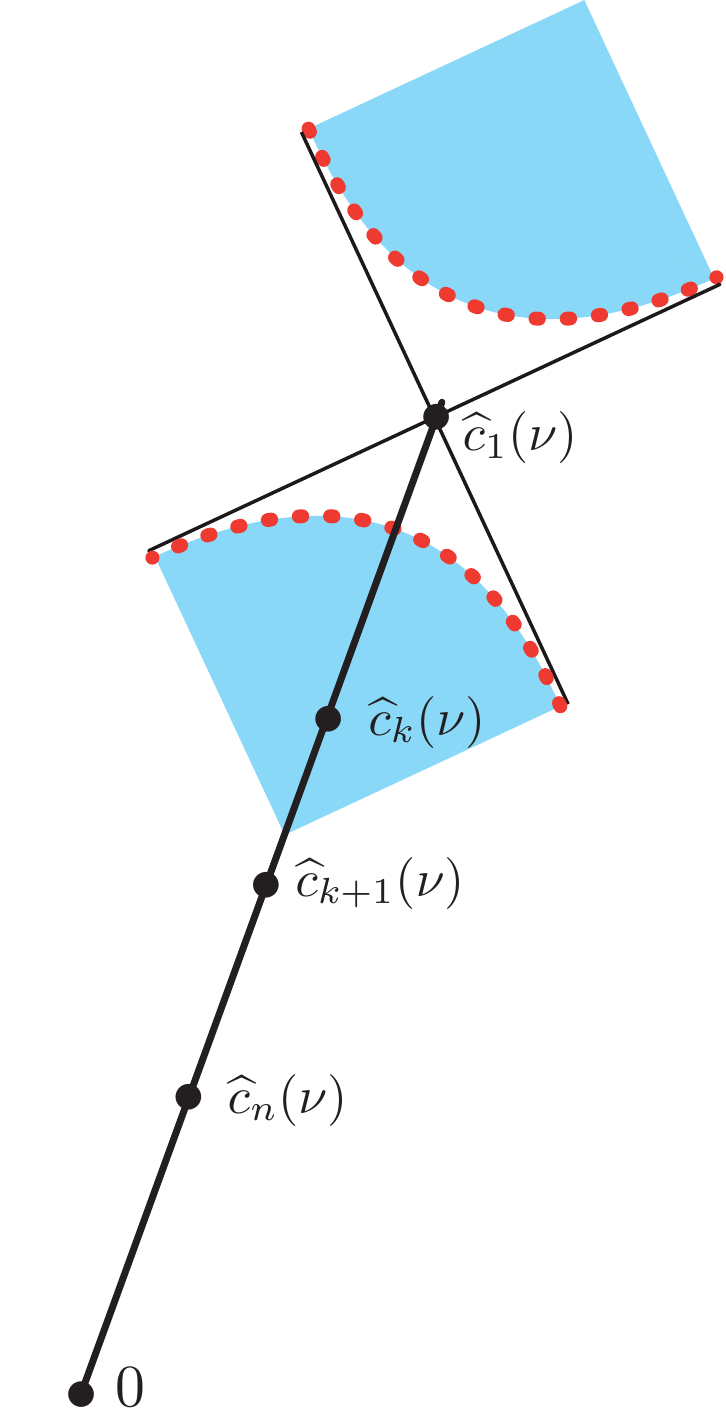}\\
\smaller$\ov\Delta_w^\rrd(\wh\theta_o^{(\nu)},c_n,\gamma)$&
\smaller$\ov\Delta_w^\rrd(\wh\theta_o^{(\nu)},c_{k+1},\gamma)$&\smaller$\ov\Delta_w^\rrd(\wh\theta_o^{(\nu)},c_k,\gamma)$&
\smaller$\ov\Delta_w^\rrd(\wh\theta_o^{(\nu)},c_1,\gamma)$
\end{tabular}
\caption{$\nu$ odd, \eg $\nu=3$, $\theta_o\in(0,\pi/8)$, $\cos2\wh\theta_o^{(1)}<0$}\label{fig:Awmodnuoddcb}
\end{center}
\end{figure}

\begin{figure}[htb]
\begin{center}
\begin{tabular}{ccc}
\includegraphics[scale=.4]{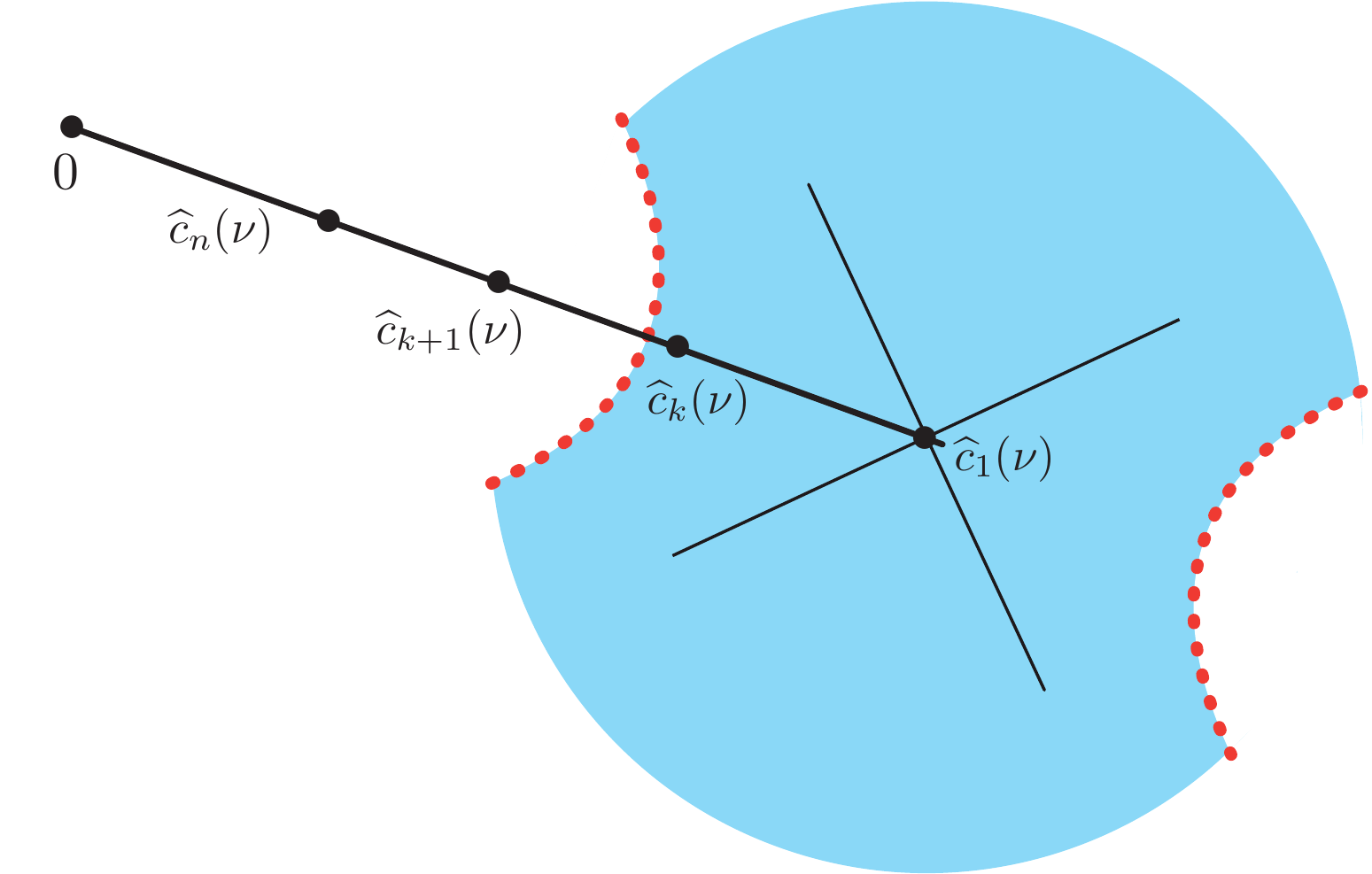}&&\includegraphics[scale=.4]{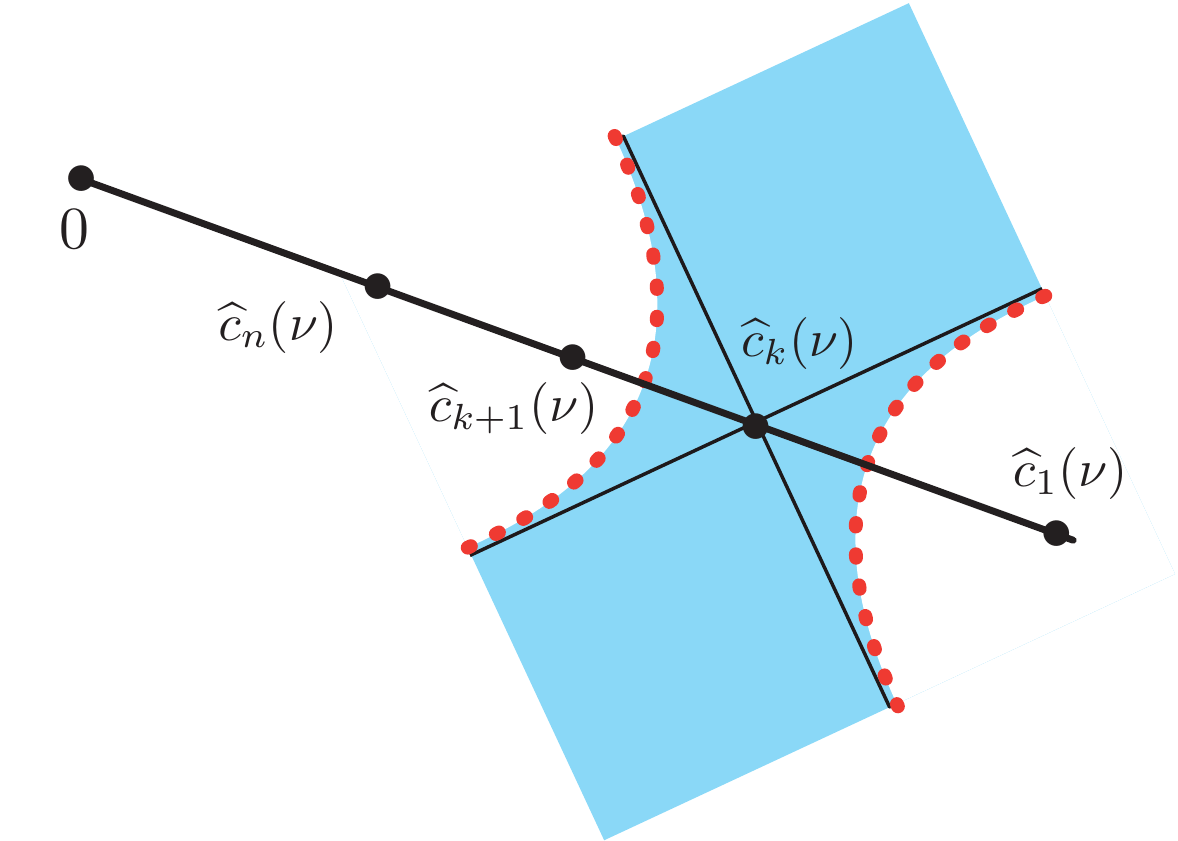}\\
\smaller$\ov\Delta_w^\rrd(\wh\theta_o^{(\nu)},c_1,\gamma)$&&\smaller$\ov\Delta_w^\rrd(\wh\theta_o^{(\nu)},c_k,\gamma)$\\
&&\smaller we could also have $\wh c_1(\nu)$ inside\\
&&\smaller the domain, but not $\wh c_{k+1}(\nu)$\\
\includegraphics[scale=.4]{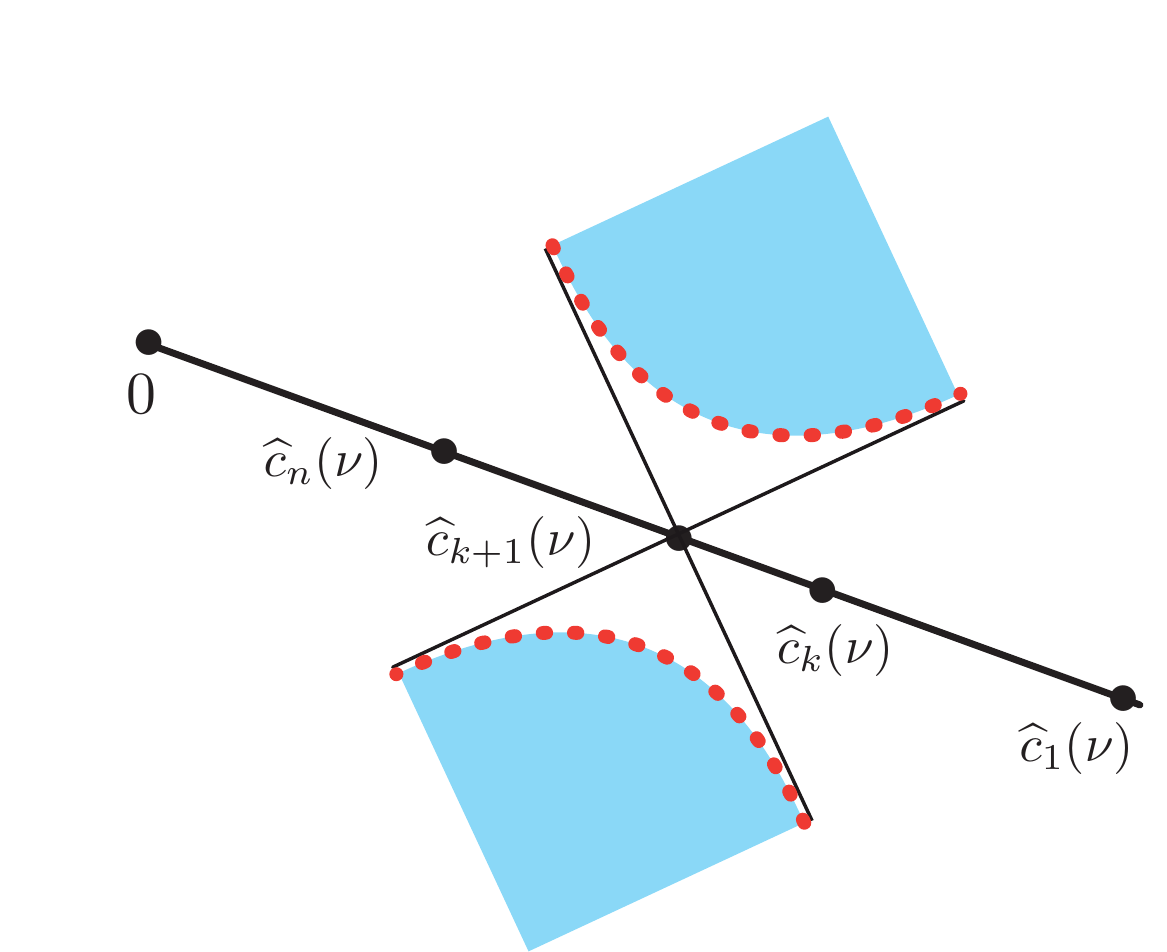}&&
\includegraphics[scale=.4]{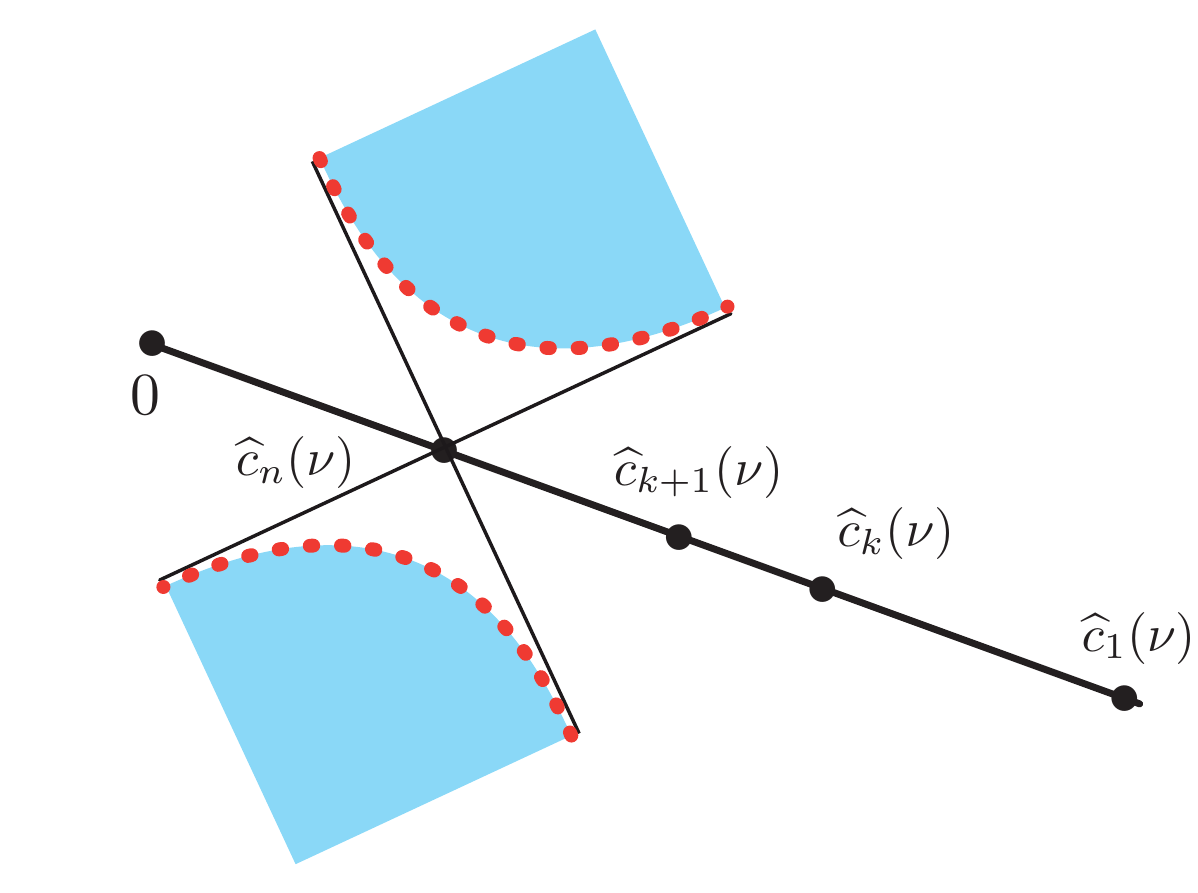}\\
\smaller$\ov\Delta_w^\rrd(\wh\theta_o^{(\nu)},c_{k+1},\gamma)$&&
\smaller$\ov\Delta_w^\rrd(\wh\theta_o^{(\nu)},c_n,\gamma)$
\end{tabular}
\caption{$\nu$ even, \eg $\nu=0$, $\theta_o\in(0,\pi/8)$, $\cos2\wh\theta_o^{(0)}>0$}\label{fig:Awmodnucb}
\end{center}
\end{figure}

Each domain \hbox{$\ov\Delta_w^\rrd(\wh\theta=\wh\theta_o^{(\nu)},c_j,\gamma)$} takes the form of one of the domains in Figure~\ref{fig:rdnueven} if~$\nu$ is even (\resp Figure \ref{fig:rdnuodd} if~$\nu$ is odd), depending whether $j\leq k$ or $j>k$. In Figures~\ref{fig:Awmodnuoddcb} and~\ref{fig:Awmodnucb}, we focus on the real half-line containing the centers $\wh c_j(\nu):=\wh c_j(\wh\theta_o^{(\nu)})$. If $\nu$ is odd, we get the pictures of Figure \ref{fig:Awmodnuoddcb} for the domain $\ov\Delta_w^\rrd(\wh\theta_o^{(\nu)},c_j,\gamma)$, according to Lemma \ref{lem:positionrelative}. If~$\nu$ is even, we get the pictures of Figure~\ref{fig:Awmodnucb}.

\begin{figure}[htb]
\begin{center}
\begin{tabular}{cc}
\includegraphics[scale=.3]{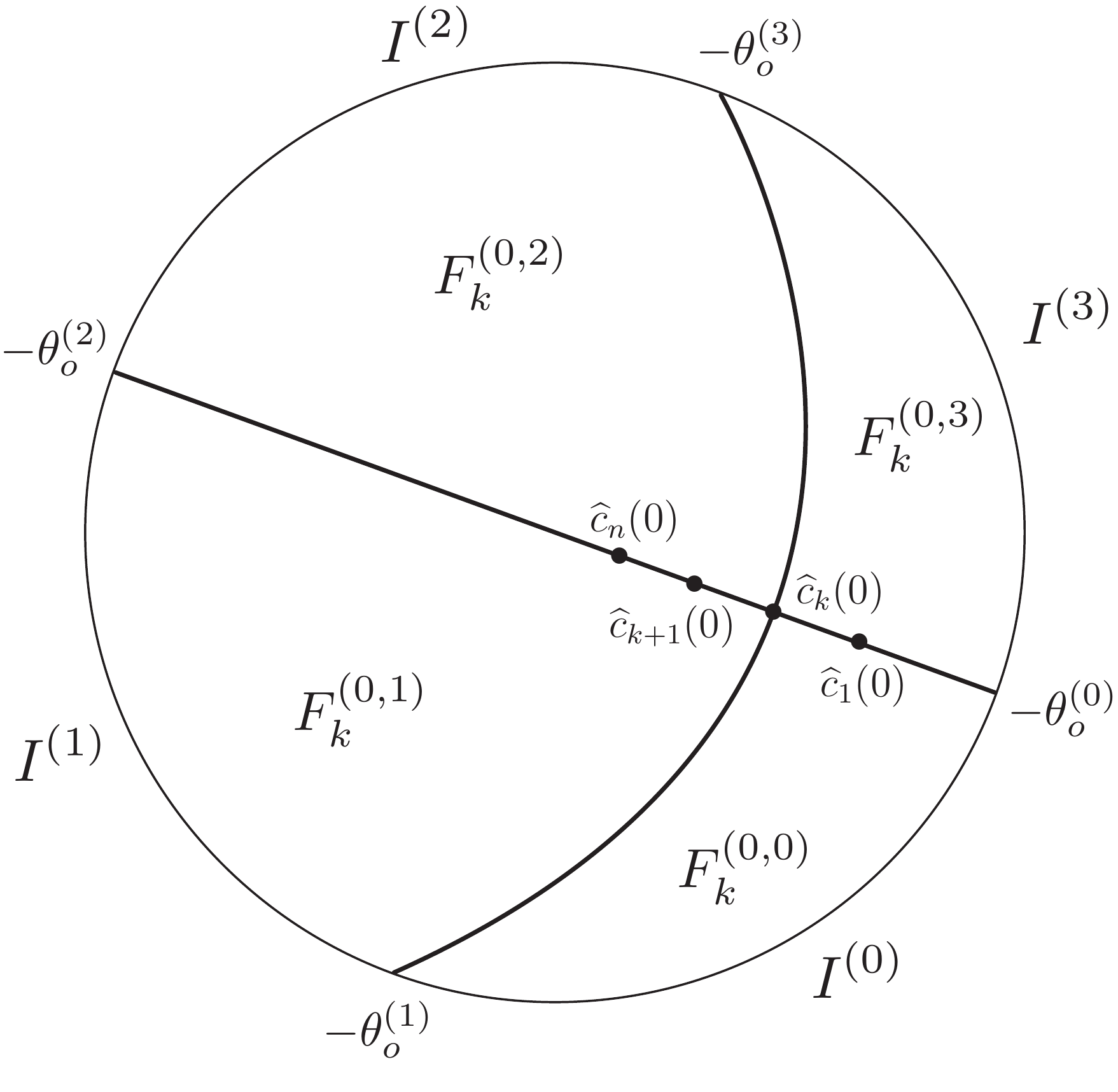}&\includegraphics[scale=.3]{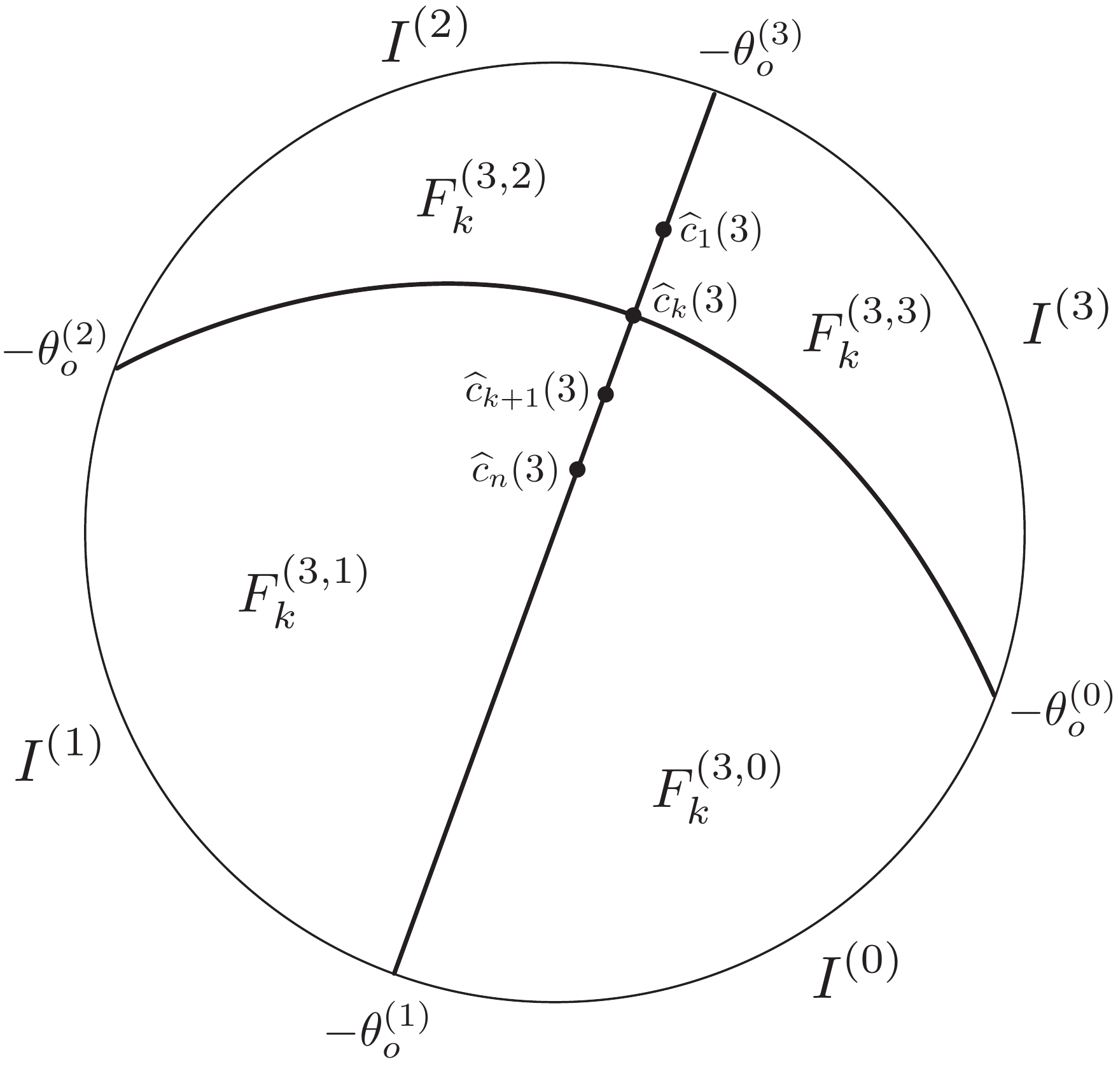}\\
$\nu$ even (\eg $\nu=0$)&$\nu$ odd (\eg $\nu=3$)
\end{tabular}
\caption{The covering $\ccF_k^{(\nu)}$ by closed subsets $F_k^{(\nu,\mu)}$ ($\mu\in\ZZ/4\ZZ$)}\label{fig:ccFknu}
\end{center}
\end{figure}

We notice that the covering $(\ov\Delta_{w,\mu})_{\mu\in\ZZ/4\ZZ}$ is a Leray covering for the sheaves $\cG_{<\gamma}$ when $\nu$ is odd, but not when $\nu$ is even. On the pictures of Figure~\ref{fig:Awmodnuoddcb} and Figure~\ref{fig:Awmodnucb}, it is induced by the four quadrants centered at the origin and one of the corresponding edges is the half-line drawn in these pictures. As soon as some domain $\ov\Delta_{w,\mu}\cap\ov\Delta_w^\rrd(\wh\theta_o^{(\nu)},c_k,\gamma)$ has two red boundary components, this produces a non-zero~$H^1$ for $\cG_{<\gamma}$. This occurs in pictures like $\ov\Delta_w^\rrd(\wh\theta_o^{(\nu)},c_j,\gamma)$ ($j=1,\dots,k$) in Figure \ref{fig:Awmodnucb}. However, for $\gamma$ (and thus $k$) fixed as above, we will consider a slightly different closed covering $\ccF_k^{(\nu)}$, as Figure \ref{fig:ccFknu}. Since $\cG$ is constant in the interior of $\ov\Delta_w$, we can as well recover $\cG_{\wh\theta_o^{(\nu)}}$ and $\cG_{<\gamma,\wh\theta_o^{(\nu)}}$ by formulas like in Lemmas \ref{lem:cG'} and \ref{lem:cGgamma}.

\begin{proof}[Proof of Theorem \ref{th:main}]
Let us start by analyzing the sheaf $\cG_{\wh\theta_o^{(\nu)}}$, following Remark \ref{rem:restrGgammafibre}. The disc $\ov\Delta_w\times\{\wh\theta_o^{(\nu)}\}$ is represented on Figure~\ref{fig:ccFknu}, together with its closed Leray covering $(F^{(\mu)})_{\mu\in\ZZ/4\ZZ}$ (here $k$ and $\nu$ are fixed and we forget them in the notation).
On each $F^{(\mu)}$, $\cG_{\wh\theta_o^{(\nu)}}$ is decomposed as the direct sum $\bigoplus_{c\in C}\cG_{c,\mu,\wh\theta_o^{(\nu)}}$, with $\cG_{c,\mu,\wh\theta_o^{(\nu)}}$ constant on $\Delta_w\cap F^{(\mu)}$ and the non-dashed boundary, and zero on the dashed boundary. The gluing maps are as in Lemma \ref{lem:cG'}. It is then clear that the covering $(F^{(\mu)})$ is Leray for $\cG_{\wh\theta_o^{(\nu)}}$ and that, for each $\mu$, $\Gamma(F^{(\mu)},\cG_{\wh\theta_o^{(\nu)}})=0$, so that the corresponding \v{C}ech complex starts in degree one. Let us denote by $[-\theta_o^{(\nu)}]$ the~half-line containing the centers of the hyperbolas, regarded as the intersection of two closed subsets of the covering $(F^{(\mu)})$.

\begin{lemme}\label{lem:quasiisoCech}
The following morphism of complexes is a quasi-isomorphism:
\[
\xymatrix@R=.5cm{
0\ar[r]&\cC^1(F^{(\cbbullet)},\cG_{\wh\theta_o^{(\nu)}})\ar[r]^{\delta_1}\ar[d]&\cC^2(F^{(\cbbullet)},\cG_{\wh\theta_o^{(\nu)}})\ar[r]\ar[d]&\cdots\\
0\ar[r]&\Gamma([-\theta_o^{(\nu)}],\cG_{\wh\theta_o^{(\nu)}})\ar[r]&0\ar[r]&\cdots
}
\]
\end{lemme}

\begin{proof}
We already know that the upper complex has cohomology in degree one only, as a consequence of Theorem \ref{th:mochi}, but it can be proved by direct arguments similar to those given below. We are thus reduced to proving that the projection $\ker\delta_1\to\Gamma([-\theta_o^{(\nu)}],\cG_{\wh\theta_o^{(\nu)}})$ is an isomorphism.
 
Assume that $\nu$ is odd, \eg $\nu=3$ as in Figure \ref{fig:ccFknu} (the case $\nu$ even is similar). In the \v{C}ech complex, we identify $\cC^1(F^{(\cbbullet)},\cG_{\wh\theta_o^{(\nu)}})$ with $\cL_{\theta_o^{(1)}}\oplus\cL_{\theta_o^{(3)}}\oplus L\oplus L$ (we keep the notation of Diagram \eqref{eq:diagL}) and $\cC^2(F^{(\cbbullet)},\cG_{\wh\theta_o^{(\nu)}})$ with $L\oplus L\oplus L\oplus L$. Let us write an element $\alpha$ of the former as $\alpha(01)\oplus\alpha(23)\oplus\alpha(13)\oplus\alpha(02)$, and an element $\beta$ of the latter as $\beta(012)\oplus\beta(123)\oplus\beta(230)\oplus\beta(301)$. We then have
\begin{align*}
\delta_1(\alpha)(012)&=-\alpha(02)+b^{(1)}{}^{-1}(\alpha(01))\\
\delta_1(\alpha)(123)&=b^{(3)}{}^{-1}(\alpha(23))-\alpha(13)\\
\delta_1(\alpha)(230)&=-\alpha(20)+b^{(3)}{}^{-1}(\alpha(23))\\
\delta_1(\alpha)(301)&=b^{(1)}{}^{-1}(\alpha(01))-\alpha(13).
\end{align*}
Therefore, the map $\alpha\mto\alpha(23)$ induces an isomorphism from $\Ker\delta_1$ to $\cL_{\theta_o^{(3)}}$, as wanted.
\end{proof}

It follows that $\wh\cL_{\wh\theta_o^{(\nu)}}=H^1(\ov\Delta_w,\cG_{\wh\theta_o^{(\nu)}})$ is identified with $\Gamma([-\theta_o^{(\nu)}],\cG_{\wh\theta_o^{(\nu)}})$, which is nothing but $\cL_{\theta_o^{(\nu)}}\simeq L$.

\begin{figure}[htb]
\begin{center}
\begin{tabular}{ccccccc}
\includegraphics[scale=.22]{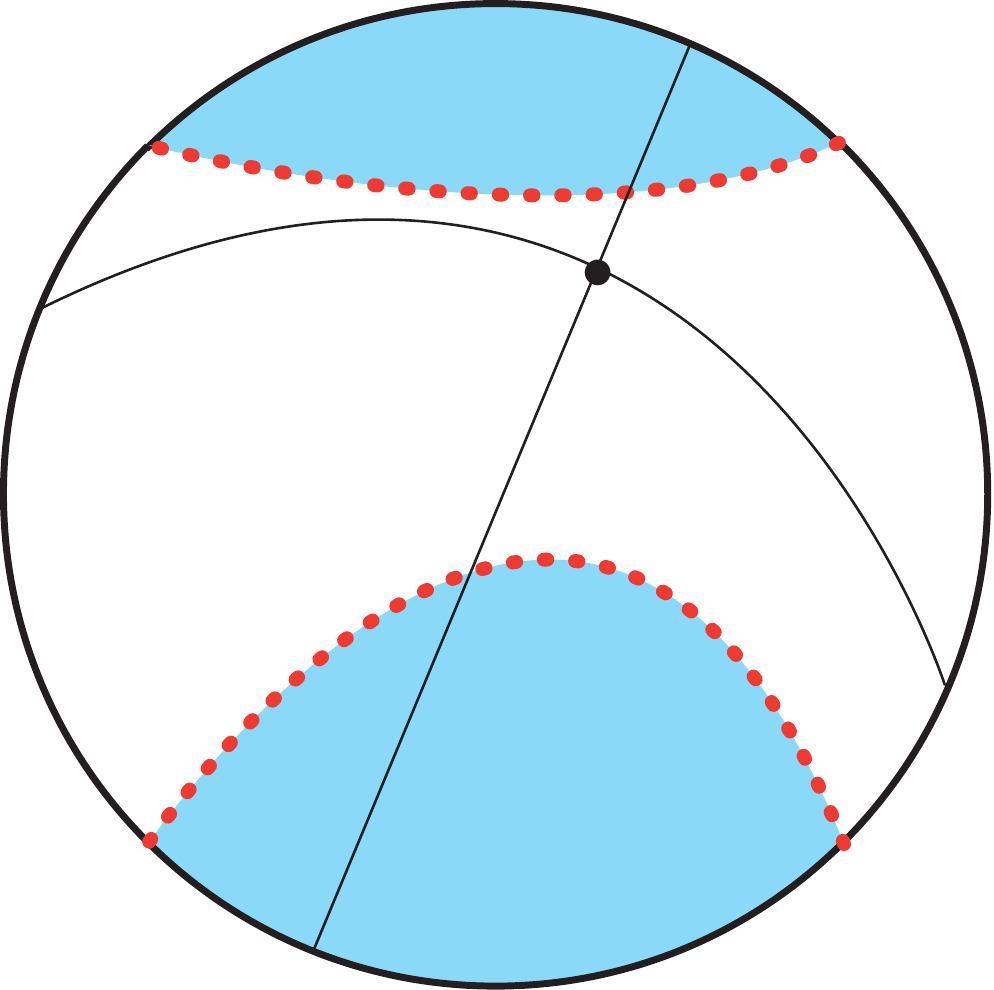}&&
\includegraphics[scale=.22]{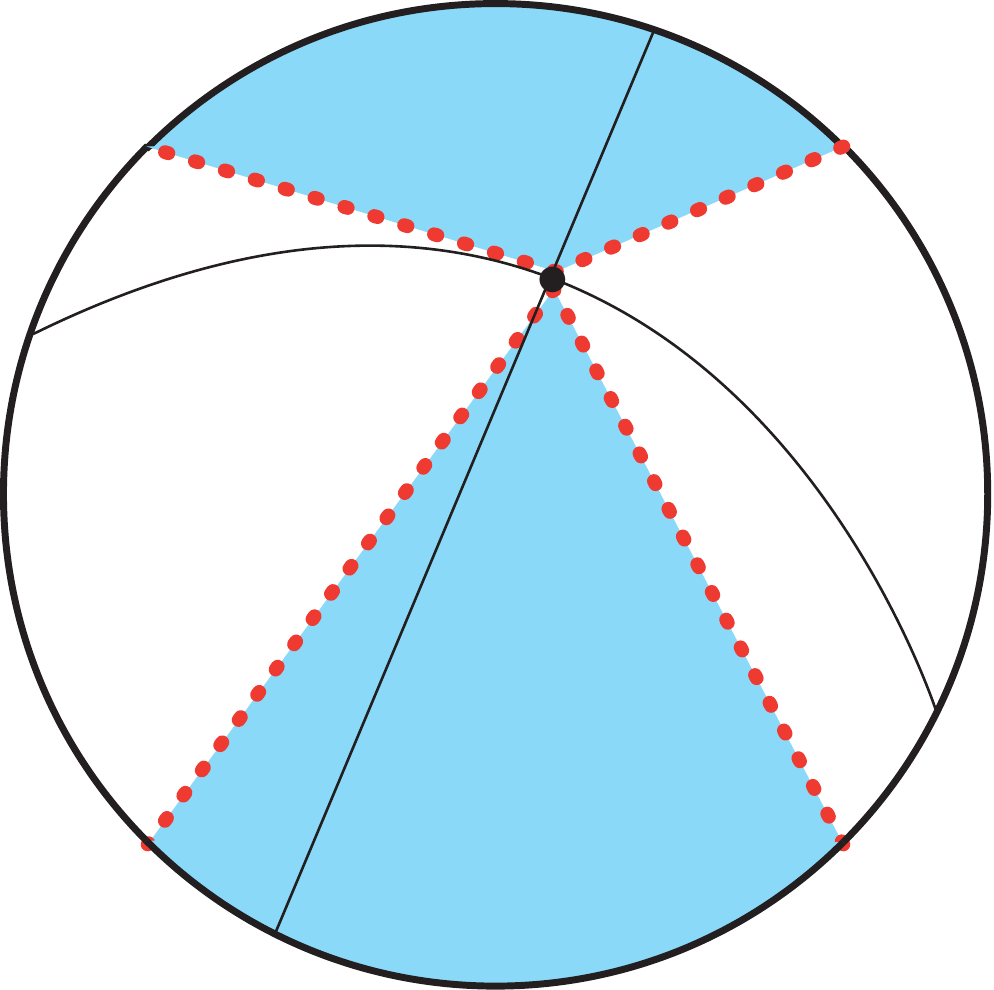}&&\includegraphics[scale=.22]{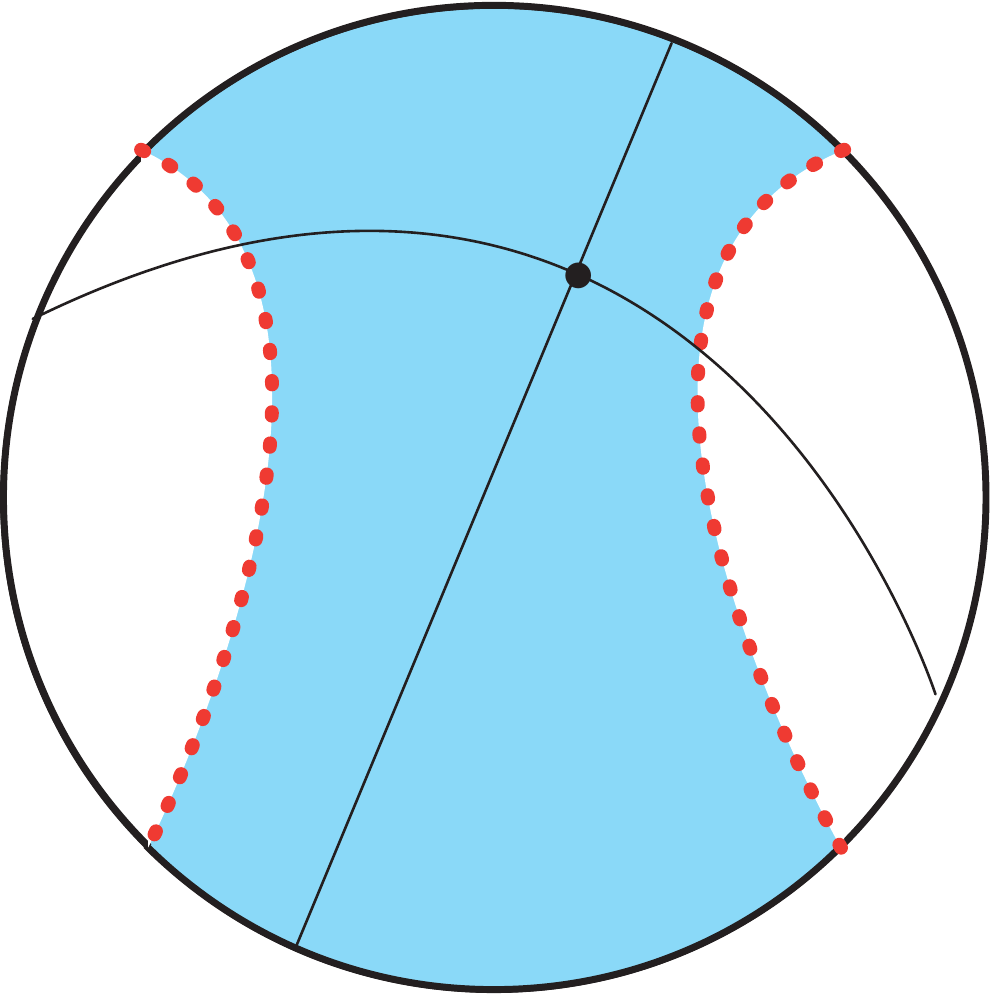}\\
$j<k$&&
$j=k$&&
$j>k$
\end{tabular}
\caption{$\ov\Delta_w^\rrd(\wh\theta,c_j,\gamma)$ if $\nu$ is odd, \eg $\nu=3$}\label{fig:nuoddcj}
\end{center}
\end{figure}

\begin{figure}[htb]
\begin{center}
\begin{tabular}{ccccccc}
\includegraphics[scale=.22]{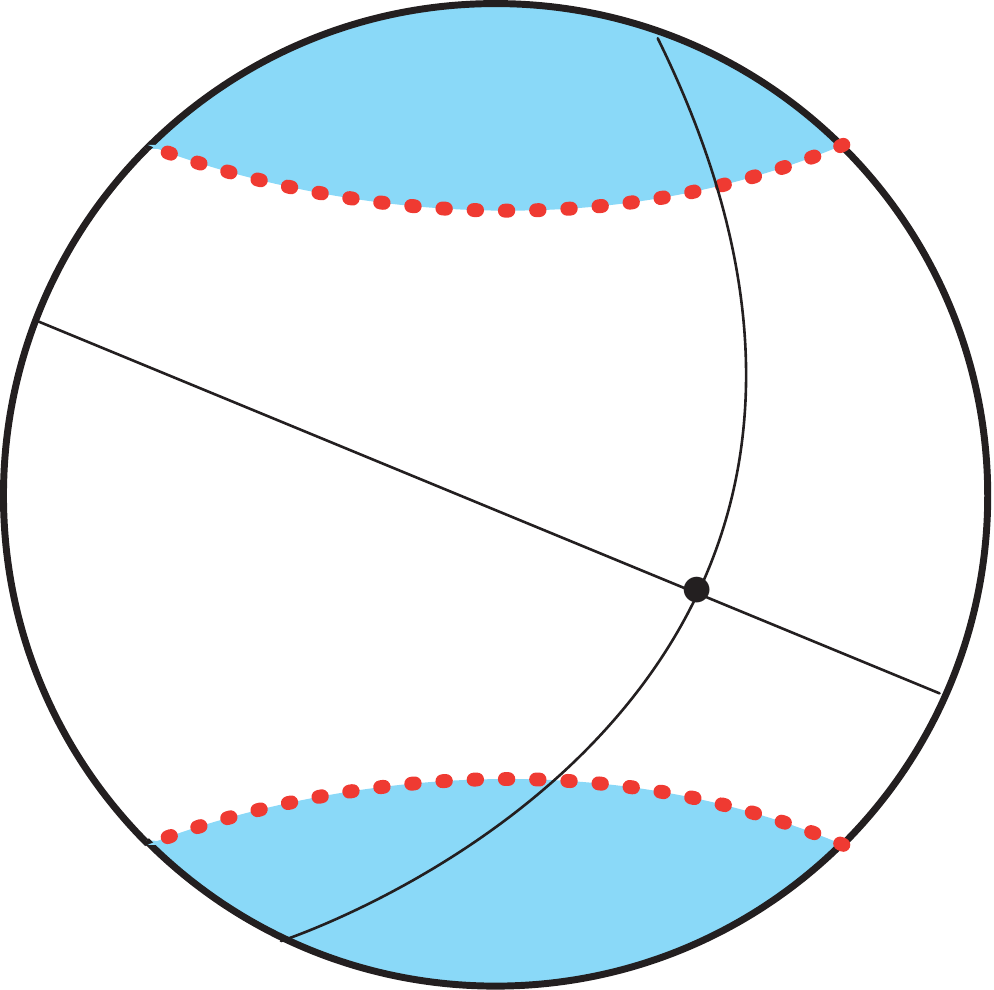}&&
\includegraphics[scale=.22]{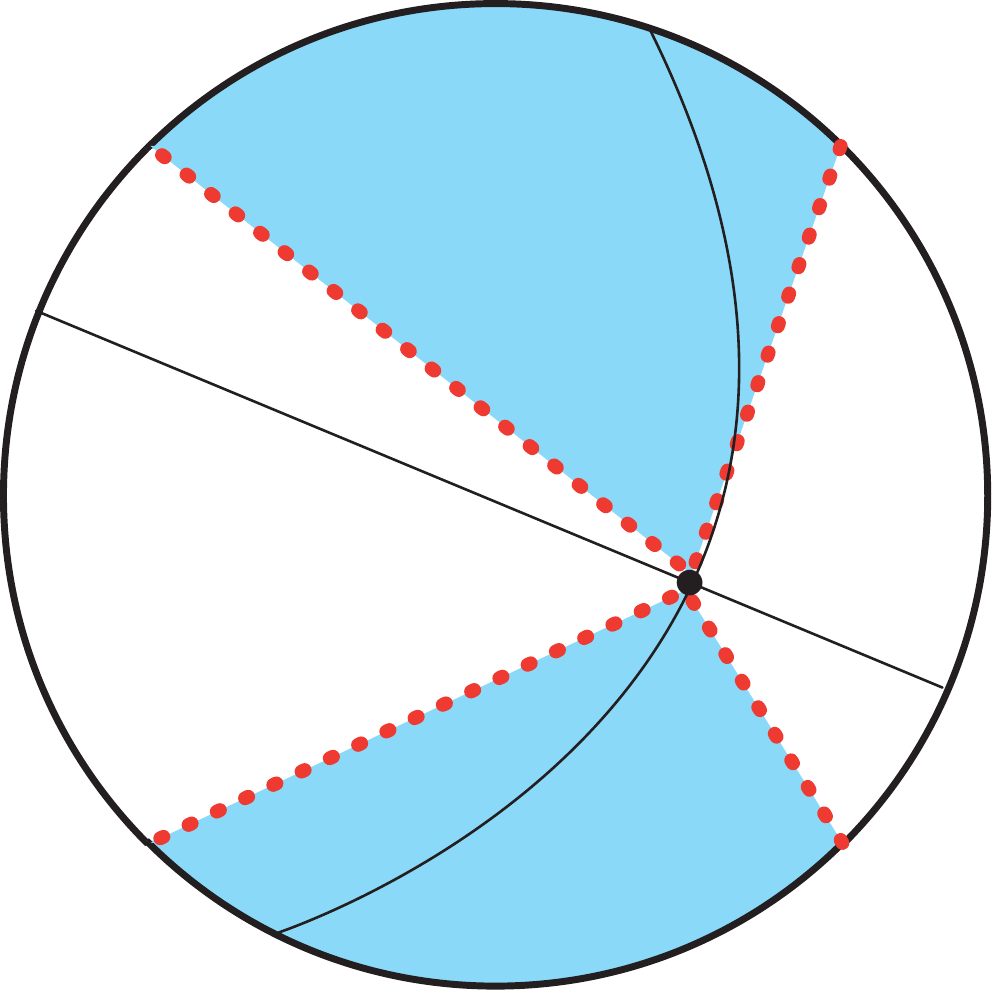}&&
\includegraphics[scale=.22]{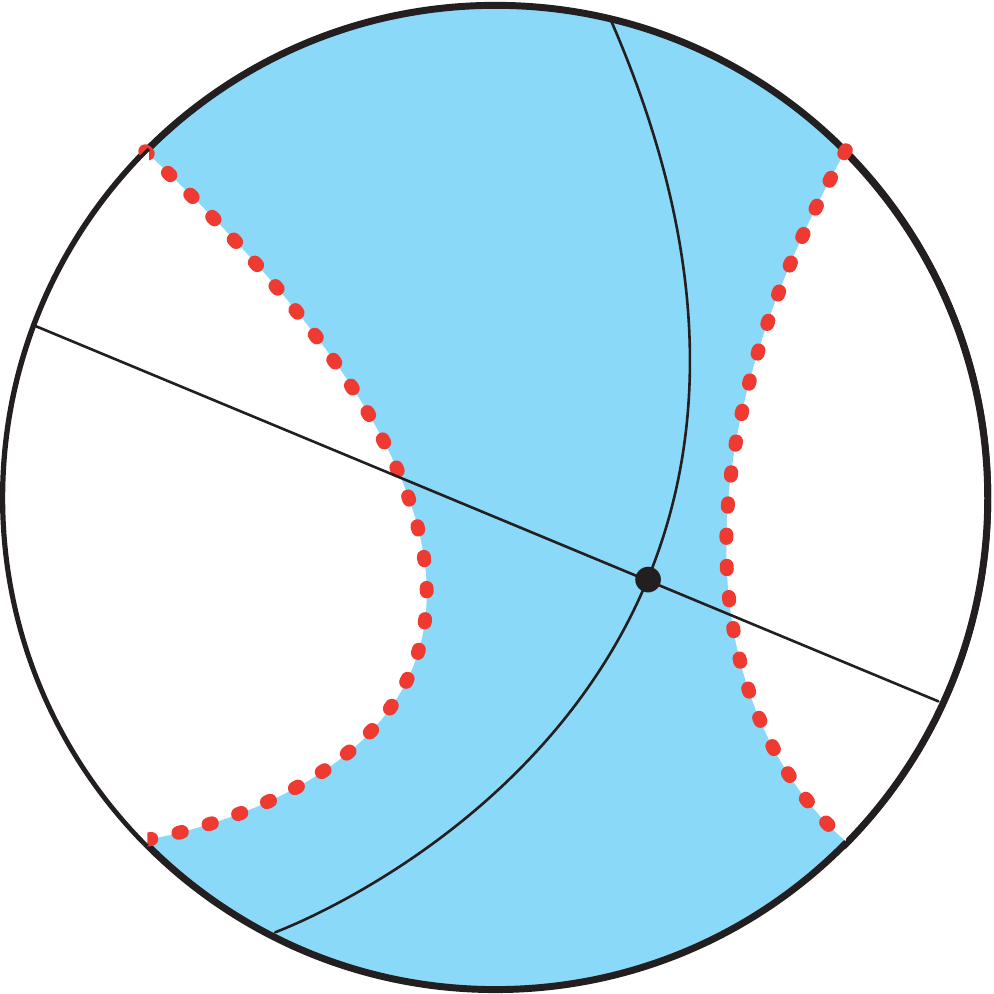}\\
$j>k$&&
$j=k$&&
$j<k$
\end{tabular}
\caption{$\ov\Delta_w^\rrd(\wh\theta,c_j,\gamma)$ if $\nu$ is even, \eg $\nu=0$}\label{fig:nuevencj}
\end{center}
\end{figure}

For each $\gamma\in\CC$, we have $\wh\cL_{<\gamma,\wh\theta_o^{(\nu)}}=H^1(\ov\Delta_w,\cG_{<\gamma,\wh\theta_o^{(\nu)}})$, that we compute with the corresponding \v{C}ech complex. According to the description of Lemma \ref{lem:cGgamma}, which is more convenient to use with a single $G_c$ for each $c$, as in Definition \ref{def:catStokesdatavariant}, and to Figures \ref{fig:Awmodnuoddcb} and \ref{fig:Awmodnucb}, we have $\cC^0(F^{(\cbbullet)},\cG_{<\gamma,\wh\theta_o^{(\nu)}})=0$. We note also that the $c_j$-component ($j=1,\dots,k$ for $\nu$ odd and $j=k+1,\dots,n$ for $\nu$ even) of each $\cC^k(F^{(\cbbullet)},\cG_{<\gamma,\wh\theta_o^{(\nu)}})$ is zero (\cf Figures \ref{fig:nuoddcj} and \ref{fig:nuevencj}), and we have a description of the complex similar to that in the proof of Lemma \ref{lem:quasiisoCech}. The map $H^1(\ov\Delta_w,\cG_{<\gamma,\wh\theta_o^{(\nu)}})\to H^1(\ov\Delta_w,\cG_{\wh\theta_o^{(\nu)}})$ is isomorphic to the map by the projection $\Gamma([-\theta_o^{(\nu)}],\cG_{<\gamma,\wh\theta_o^{(\nu)}})\to\Gamma([-\theta_o^{(\nu)}],\cG_{\wh\theta_o^{(\nu)}})$, which is nothing but the inclusion $L_{\lenu\gamma}\hto L$.
\end{proof}

\begin{remarque}\label{rem:Fpm}
The formulation of Theorem \ref{th:main} makes clear the property $\wh{{\wh M\hphantom{\,}}}\!\simeq\iota^+M$, where $\iota$ denotes here the involution $t\mto-t$, since $\wh{\!\wh C}=C$. Similarly, considering the Laplace transformation $\exp(t\tau)$ (\cf Remark \ref{rem:Lpm}) amounts to changing $-t\tau$ to $+t\tau$ in the formulas. The centers of the hyperbolas are now defined by $w'=-\wh c(\wh\theta)$, and we have $\arg_{w'}(-\wh c(\wh\theta_o^{(\nu)}))=-\theta_o^{(\nu+1)}$. It is then clear that the composition of both Laplace transformations is equal to $\id$, since the rotation by $-\pi/2$ for the first one is annihilated by the rotation by $+\pi/2$ for the second one.
\end{remarque}

\setcounter{section}{1}
\setcounter{subsection}{0}
\setcounter{equation}{0}
\def\thesection{\Alph{section}}
\section*{Appendix. Topological computation of moderate and rapid \texorpdfstring{\hbox{decay}}{decay} de~Rham complexes}\label{sec:prelimcomp}

Let $R$ be a free $\CC[u,u^{-1},v]$-module of finite rank with a flat connection (hence with poles along the divisor $D:=\{u=0\}$) having regular singularities along $u=0$. Recall that $E^{v/u^2}=\big(\CC[u,u^{-1},v]\rd+\rd(v/u^2)\big)$ and $E^{v^2/u^2}=\big(\CC[u,u^{-1},v]\rd+\rd(v^2/u^2)\big)$. We will mainly interested in the behaviour at $u=v=0$, and we will denote by the same letters the corresponding meromorphic germs at the origin, over the ring $\cO_{u,v}[1/u]$. We thus consider a germ at $u=v=0$ of the form $E^{v/u^2}\otimes\nobreak R$ or $E^{v^2/u^2}\otimes R$.

\subsubsection*{Geometry}
Let $X$ denote a neighbourhood of the origin in $\CC^2_{(u,v)}$ and let $\varpi=\varpi_X:\wt X=\wt X(D)\to X$ denote the real blowing up map of $X$ along $D$. The boundary $\partial\wt X(D)$ of $\wt X(D)$ is identified with $D\times S^1_u$, with coordinates $v$ on $D$ and $\theta:=\arg u$ on~$S^1_u$. Let $\cL$ denote the local system determined by $R$ on $D\times S^1_u$.

Let $j_0:D\moins\{0\}\hto D$ denote the inclusion, and let $\wtj_0:(D\moins\{0\})\times S^1_u\hto D\times S^1_u$ denote the corresponding inclusion. Let $i_0:\{0\}\hto D$ \resp $\wti_0:\{0\}\times S^1_u\hto D\times S^1_u$ denote the complementary inclusions.

We denote by $\wt D$ the real blow-up space of $D$ at the origin $v=0$, so that we can identify $\wt D$ with $[0,\epsilon)\times S^1_v$. We can fill the hole by gluing a disc along the boundary~$\partial\wt D$. We thus get a space $\underline{\wt D}$, with a map $\underline\varpi_D:\underline{\wt D}\to D$, which contracts the closure of the filling disc to the origin in $D$. The restriction $\varpi_D$ to $\wt D$ is the real blowing up map, which contracts the boundary $\partial\wt D$ to the origin.

On $\wt D\times S^1_u$, we denote by $\wt L_{1,+}$ (\resp $\wt L_{2,+}$) the open subset defined by $\reel(v/u^2)>0$ (\resp $\reel(v^2/u^2)>0$), that is, $\arg v-2\theta\in(-\pi/2,\pi/2)\bmod2\pi$ (\resp $\arg v-\theta\in(-\pi/4,\pi/4)\bmod\pi$). We also denote by $L_{i,+}$ their restriction to $(D\moins\{0\})\times S^1_u$ and by $\wtj_i:\wt L_{i,+}\hto\wt D\times S^1_u$ $(i=1,2$) the open inclusion. We will also consider the subspaces $\underline{\wt L}_{i,+}$ in $\underline{\wt D}\times S^1_u$ defined as the union of the two sets $\wt L_{i,+}$ and $(\underline{\wt D}\moins\wt D)\times S^1_u$, and the corresponding open inclusion $\underline{\wtj}_i$. We will then denote by $\underline{\wt\beta}_i$ the functor $\underline\wtj_{i,!}\underline\wtj^{-1}_i$, etc.

The sets $\wt L_{i,+},\underline{\wt L}_{i,+}$ are topological fibrations over $S^1_u$. We illustrate below a typical fibre of this fibration (contained in $\wt D$ and $\underline{\wt D}$), and the fibration is obtained by rotating the picture around the center of the (empty or full) disc. The map $\varpi_D$ (\resp $\underline\varpi_D$) contracts the boundary circle (\resp the disc) to the origin.

\begin{figure}[htb]
\begin{center}
\includegraphics[scale=.6]{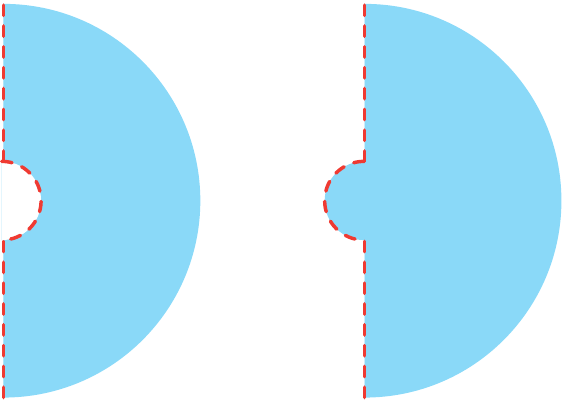}
\caption{Restriction to $\arg u=\theta_o$ of $L_{1,+}$ and $\underline{\wt L}_{1,+}$}\label{fig:L1}
\end{center}
\end{figure}

\begin{figure}[htb]
\begin{center}
\includegraphics[scale=.6]{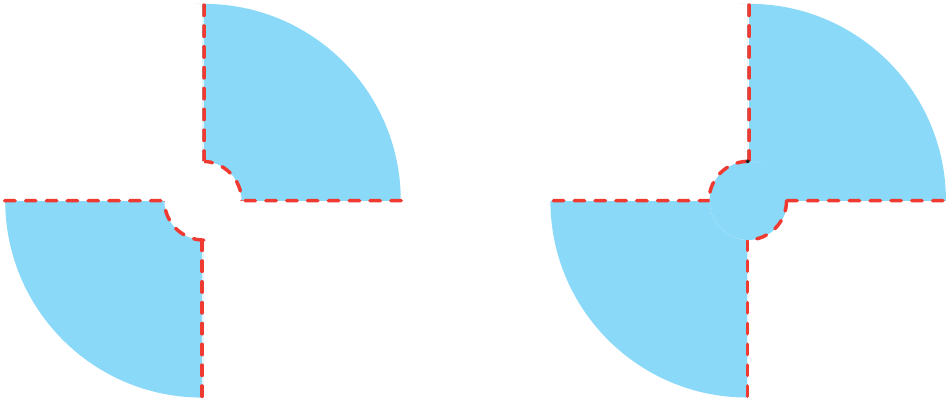}
\caption{Restriction to $\arg u=\theta_o$ of $L_{2,+}$ and $\underline{\wt L}_{2,+}$}\label{fig:L2}
\end{center}
\end{figure}

\subsubsection*{Analysis}
The space $\wt X(D)$ is equipped with the sheaves
\begin{itemize}
\item
$\cA_{\wt X(D)}^\modD$ of holomorphic functions on $X^*:=X\moins D$ having moderate growth along $\partial\wt X(D)$,
\item
$\cA_{\wt X(D)}^\rdD$ of holomorphic functions on $X^*:=X\moins D$ having rapid decay along $\partial\wt X(D)$.
\end{itemize}
We will be mainly interested to their restriction to $\partial\wt X(D)=D\times S^1_u$.

\subsubsection*{Moderate and rapid decay de~Rham complexes}

Given a free $\cO_X(*D)$-module $\cM$ with flat connection, we consider the corresponding de~Rham complexes $\DR^\modD(\cM)$ and $\DR^\rdD(\cM)$ on $\wt X(D)$. These are complexes on $\wt X(D)$ which are both equal to the holomorphic de~Rham complex $\DR(\cM)$ when restricted to $X^*$.

\begin{lemme}\label{lem:DRmodv/u}
We have
\begin{starequation}\label{eq:DRmodv/u1}
\DR^\rdD(E^{v/u^2}\otimes R)=\DR^\modD(E^{v/u^2}\otimes R)=\wtj_{0,!}\beta_1\cL,
\end{starequation}%
and
\begin{starstarequation}\label{eq:DRmodv/u2}
\begin{aligned}
\DR^\rdD(E^{v^2/u^2}\otimes R)&=\bR\varpi_{D,*}\wt\beta_2\cL=\wtj_{0,!}\beta_2\cL,\\
\DR^\modD(E^{v^2/u^2}\otimes R)&=\bR\underline\varpi_{D,*}\underline{\wt\beta}_2\cL.
\end{aligned}
\end{starstarequation}
\end{lemme}

As a consequence, the germ at $\{0\}\times S^1_u$ of the complexes in \eqref{eq:DRmodv/u1} and in the first line of \eqref{eq:DRmodv/u2} are identically zero. On the other hand,
\begin{equation}\label{eq:H1DRmod}
\cH^0\DR^\modD(E^{v^2/u^2}\otimes R)=0\quad\text{and}\quad\cH^1\DR^\modD(E^{v^2/u^2}\otimes R)\simeq\cL_{|\{0\}\times S^1_u}.
\end{equation}

\begin{remarque}
In both examples, we observe that $\DR^\rdD$ commutes with the restriction to $v=0$, while $\DR^\modD$ does not: indeed, the restriction to $v=0$ of $E^{v/u^2}\otimes R$ or of $E^{v^2/u^2}\otimes R$ is equal to the restriction of $R$ to $v=0$, hence is regular; thus its moderate de~Rham complex is a locally constant sheaf, while its rapid decay de~Rham complex is zero.
\end{remarque}

\begin{proof}[Proof of \eqref{eq:DRmodv/u1}]
Let us first prove the result with $E^{v/u}$ instead of $E^{v/u^2}$. It is enough to check the result at the origin of $D$, since it is clear away from the origin. We blow up the origin $e:X'\to X$, so that $X'$ comes equipped with two charts with coordinates $(u',v')$ and $(u'',v'')$ such that $u=u'v'$, $v=v'$ and $u=u''$, $v=u''v''$. The pull-back $e^{-1}(D)$ consists of the strict transform $D'=\{v=0\}$ of $D$ and the exceptional divisor $E=\{v'=0\}\cup\{u''=0\}$. We have a map between the real blow-up spaces $\wt e:\wt X'(D'\cup E)\to\wt X(D)$.

The reason for using such a complex blowing up and the associated real blowing up is that the moderate or rapid decay de~Rham complexes we are interested in can be computed as the push-forward by $\wt e$ of the corresponding complexes on $\wt X'$, and that, on $\wt X'$, these complexes have cohomology in degree zero at most, and their $\cH^0$ is easily computed (\cf\eg\cite[Chap\ptbl8]{Bibi10}).

\subsubsection*{Chart $(u'',v'')$}
Above this chart, we identify $\partial\wt X'$ with $\Afu_{v''}\times S^1_u$ by identifying $\theta''=\arg u''$ with $\theta=\arg u$. The pull-back $e^+(E^{v/u}\otimes R)$ has regular singularities along $E=\{u''=0\}$, so $\DR^\rdE(e^+(E^{v/u}\otimes R))=0$ on $\partial\wt X'$ in this chart, and $\DR^\modE(e^+(E^{v/u}\otimes\nobreak R))$ is the pull-back sheaf $\wt e^{-1}\wti_0^{-1}\cL$ by the map $\wt e:(\theta'',v'')\mto(\theta'',v''\re^{\ri\theta''})=(\theta,v)$.

\subsubsection*{Chart $(u',v')$}
In this chart, we identify $\wt X'$ with $(\RR_+)^2\times(S^1)^2$ with coordinates $(|v'|,|u'|,\arg v',\arg u')$ and $\wt e$ is given by
\[
(|v'|,|u'|,\arg v',\arg u')\mto(|v'|\re^{\ri\arg v'},\arg v'+\arg u').
\]
Since we already computed the de~Rham complexes away from the strict transform of $D$, let us restrict above the strict transform $v'=0$, and thus, above the origin $u'=v'=0$, which is equal to $S^1_{v'}\times S^1_{u'}$. In the neighbourhood of $u'=v'=0$, the pull-back of $E^{v/u}\otimes R$ is equal to $E^{1/u'}\otimes e^+R$.

Let us identify $S^1_{v'}\times S^1_{u'}$ with $S^1_v\times S^1_u$ by the isomorphism $(\arg v',\arg u')\mto(\arg v',\arg v'+\arg u')$. Then, the restriction of $\DR^{\rrd(E\cup D')}(e^+(E^{v/u}\otimes R))$ to this set is zero since $\re^{-1/u'}$ does not have rapid decay at points above $E\moins(D'\cap E)$. That of $\DR^{\rmod(E\cup D')}(e^+(E^{v/u}\otimes R))$ is clearly identified with the restriction of $\wt\beta_1\cL$ to $S^1_v\times S^1_u$.

At this point, we can conclude that
\[
\DR^\rdD(E^{v/u}\otimes R)=\wtj_{0,!}\wtj_0^{-1}\DR^\rdD(E^{v/u}\otimes R)=\wtj_{0,!}\beta_1\cL.
\]

\subsubsection*{Gluing the two charts}
Let us identify topologically $\Afu_{v''}$ with an open disc $B_{v''}$ with coordinate $v''$ (of radius $1$, say). We have a homeomorphism $B_{v''}\times S^1_u\isom B_{w}\times\nobreak S^1_u$ by sending $(v'',\arg u)$ to $(w\!=\!v''\re^{\ri\arg u},\arg u)$. We regard $B_w$ as the filling disc in~$\underline{\wt D}$. Indeed, on $\partial\ov B_w$, we have $\arg w=\arg v$. We also identify \hbox{$(\RR_+)_{v'}\times S^1_{v'}\times S^1_{u'}$} with $(\RR_+)_v\times S^1_v\times S^1_u$ via $\wt e$, sending $(|v'|,\arg v',\arg u')$ to $(|v'|,\arg v',\arg u'+\arg v')$ as above.

Then $\DR^{\rmod(E\cup D')}(e^+(E^{v/u}\otimes R))$ is identified with $\underline{\wt\beta}_1\cL$, and $\DR^\modD(E^{v/u}\otimes R)$ with $\bR\underline\varpi_{D,*}\underline{\wt\beta}_1\cL$.

It remains to check that the latter complex is zero when restricted to any point $(0,\theta)$ of $D\times S^1_u$. At such a point, the germ of $\DR^\modD(E^{v/u}\otimes R)$ has cohomology equal to the cohomology with compact support of the union of an open disc and an open interval in its boundary, which is easily seen to be equal to zero.

By choosing a square root of the monodromy of $\cL$, one expresses $E^{v/u^2}\otimes R$ as the pull-back by the ramification $u\mto u^2$ of a meromorphic connection $E^{v/u}\otimes R'$. Similarly, $\DR^\rdD(E^{v/u^2}\otimes R)$ and $\DR^\modD(E^{v/u^2}\otimes R)$ are the corresponding pull-back complexes. Then \eqref{eq:DRmodv/u1} follows.

For \eqref{eq:DRmodv/u2} we do not need to use a covering with respect to $u$ and we can argue with $v^2/u^2$ as we did with $v/u$. The proof is similar, except the conclusion on the vanishing of the germ of $\DR^\modD(E^{v/u}\otimes R)$ at $(0,\theta)\in D\times S^1$, since the cohomology is now equal to the cohomology with compact support of the union of an open disc and two disjoint open intervals in its boundary. Such a cohomology vanishes in degree $\neq1$, and has rank one in degree one.
\end{proof}

Let now $\cM$ be a locally free $\cO_X(*D)$-module with flat connection, which satisfies, locally on $\wt X(D)$,
\[
\cA_{\wt X}^\modD\otimes_{\varpi^{-1}\cO_X}\cM\simeq\bigoplus_{\lambda\in\Lambda}\cA_{\wt X}^\modD\otimes_{\varpi^{-1}\cO_X}(E^{\lambda v/u^2}\otimes R_ \lambda),
\]
where $\Lambda$ is a finite subset of $\CC^*$.

\begin{corollaire}\label{cor:comprdmod}
With these assumptions, the natural morphism
\[
\DR^\rdD(\cM)\to\DR^\modD(\cM)
\]
is a quasi-isomorphism.
\end{corollaire}

\begin{proof}
We can argue locally on $\partial\wt X(D)$ and, according to the assumption on $\cM$, we can replace $\cM$ with $\bigoplus_{\lambda\in\Lambda}(E^{\lambda v/u^2}\otimes R_ \lambda)$, so that we can apply \eqref{eq:DRmodv/u1}.
\end{proof}

\backmatter
\newcommand{\eprint}[1]{\href{http://arxiv.org/abs/#1}{\texttt{arXiv\string:\allowbreak#1}}}
\providecommand{\bysame}{\leavevmode ---\ }
\providecommand{\og}{``}
\providecommand{\fg}{''}
\providecommand{\smfandname}{\&}
\providecommand{\smfedsname}{\'eds.}
\providecommand{\smfedname}{\'ed.}
\providecommand{\smfmastersthesisname}{M\'emoire}
\providecommand{\smfphdthesisname}{Th\`ese}


\begin{thebibliography}{Moc14}

\bibitem[Ari10]{Arinkin08}
{\scshape D.~Arinkin} -- {\og {Rigid irregular connections on
  $\mathbb{P}^1$}\fg}, \emph{Compositio Math.} \textbf{146} (2010), no.~5,
  p.~1323--1338, \eprint{0808.0742}.

\bibitem[BV89]{B-V89}
{\scshape {\relax D.G}.~Babbitt {\normalfont \smfandname} {\relax
  V.S}.~Varadarajan} -- \emph{Local moduli for meromorphic differential
  equations}, Ast{\'e}risque, vol. 169-170, Soci{\'e}t{\'e} Math{\'e}matique de
  France, Paris, 1989.

\bibitem[DK13]{D-K13}
{\scshape A.~D'Agnolo {\normalfont \smfandname} M.~Kashiwara} -- {\og
  {Riemann-Hilbert correspondence for holonomic D-modules}\fg},
  \eprint{1311.2374}, 2013.

\bibitem[Del07]{Deligne78}
{\scshape P.~Deligne} -- {\og {Lettre \`a B.~Malgrange du 19/4/1978}\fg}, in
  \emph{{Singularit\'es irr\'eguli\`eres, Correspondance et documents}},
  Documents math\'ematiques, vol.~5, Soci{\'e}t{\'e} Math{\'e}matique de
  France, Paris, 2007, p.~25--26.

\bibitem[God64]{Godement64}
{\scshape R.~Godement} -- \emph{Topologie alg\'ebrique et th\'eorie des
  faisceaux}, Hermann, Paris, 1964.

\bibitem[HS11]{H-S09}
{\scshape C.~Hertling {\normalfont \smfandname} C.~Sabbah} -- {\og {Examples of
  non-commutative Hodge structures}\fg}, \emph{Journal de l'Institut
  mathématique de Jussieu} \textbf{10} (2011), no.~3, p.~635--674,
  \eprint{0912.2754}.

\bibitem[HS14]{H-S14}
{\scshape M.~Hien {\normalfont \smfandname} C.~Sabbah} -- {\og {The local
  Laplace transform of an elementary irregular meromorphic connection}\fg},
  \eprint{1405.5310}, 2014.

\bibitem[Maj84]{Majima84}
{\scshape H.~Majima} -- \emph{Asymptotic analysis for integrable connections
  with irregular singular points}, Lect. Notes in Math., vol. 1075,
  Springer-Verlag, 1984.

\bibitem[Mal83]{Malgrange83bb}
{\scshape B.~Malgrange} -- {\og {La classification des connexions
  irr\'eguli\`eres {\`a} une variable}\fg}, in \emph{{S\'eminaire E.N.S.
  Math\'ematique et Physique}} (L.~Boutet~{de Monvel}, A.~Douady {\normalfont
  \smfandname} J.-L. Verdier, \smfedsname), Progress in Math., vol.~37,
  Birkh{\"a}user, Basel, Boston, 1983, p.~381--399.

\bibitem[Mal91]{Malgrange91}
\bysame , \emph{{\'E}quations diff\'erentielles {\`a} coefficients
  polynomiaux}, Progress in Math., vol.~96, Birkh{\"a}user, Basel, Boston,
  1991.

\bibitem[Moc10]{Mochizuki09b}
{\scshape T.~Mochizuki} -- {\og {Note on the Stokes structure of the Fourier
  transform}\fg}, \emph{Acta Math. Vietnam.} \textbf{35} (2010), no.~1,
  p.~101--158.

\bibitem[Moc14]{Mochizuki10}
\bysame , \emph{{Holonomic $\mathcal D$-modules with Betti structure}},
  M{\'e}m. Soc. Math. France (N.S.), vol. 138--139, Soci{\'e}t{\'e}
  Math{\'e}matique de France, 2014, \eprint{1001.2336}.

\bibitem[Sab93]{Bibi93}
{\scshape C.~Sabbah} -- {\og {{\'E}quations diff\'erentielles \`a points
  singuliers irr\'eguliers en dimension~$2$}\fg}, \emph{Ann. Inst. Fourier
  (Grenoble)} \textbf{43} (1993), p.~1619--1688.

\bibitem[Sab00]{Bibi97}
\bysame , \emph{{{\'E}quations diff\'erentielles {\`a} points singuliers
  irr\'eguliers et ph\'enom\`ene de Stokes en dimension {$2$}}},
  Ast{\'e}risque, vol. 263, Soci{\'e}t{\'e} Math{\'e}matique de France, Paris,
  2000.

\bibitem[Sab06]{Bibi05b}
\bysame , {\og {Monodromy at infinity and Fourier~transform~II}\fg},
  \emph{Publ. RIMS, Kyoto Univ.} \textbf{42} (2006), p.~803--835.

\bibitem[Sab07]{Bibi00b}
\bysame , \emph{{Isomonodromic deformations and Frobenius manifolds}},
  Universitext, Springer \& EDP~Sciences, 2007, (in French: 2002).

\bibitem[Sab08]{Bibi07a}
\bysame , {\og {An explicit stationary phase formula for the~local formal
  Fourier-Laplace~transform}\fg}, in \emph{{Singularities, vol.~1}}, Contemp.
  Math., American Mathematical Society, Providence, RI, 2008, p.~300--330,
  \eprint{0706.3570}.

\bibitem[Sab13]{Bibi10}
\bysame , \emph{{Introduction to Stokes structures}}, Lect. Notes in Math.,
  vol. 2060, Springer-Verlag, 2013,
  \url{http://dx.doi.org/10.1007/978-3-642-31695-1} \& \eprint{0912.2762}.

\end{thebibliography}
\end{document}